\documentclass[12pt,a4paper,twoside]{article}
\usepackage[french] {babel}
\usepackage{amssymb,amsmath,amsfonts,graphics}
\usepackage{stmaryrd}
\usepackage[all]{xy}

\def\rondI{\newbox\boxx{\hbox{$1$}}\hskip-11pt\raisebox{0pt}{$\displaystyle{\textrm{\Large$\bigcirc$}}$}}
\def\rondII{\newbox\boxx{\hbox{$2$}}\hskip-11pt\raisebox{0pt}{$\displaystyle{\textrm{\Large$\bigcirc$}}$}}
\def\rondIII{\newbox\boxx{\hbox{$3$}}\hskip-11pt\raisebox{0pt}{$\displaystyle{\textrm{\Large$\bigcirc$}}$}}
\def\rondIV{\newbox\boxx{\hbox{$4$}}\hskip-11pt\raisebox{0pt}{$\displaystyle{\textrm{\Large$\bigcirc$}}$}}
\def\rondV{\newbox\boxx{\hbox{$5$}}\hskip-11pt\raisebox{0pt}{$\displaystyle{\textrm{\Large$\bigcirc$}}$}}
\def\rondVI{\newbox\boxx{\hbox{$6$}}\hskip-11pt\raisebox{0pt}{$\displaystyle{\textrm{\Large$\bigcirc$}}$}}

\begin{document}

\newcommand{\hooklongrightarrow}{\lhook\joinrel\longrightarrow}

\begin{center}
\textbf{Images directes I: Espaces rigides analytiques et images directes}

\vskip20mm

Jean-Yves ETESSE
  \footnote{(CNRS - IRMAR, Universit\'e de Rennes 1, Campus de Beaulieu - 35042 RENNES Cedex France)\\
E-mail : Jean-Yves.Etesse@univ-rennes1.fr}
\end{center}
 
 \vskip60mm
 \noindent\textbf{Sommaire}\\

		\begin{enumerate}
		\item[0.] Introduction
		\item[1.] Changement de base pour un morphisme propre
		\item[2.] Sorites sur les voisinages stricts
		\item[3.] Images directes d'isocristaux
			\begin{enumerate}
			\item[3.1.] Sections surconvergentes
			\item[3.2.] D\'efinition des images directes
			\item[3.3.] Changement de base
			\item[3.4.] Surconvergence des images directes
			\end{enumerate}
		\end{enumerate}
\newpage

\noindent\textbf{R\'esum\'e}\\

Cet article est le premier d'une s\'erie de trois articles consacr\'es aux images directes d'isocristaux: ici nous consid\'erons des isocristaux sans structure de Frobenius; dans le deuxi\`eme [Et 6] (resp. le troisi\`eme [Et 7]), nous introduirons une structure de Frobenius dans le contexte convergent (resp. surconvergent).\\

Pour un morphisme propre et lisse relevable nous \'etablissons la surconvergence des images directes, gr\^ace \`a un th\'eor\`eme de changement de base pour un morphisme propre entre espaces rigides analytiques. Ce r\'esultat r\'epond partiellement \`a une conjecture de Berthelot sur la surconvergence des images directes par un morphisme propre et lisse.\\
  
\noindent\textbf{Abstract}\\

This article is the first one of a series of three articles devoted to direct images of isocrystals: here we consider isocrystals without Frobenius structure; in the second one [Et 6] (resp. the third one [Et 7]), we will introduce a Frobenius structure in the convergent (resp. overconvergent) context.\\

For a liftable proper smooth morphism we establish the overconvergence of direct images, owing to a base change theorem for a proper morphism between rigid analytic spaces. This result partially answers  a conjecture of Berthelot on the overconvergence of direct images under a proper smooth morphism.\\


Mots cl\'es: alg\`ebres de Monsky-Washnitzer, sch\'emas formels, espaces rigides analytiques, th\'eor\`eme de changement de base propre, cohomologie rigide, isocristaux surconvergents.\\


\newpage 

\section*{0. Introduction}
 Soit $\mathcal{V}$ un anneau de valuation discr\`ete complet, de corps r\'esiduel $k = \mathcal{V}/\mathfrak{m}$ de caract\'eristique $p > 0$ et de corps des fractions $K$ de caract\'eristique $0$.\\
 
 Cet article est le premier d'une s\'erie de trois consacr\'es aux images directes d'isocristaux.
Ici on rappelle la d\'efinition des images directes d'isocristaux surconvergents par un morphisme $f: X\rightarrow S$ (cf [B 5], ou [C-T], [LS]) de $k$-sch\'emas et on prouve leur surconvergence dans le cas o\`u $f$ est propre et lisse relevable, ou $X$ est une intersection compl\`ete relative dans des espaces projectifs sur $S$.\\

L'outil essentiel dans la preuve de ce r\'esultat est le th\'eor\`eme de changement de base pour un morphisme propre du \S1: on \'etablit celui-ci d'abord dans le cadre des sch\'emas formels [th\'eo (1.1)], puis dans le cadre des espaces rigides analytiques [th\'eo (1.2)] gr\^ace aux travaux de Bosch et L\" utkebohmert [Bo-L\" u 1],[Bo-L\" u 2],[L\" u].\\

La notion de voisinage strict \'etant \'etroitement li\'ee \`a celle de surconvergence des isocristaux, on d\'eveloppe au \S2 quelques propri\'et\'es de ces voisinages stricts dans le cas plongeable: dans le cas cart\'esien, l'image inverse d'un syst\`eme fondamental de voisinages stricts est un syst\`eme fondamental de voisinages stricts [prop (2.2.3)]; le m\^eme r\'esultat vaut pour un morphisme fini et plat, ou fini \'etale, ou fini \'etale galoisien [prop (2.3.1)], ce qui nous servira dans un article ult\'erieur.\\

Apr\`es avoir rappel\'e au \S3 la d\'efinition des images directes d'isocristaux surconvergents donn\'ee par Berthelot dans une note non publi\'ee [B 5] (voir les articles [C-T], [LS] de Chiarellotto-Tsuzuki et Le Stum pour la publication des d\'etails), on \'etablit leur surconvergence pour un morphisme propre et lisse relevable, en m\^eme temps que le th\'eor\`eme de changement de base [th\'eo (3.4.4)]: un ingr\'edient essentiel est l'extension aux voisinages stricts du th\'eor\`eme de changement de base pour un morphisme propre [th\'eo (3.3.2)].\\

Des r\'esultats de rel\`evement de [Et 5] on d\'eduit alors la surconvergence des images directes d'un isocristal surconvergent par un morphisme $f: X\rightarrow S$ projectif et lisse, o\`u $S$ est lisse sur le corps de base $k$ et $X$ relevable en un $\mathcal{V}$-sch\'ema plat [th\'eo (3.4.8.2)] (resp $X$ est une intersection relative dans des espaces projectifs sur $S$ [cor (3.4.8.6)]). Une variante, \'etudi\'ee dans le th\'eor\`eme (3.4.9), ram\`ene la preuve de la surconvergence du cas projectif lisse au cas o\`u la base $S$ est affine et lisse sur $k$. Ici nous avons laiss\'e de c\^ot\'e le cas d'un morphisme  fini \'etale qui sera trait\'e dans [Et 6] et [Et 7] (cf [Et 4, chap III, et IV]).\\

Ces th\'eor\`emes (3.4.8.2), (3.4.8.6) et (3.4.9) [resp. le th\'eo (3.4.4)] r\'esolvent partiellement une conjecture de Berthelot [B 2] sur la surconvergence des images directes par un morphisme propre et lisse: dans [Et 7] ( cf [Et 4, chap IV]) nous en donnons une version avec structure de Frobenius.\\

Dans le cas o\`u la base $S$ est une courbe affine et lisse sur un corps alg\'ebriquement clos $k$ (resp. $S$ est une courbe lisse sur un corps parfait $k$) la conjecture a \'et\'e prouv\'ee pour le faisceau structural par Trihan [Tri] (resp. Matsuda et Trihan [M-T]). Dans le cas relevable, la conjecture a aussi \'et\'e prouv\'ee ind\'ependamment par Tsuzuki par des voies diff\'erentes [Tsu]: en particulier il n'a pas \`a sa disposition le th\'eor\`eme de changement de base pour un morphisme propre entre espaces rigides analytiques, alors que pour nous celui-ci tient une place centrale. Une autre approche est celle de Shiho qui utilise les log-sch\'emas [Shi 1], [Shi 2], [Shi 3].\\

\textbf {Notations:}
 Pour les notions sur les espaces rigides analytiques et la cohomologie rigide nous renvoyons le lecteur \`a [B 3], [B 4], [B-G-R], [C-T] et [LS].
Sauf mention contraire, dans tout cet article on d\'esigne par $\mathcal{V}$ un anneau de valuation discr\`ete complet, de corps r\'esiduel $k = \mathcal{V}/\mathfrak{m}$ de caract\'eristique $p > 0$, de corps des fractions $K$ de caract\'eristique $0$, d'uniformisante $\pi$ et d'indice de ramification $e$.\\

On suppose donn\'e un entier $a \in \mathbb{N}^{\ast}$ et on d\'esigne par $C(k)$ un anneau de Cohen de $k$ de caract\'eristique $0$ [Bour, AC IX, $\S$ 2, $\no3$, prop 5] : $C(k)$ est un anneau de valuation discr\`ete complet d'id\'eal maximal $p\ C(k)$ [EGA $O_{IV}$, 19.8.5] et on note $K_{0}$ son corps des fractions, $K_{0} =$ Frac $(C(k))$. Il existe une injection fid\`element plate $C(k) \hookrightarrow \mathcal{V}$ qui fait de $\mathcal{V}$ un $C(k)$-module libre de rang $e$ [EGA $O_{IV}$, 19.8.6, 19.8.8] et [Bour, AC IX, $\S$ 2, $\no1$, prop 2]. On fixe un rel\`evement $\sigma : C(k) \rightarrow C(k)$ de la puissance $a^{i \grave{e}me}$ du Frobenius absolu de $k$ comme dans  [Et 2, I, 1.1]; on suppose que l'on peut \'etendre $\sigma$ en un endomorphisme de $\mathcal{V}$, encore not\'e $\sigma$, de telle sorte que $\sigma(\pi) = \pi$; on notera encore $\sigma$ l'extension naturelle de $\sigma$ \`a $K$: lorsque $k$ est parfait $C(k)$ est isomorphe \`a l'anneau $W(k)$ des vecteurs de Witt de $k$ et $\sigma$ est un automorphisme de $K$. Si $k \hookrightarrow k'$ est une extension de corps de caract\'eristique $p > 0$, $\mathcal{V'} : = \mathcal{V} \otimes_{C(k)} C(k'), K' = \operatorname{Frac} (\mathcal{V'})$, on peut relever la puissance $a^{i \grave{e}me}$ du Frobenius absolu de $k'$ en un morphisme $\sigma' : K' \rightarrow K'$ au-dessus de $\sigma : K \rightarrow K$ [Et 2, I, 1.1].

\newpage

\section*{1. Changement de base pour un morphisme propre}

 Le th\'eor\`eme suivant est un pr\'ealable pour les th\'eor\`emes de changement de base en g\'eom\'etrie rigide.
 \vskip 3mm
\noindent \textbf{Th\'eor\`eme (1.1)}.
\textit{Soient $\mathcal{V}$ un anneau noeth\'erien et $I \varsubsetneq \mathcal{V}$ un id\'eal ; on suppose $\mathcal{V}$ s\'epar\'e et complet pour la topologie $I$-adique.
Soit }

$$
\xymatrix{
\mathcal{X'}\ \ar[r]^{v} \ar[d]_{g} &\mathcal{X} \ar[d]^{f}\\
\mathcal{S'}\   \ar[r]^{u}&  \mathcal{S}
}
$$

\noindent \textit{un carr\'e cart\'esien de $\mathcal{V}$-sch\'emas formels (pour la topologie $I$-adique) de type fini, avec $f$ propre. }

\begin{enumerate}
\item[(1.1.1)] \textit{Soit $\mathcal{F}$ un $\mathcal{O}_{\mathcal{X}}$-module coh\'erent. Alors, pour tout entier $i \geqslant 0$, on a:}
	\begin{enumerate}
	\item[(1)]\textit{$R^i f_{\ast}(\mathcal{F})$ est un $\mathcal{O}_{\mathcal{S}}$-module coh\'erent. 
	\item[(2)] Supposons de plus $u$ plat; alors le morphisme de changement de base}
	$$u^{\ast}  R^{i}  f_{\ast}(\mathcal{F}) \longrightarrow R^{i}  g_{\ast} v^{\ast} (\mathcal{F})$$
	\textit{est un isomorphisme.}
	\end{enumerate}

\noindent \textit{Plus g\'en\'eralement on a:}
\item[(1.1.2)] \textit{Soient $\mathcal{E}^{\bullet}_{\mathcal{X}}$ un complexe born\'e de $\mathcal{O}_{\mathcal{X}}$-modules coh\'erents et  $\mathcal{E}^{\bullet}_{\mathcal{X'}} = \mathcal{E}^{\bullet}_{\mathcal{X}} \otimes_{\mathcal{O}_{\mathcal{S}}} \mathcal{O}_{\mathcal{S}'}$.  Alors, pour tout entier $i \geqslant 0$, on a:}
	\begin{enumerate}
	\item[(1)]\textit{$R^{i} f_{\ast} \mathcal{(E}^{\bullet}_{\mathcal{X}})$ est un $\mathcal{O}_{\mathcal{S}}$-module coh\'erent. 
	\item[(2)] Supposons de plus $u$ plat; alors le morphisme de changement de base}
	$$u^{\ast}  R^{i}  f_{\ast}(\mathcal{E}^{\bullet}_{\mathcal{X}})\  \tilde{\longrightarrow}\  R^{i}  g_{\ast}  	(\mathcal{E}^{\bullet}_{\mathcal{X'}}) $$
	\textit{est un isomorphisme, qui s'interpr\`ete aussi comme un isomorphisme}
	$$\displaystyle \mathop{\lim}_{\leftarrow\atop{n}} u^{\ast}_{n}  R^{i} f_{n \ast}  (\mathcal{E}^{\bullet}_	{X_{n}}) \tilde{\longrightarrow}\  \displaystyle \mathop{\lim}_{\leftarrow\atop{n}} R^{i} g_{n \ast}(\mathcal{E}^{\bullet}_{X'_{n}})\ ;$$
	\textit{(cf. notations plus bas).}
	\end{enumerate}
\end{enumerate}

\newpage
\vskip 3mm
\noindent \textit{D\'emonstration}.

\textbf{Pour (1.1.1)(1)}. La coh\'erence de $R^{i}  f_{\ast}(\mathcal{F})$ sur $\mathcal{O}_{\mathcal{S}}$ est rappel\'ee pour m\'emoire [EGA III, (3.4.2)]. \\

\textbf{Pour (1.1.1)(2)}. Notons $\mathcal{F}_{n} = \mathcal{F} / I^{n+1} \mathcal{F}$  et $\varphi_{n}$ le morphisme canonique $\varphi_{n} : R^{i} f_{\ast}(\mathcal{F}) \rightarrow R^{i} f_{\ast}(\mathcal{F}_{n})$.  Posons $C_{n} = \mbox{Coker}\ \varphi_{n},\  \mathcal{H} = R^{i} f_{\ast}(\mathcal{F}), \mathcal{H}_{n} = R^{i} f_{\ast}(\mathcal{F}_{n})$ et $\mathcal{H}'_{n} = \mathcal{H} / I^{n+1} \mathcal{H}$. Comme $I^{n+1} \mathcal{H} \subset \mbox{Ker}\  \varphi_{n}$ on a une surjection

$$\mathcal{H}'_{n}\    \twoheadrightarrow \mathcal{H} / \mbox{Ker}\  \varphi_{n}$$

\noindent de noyau not\'e $K_{n}$. Dans les suites exactes\\

(1.1.1.1) $\qquad 0 \longrightarrow K_{n} \longrightarrow \mathcal{H}'_{n} \longrightarrow \mathcal{H} / \mbox{Ker}\ \varphi_{n} \longrightarrow 0$\\

(1.1.1.2) $\qquad 0 \longrightarrow \mathcal{H} / \mbox{Ker}\ \varphi_{n} \longrightarrow \mathcal{H}_{n}
\longrightarrow C_{n} \longrightarrow 0$,\\

\noindent les syst\`emes projectifs $(\mathcal{H}'_{n})_{n}$ et $(\mathcal{H} / \mbox{Ker}\ \varphi_{n})_{n}$ v\'erifient la condition de Mittag-Leffler (not\'ee M-L) car les fl\`eches de transition sont surjectives. De plus $(\mathcal{H}_{n})_{n}$ v\'erifie M-L d'apr\`es [EGA III, (3.4.3)], donc $(C_{n})_{n}$ aussi [EGA ${O_{III}}$, (13.2.1)]. D'o\`u l'exactitude de la suite\\

(1.1.1.3) $\qquad  0 \longrightarrow \displaystyle \mathop{\lim}_{\leftarrow\atop{n}}\ \mathcal{H} / \mbox{Ker}\ \varphi_{n} \longrightarrow \displaystyle \mathop{\lim}_{\leftarrow\atop{n}}\ \mathcal{H}_{n}\
 \longrightarrow \displaystyle \mathop{\lim}_{\leftarrow\atop{n}}\  C_{n} \longrightarrow 0\ .$
 
 \noindent D'autre part on a un isomorphisme [EGA III, (3.4.3)]
 
 $$R^{i} f_{\ast}(\mathcal{F}) \tilde{\longrightarrow}  \displaystyle \mathop{\lim}_{\leftarrow\atop{n}}\  \mathcal{H}_{n}\  ;$$
 comme $R^{i} f_{\ast}(\mathcal{F})$ est un $\mathcal{O}_{\mathcal{S}}$-module coh\'erent [EGA III, (3.4.2)], il est s\'epar\'e et complet pour la topologie $I$-adique, donc on a aussi un isomorphisme
 
 $$R^{i} f_{\ast}(\mathcal{F}) \tilde{\longrightarrow}  \displaystyle \mathop{\lim}_{\leftarrow\atop{n}}\  \mathcal{H}'_{n}. $$ 
 Ainsi le morphisme compos\'e
 
 $$\displaystyle \mathop{\lim}_{\leftarrow\atop{n}}\ \mathcal{H}'_{n} \longrightarrow \displaystyle \mathop{\lim}_{\leftarrow\atop{n}} \mathcal{H} /  \mbox{Ker}\ \varphi_{n} \hookrightarrow   \displaystyle \mathop{\lim}_{\leftarrow\atop{n}}\ \mathcal{H}_{n}$$
 est un isomorphisme, donc chacune des fl\`eches est un isomorphisme (car la seconde est injective par (1.1.1.3)). Par suite $\displaystyle \mathop{\lim}_{\leftarrow\atop{n}}\ C_{n} = 0$ ; comme $(C_{n})_{n}$ v\'erifie M-L, on en d\'eduit que le pro-objet $\ll (C_{n})_{n} \gg$ associ\'e est le pro-objet nul [G , 195-03, \S 2], i.e. \\
 
 (1.1.1.4)\qquad $\forall n$,  $\exists n' \geqslant n$ tel que $C_{n'} \rightarrow C_{n}$ soit la fl\`eche nulle\\
 
 \noindent (cf aussi [Bour, TG II, \S3, \no 5, th\'eo 1]).\\
 D'apr\`es [EGA III, (3.4.4)] les Ker $\varphi_{n}$ d\'efinissent sur $\mathcal{H}$ une filtration $I$-bonne, i.e. [EGA $0_{III}$, (13.7.7)], on a $I \mbox{Ker}\ \varphi_{n}\subset \mbox{Ker}\ \varphi_{n+1}$ avec \'egalit\'e pour $n$ assez grand; en particulier on a $\mbox{Ker}\ \varphi_{n}\subset I\mathcal{H}$ pour $n$ assez grand et la topologie sur $\mathcal{H}$ d\'efinie par les $\mbox{Ker}\ \varphi_{n}$ co\"{\i}ncide avec la topologie $I$-adique: par platitude de $u$, la topologie sur $\mathcal{O}_{\mathcal{S'}}\otimes_{\mathcal{O}_{\mathcal{S}}}\mathcal{H}$ d\'efinie par les $\mathcal{O}_{\mathcal{S'}}\otimes_{\mathcal{O}_{\mathcal{S}}}\mbox{Ker}\ \varphi_{n}$  co\"{\i}ncide donc avec la topologie $I$-adique. Comme $\mathcal{H}$ est coh\'erent sur $\mathcal{O}_{\mathcal{S}}$ on en d\'eduit des isomorphismes\\
 
$\mathcal{O}_{\mathcal{S'}}\otimes_{\mathcal{O}_{\mathcal{S}}}\mathcal{H}  \tilde{\longrightarrow}\  \displaystyle \mathop{\lim}_{\leftarrow\atop{n}}\ ( \mathcal{O}_{\mathcal{S'}}\otimes_{\mathcal{O}_{\mathcal{S}}}\mathcal{H}/\  I^{n+1}(\mathcal{O}_{\mathcal{S'}}\otimes_{\mathcal{O}_{\mathcal{S}}}\mathcal{H}))$

$   \quad  \qquad \qquad   \tilde{\longrightarrow}\  \displaystyle \mathop{\lim}_{\leftarrow\atop{n}}\  (\mathcal{O}_{\mathcal{S'}}\otimes_{\mathcal{O}_{\mathcal{S}}}\mathcal{H}/\  \mathcal{O}_{\mathcal{S'}}\otimes_{\mathcal{O}_{\mathcal{S}}}\mbox{Ker}\ \varphi_{n})$

$   \quad  \qquad \qquad   \tilde{\longrightarrow}\  \displaystyle \mathop{\lim}_{\leftarrow\atop{n}}\  (\mathcal{O}_{\mathcal{S'}}\otimes_{\mathcal{O}_{\mathcal{S}}}(\mathcal{H}/\mbox{Ker}\ \varphi_{n}))\ .$

 \noindent Or il r\'esulte de (1.1.1.4) que\\
 
\noindent $\forall n$,  $\exists n' \geqslant n$ tel que $\mathcal{O}_{\mathcal{S'}}\otimes_{\mathcal{O}_{\mathcal{S}}}C_{n'} \rightarrow \mathcal{O}_{\mathcal{S'}}\otimes_{\mathcal{O}_{\mathcal{S}}}C_{n}$ soit la fl\`eche nulle,\\
 
 \noindent donc $\displaystyle \mathop{\lim}_{\leftarrow\atop{n}}\ ( \mathcal{O}_{\mathcal{S'}}\otimes_{\mathcal{O}_{\mathcal{S}}}C_{n})=0$; comme les fl\`eches de transition du syst\`eme projectif $ \{\mathcal{O}_{\mathcal{S'}}\otimes_{\mathcal{O}_{\mathcal{S}}}(\mathcal{H}/\mbox{Ker}\ \varphi_{n})\}_{n}$ sont surjectives, ce syst\`eme projectif v\'erifie M-L: la suite exacte (1.1.1.2) le reste apr\`es tensorisation par $\mathcal{O}_{\mathcal{S'}}$ sur $\mathcal{O}_{\mathcal{S}}$ et en passant \`a la $\displaystyle \mathop{\lim}_{\leftarrow\atop{n}}$ on d\'eduit de ce qui pr\'ec\`ede un isomorphisme
 $$
  \displaystyle \mathop{\lim}_{\leftarrow\atop{n}}\  (\mathcal{O}_{\mathcal{S'}}\otimes_{\mathcal{O}_{\mathcal{S}}}(\mathcal{H}/\mbox{Ker}\ \varphi_{n}))\  \tilde{\longrightarrow}\  \displaystyle \mathop{\lim}_{\leftarrow\atop{n}}\  (\mathcal{O}_{\mathcal{S'}}\otimes_{\mathcal{O}_{\mathcal{S}}}\mathcal{H}_{n})\ . 
 $$
\noindent D'o\`u l'isomorphisme 
$$
u^{\ast}(\mathcal{H}) \tilde{\longrightarrow}\  \displaystyle \mathop{\lim}_{\leftarrow\atop{n}}\  u^{\ast}_{n}(\mathcal{H}_{n})\ ,
\leqno (1.1.1.5)
$$

\noindent o\`u l'on note \\

$$
\xymatrix{
X'_{n}\ \ar[r]^{v_{n}} \ar[d]_{g_{n}} & X_{n} \ar[d]^{f_{n}}\\
S'_{n}\   \ar[r]_{u_{n}} &  S_{n}
}
$$

\noindent le carr\'e cart\'esien d\'eduit de celui de la proposition par r\'eduction mod $I^{n+1}$.\\
\noindent Or, $u_{n}$ \'etant plat, [EGA III, (1.4.15)] fournit un isomorphisme de changement de base

$$u^{\ast}_{n}(\mathcal{H}_{n}) =  u^{\ast}_{n}(R^{i} f_{n \ast}(\mathcal{F}_{n})) \tilde{\longrightarrow}\ R^{i} g_{n \ast} (v^\ast_{n}(\mathcal{F}_{n})),$$
d'o\`u, par passage \`a la limite et via (1.1.1.5), des isomorphismes

$$u^{\ast}(R^{i} f_{\ast}(\mathcal{F}))  \tilde{\longrightarrow}\  \displaystyle \mathop{\lim}_{\leftarrow\atop{n}}\ R^{i} g_{n \ast} (v^\ast_{n}(\mathcal{F}_{n})) = \displaystyle \mathop{\lim}_{\leftarrow\atop{n}}\ R^{i} g_{\ast} (v^\ast_{n}(\mathcal{F}_{n}))$$

$\qquad\qquad\ \  \quad \qquad \quad \tilde{\longrightarrow}\ R^{i} g_{\ast} (v^{\ast}(\mathcal{F})) $ [EGA III, (3.4.3)],\\

\noindent car $ v^\ast(\mathcal{F}) = \displaystyle \mathop{\lim}_{\leftarrow\atop{n}}\ v^\ast_{n}(\mathcal{F}_{n})$. D'o\`u le (1.1.1) du th\' eor\`eme.\\

 \textbf{Pour (1.1.2)(1)}. On note $\mathcal{V}_{n} = \mathcal{V}/I^{n+1}$, $\mathcal{E}^{\bullet}_{X_{n}} = \mathcal{E}^{\bullet}_{\mathcal{X}} / I^{n+1} \mathcal{E}^{\bullet}_{\mathcal{X}}$ et $\mathcal{E}^{\bullet}_{X'_{n}} = \mathcal{E}^{\bullet}_{\mathcal{X}'} / I^{n+1}\ \mathcal{E}^{\bullet}_{\mathcal{X}'}$. Il suffit d'\'etablir le lemme suivant:\\

\noindent \textbf{Lemme (1.1.2.1)}.
\textit{Sous les hypoth\`eses  (1.1.2) et avec les notations ci-dessus on a un isomorphisme de $\mathcal{O}_{\mathcal{S}}$-modules coh\'erents:}
$$
  R^{i}  f_{\ast}(\mathcal{E}^{\bullet}_{\mathcal{X}})\  \tilde{\longrightarrow}\  
	\displaystyle \mathop{\lim}_{\leftarrow\atop{n}}  R^{i} f_{n \ast}  (\mathcal{E}^{\bullet}_{X_{n}})\ . 
$$
\vskip 3mm
\noindent \textit{D\'emonstration du lemme}. Il s'agit d'\'etendre la preuve de [EGA III, (3.4.4)] des $\mathcal{O}_{\mathcal{X}}$-modules coh\'erents aux complexes de $\mathcal{O}_{\mathcal{X}}$-modules coh\'erents.\\
  Comme $R^{i+j} f_{n \ast}\ (\mathcal{E}^{\bullet}_{{X}_{n}})$ est l'aboutissement d'une suite spectrale de terme $E^{i,j}_{1}$ donn\'e par
  $$E^{i,j}_{1} = R^{j}  f_{n \ast} (\mathcal{E}^i_{{X}_{n}})$$
  \noindent et que $E^{i,j}_{1}$ est coh\'erent sur $\mathcal{O}_{{S}_{n}}$ puisque $f_{n}$ est propre, on en d\'eduit que $R^{i+j} f_{n \ast}\ (\mathcal{E}^{\bullet}_{{X}_{n}})\ $ est coh\'erent sur $\mathcal{O}_{{S}_{n}}$: plus g\'en\'eralement, par la m\^eme m\'ethode que pour la preuve de [EGA III, (3.4.4)], on prouve que la condition $(F_{i})$ de [EGA $0_{III}$, (13.7.7)] est v\'erifi\'ee pour tout $i$, donc par [loc. cit.] on a:
 \begin{enumerate}
 \item[(1.1.2.2)] Pour tout $i$ le syst\`eme projectif $(R^{i} f_{n \ast}  (\mathcal{E}^{\bullet}_{X_{n}}))_{n}$ v\'erifie M-L.
 \item[(1.1.2.3)]$ \displaystyle \mathop{\lim}_{\leftarrow\atop{n}}  R^{i} f_{n \ast}  (\mathcal{E}^{\bullet}_{X_{n}})$ est un $\mathcal{O}_{\mathcal{S}}$-module coh\'erent.
 \item[(1.1.2.4)] La filtration sur $R^{i}  f_{\ast}(\mathcal{E}^{\bullet}_{\mathcal{X}})$ d\'efinie par
 $$
 \mbox{Ker}\{ \psi_{n}: R^{i}  f_{\ast}(\mathcal{E}^{\bullet}_{\mathcal{X}})\rightarrow   R^{i} f_{ \ast}  (\mathcal{E}^{\bullet}_{X_{n}})  \} 
 $$
 \noindent est $I$-bonne.

 \end{enumerate}
 
 \noindent Or $R^{i+j}  f_{\ast}(\mathcal{E}^{\bullet}_{\mathcal{X}})$ est l'aboutissement d'une suite spectrale de terme $E^{i,j}_{1} = R^{j}  f_{ \ast} (\mathcal{E}^i_{\mathcal{X}})$ qui est un $\mathcal{O}_{\mathcal{S}}$-module coh\'erent via (1.1.1), donc $R^{i+j}  f_{\ast}(\mathcal{E}^{\bullet}_{\mathcal{X}})$ est un $\mathcal{O}_{\mathcal{S}}$-module coh\'erent; comme la topologie $I$-adique sur $R^{i}  f_{\ast}(\mathcal{E}^{\bullet}_{\mathcal{X}})$ co\"{\i}ncide d'apr\`es (1.1.2.4) avec la topologie d\'efinie par les $ \mbox{Ker}\{ \psi_{n}\}$, on en d\'eduit des isomorphismes
 $$
  R^{i}  f_{\ast}(\mathcal{E}^{\bullet}_{\mathcal{X}})\  \tilde{\longrightarrow}\  
	\displaystyle \mathop{\lim}_{\leftarrow\atop{n}} ( R^{i} f_{ \ast}  (\mathcal{E}^{\bullet}_{\mathcal{X}})/I^{n+1}R^{i} f_{ \ast}  (\mathcal{E}^{\bullet}_{\mathcal{X}}))  \tilde{\longrightarrow}\ \displaystyle \mathop{\lim}_{\leftarrow\atop{n}}  (R^{i} f_{ \ast}  (\mathcal{E}^{\bullet}_{\mathcal{X}})/\mbox{Ker}\{ \psi_{n}\})\ .
$$

\noindent En passant \`a la limite projective sur les injections canoniques
$$
(R^{i} f_{ \ast}  (\mathcal{E}^{\bullet}_{\mathcal{X}})/\mbox{Ker}\{ \psi_{n}\}) \hookrightarrow R^{i} f_{ \ast}  (\mathcal{E}^{\bullet}_{X_{n}})
$$
\noindent il en r\'esulte une injection canonique
$$
R^{i} f_{ \ast}  (\mathcal{E}^{\bullet}_{\mathcal{X}}) \hookrightarrow \displaystyle \mathop{\lim}_{\leftarrow\atop{n}}R^{i} f_{ \ast}  (\mathcal{E}^{\bullet}_{X_{n}})\ .
$$

 Montrons \`a pr\'esent la surjectivit\'e de cette fl\`eche. Puisque $R^{i}  f_{\ast}(\mathcal{E}^{\bullet}_{\mathcal{X}})$ (resp. $\displaystyle \mathop{\lim}_{\leftarrow\atop{n}}R^{i} f_{ \ast}  (\mathcal{E}^{\bullet}_{X_{n}})$) est le faisceau associ\'e au pr\'efaisceau $U \mapsto H^{i}(f^{-1}(U), \mathcal{E}^{\bullet}_{\mathcal{X}})$ (resp. $U \mapsto \displaystyle \mathop{\lim}_{\leftarrow\atop{n}}H^{i}(f^{-1}(U), \mathcal{E}^{\bullet}_{{X}_{n}})$), il suffit de montrer que pour tout ouvert $U$ de $\mathcal{S}$ on a une surjection canonique
 
 $$
  H^{i}(f^{-1}(U), \mathcal{E}^{\bullet}_{\mathcal{X}}) \twoheadrightarrow  \displaystyle \mathop{\lim}_{\leftarrow\atop{n}}H^{i}(f^{-1}(U), \mathcal{E}^{\bullet}_{{X}_{n}})\ .
 $$
 
\noindent Pour $m\geqslant n$, notons $\mathcal{N}^{\bullet}_{m,n}$ le noyau de la projection $\mathcal{E}^{\bullet}_{X_{m}} \rightarrow \mathcal{E}^{\bullet}_{X_{n}}$ et choisissons des r\'esolutions injectives [C-E, chap. XVII, \S1] $\mathcal{I}^{\bullet \bullet}_{m,n}\ , \ \mathcal{I}^{\bullet\bullet}_{m}\ , \ \mathcal{I}^{\bullet\bullet}_{n}$ de $\mathcal{N}^{\bullet}_{m,n}\ ,\  \mathcal{E}^{\bullet}_{X_{m}}\ ,\  \mathcal{E}^{\bullet}_{X_{n}}$ respectivement telles que l'on ait un diagramme commutatif \`a lignes et colonnes exactes [C-E, chap. V, \S2] et morphismes d'augmentation les $\varepsilon$

$$
\begin{array}{c}
\xymatrix{
&0\ar[d]&0\ar[d]&0\ar[d]&\\
0\ar[r]&\mathcal{N}^{\bullet}_{m,n} \ar[r] \ar[d]^{\varepsilon_{m,n}}  &\mathcal{E}^{\bullet}_{X_{m}} \ar[d]^{\varepsilon_{m}} \ar[r]& \mathcal{E}^{\bullet}_{X_{n}}\ar[d]^{\varepsilon_{n}}\ar[r]&0\\
0\ar[r]&\mathcal{I}^{\bullet \bullet}_{m,n} \ar [r] & \mathcal{I}^{\bullet \bullet}_{m} \ar[r]&\mathcal{I}^{\bullet \bullet}_{n}\ar[r]&0 & .
} 
\end{array}
$$
 \noindent Alors, pour $U$ ouvert de $\mathcal{S}$, la suite 
 $$
\begin{array}{c}
\xymatrix{
0\ar[r]&\Gamma(f^{-1}(U), \mathcal{I}^{\bullet \bullet}_{m,n}) \ar[r]  &\Gamma(f^{-1}(U), \mathcal{I}^{\bullet \bullet}_{m}) \ar[r]&\Gamma(f^{-1}(U), \mathcal{I}^{\bullet \bullet}_{n})\ar[r]&0\\
} 
\end{array}
$$
 est exacte; en particulier le syst\`eme projectif $\{\Gamma(f^{-1}(U), \mathcal{I}^{\bullet \bullet}_{n})\}_{n}$ a des fl\`eches de transition surjectives, d'o\`u $R^{j} \displaystyle \mathop{\lim}_{\leftarrow\atop{n}}\Gamma(f^{-1}(U), \mathcal{I}^{\bullet \bullet}_{n})=0$ pour $j>0$ [J, prop 2.1 et th\'eo 1.8]. Pour la m\^eme raison on a aussi $R^{j} \displaystyle \mathop{\lim}_{\leftarrow\atop{n}} \mathcal{E}^{\bullet}_{X_{n}}=0$ pour $j>0$. On en d\'eduit donc des quasi-isomorphismes
  $$
 \mathbb{R}\Gamma(f^{-1}(U),\mathcal{E}^{\bullet}_{\mathcal{X}})\simeq \mathbb{R}\Gamma(f^{-1}(U), \mathbb{R} \displaystyle \mathop{\lim}_{\leftarrow\atop{n}} \mathcal{E}^{\bullet}_{X_{n}})\simeq \mathbb{R} \displaystyle \mathop{\lim}_{\leftarrow\atop{n}}\mathbb{R}\Gamma(f^{-1}(U),  \mathcal{E}^{\bullet}_{X_{n}})$$
$$ \simeq\mathbb{R}  \displaystyle \mathop{\lim}_{\leftarrow\atop{n}}\Gamma(f^{-1}(U),  \mathcal{I}^{\bullet\bullet}_{n})\simeq\displaystyle \mathop{\lim}_{\leftarrow\atop{n}}\Gamma(f^{-1}(U),  \mathcal{I}^{\bullet\bullet}_{n})\ .
   $$
\noindent Ainsi, pour tout entier $i\geqslant0$, le morphisme canonique $\theta$ \\
$$
H^{i}(f^{-1}(U),\mathcal{E}^{\bullet}_{\mathcal{X}})\overset{\sim}{\longrightarrow}\mathcal{H}^{i}( \mathbb{R}\Gamma(f^{-1}(U),\mathcal{E}^{\bullet}_{\mathcal{X}}))\overset{\sim}{\longrightarrow}$$
$$\mathcal{H}^{i}(\displaystyle \mathop{\lim}_{\leftarrow\atop{n}}\Gamma(f^{-1}(U),  \mathcal{I}^{\bullet\bullet}_{n}))\overset{\theta}{\longrightarrow} \displaystyle \mathop{\lim}_{\leftarrow\atop{n}}\mathcal{H}^{i}(\Gamma(f^{-1}(U),  \mathcal{I}^{\bullet\bullet}_{n}))\overset{\sim}{\longrightarrow}\displaystyle \mathop{\lim}_{\leftarrow\atop{n}}H^{i}(f^{-1}(U),\mathcal{E}^{\bullet}_{X_{n}})$$
  \noindent est surjectif d'apr\`es [EGA $0_{III}$, (13.2.3)]; en fait [loc. cit.] fournit aussi l'injectivit\'e de $\theta$ puisque le syst\`eme projectif $\{H^{i}(f^{-1}(U),\mathcal{E}^{\bullet}_{X_{n}})\}_{n}$ v\'erifie M-L d'apr\`es (1.1.2.2);
  d'o\`u le lemme. $\square$\\
 
 \textbf{Pour (1.1.2)(2)}.  Au cran fini $n$ le morphisme de changement de base par $u_{n}$\\
 
(1.1.2.5) \qquad\qquad$
 u^{\ast}_{n}  R^{i+j} f_{n \ast}  (\mathcal{E}^{\bullet}_	{X_{n}}) \longrightarrow\   R^{i+j} g_{n \ast}(\mathcal{E}^{\bullet}_{X'_{n}})
 $\\
 
\noindent est un isomorphisme car c'en est un au niveau des termes $E^{i,j}_{1}$ de chaque membre puisque $u_{n}$ est plat. Compte tenu de (1.1.2.4) le passage \`a la $\displaystyle \mathop{\lim}_{\leftarrow\atop{n}} $ sur (1.1.2.5) fournit, par la m\^eme m\'ethode qu'au cas (1.1.1)(2) du th\'eor\`eme (1.1), l'isomorphisme de changement de base 
$$u^{\ast}  R^{i}  f_{\ast}(\mathcal{E}^{\bullet}_{\mathcal{X}})\  \tilde{\longrightarrow}\  R^{i}  g_{\ast}  	(\mathcal{E}^{\bullet}_{\mathcal{X'}}) $$
	\noindent qui s'interpr\`ete aussi comme un isomorphisme
	$$\displaystyle \mathop{\lim}_{\leftarrow\atop{n}} u^{\ast}_{n}  R^{i} f_{n \ast}  (\mathcal{E}^{\bullet}_	{X_{n}}) \tilde{\longrightarrow}\  \displaystyle \mathop{\lim}_{\leftarrow\atop{n}} R^{i} g_{n \ast}(\mathcal{E}^{\bullet}_{X'_{n}})\ . \ \square$$

\vskip 3mm
\noindent \textbf{Th\'eor\`eme (1.2)}.
\textit{Soient $\mathcal{V}$ un anneau de valuation discr\`ete s\'epar\'e et complet pour la topologie $\mathfrak{m}$-adique, o\`u $\mathfrak{m}$ est son id\'eal maximal, $k = \mathcal{V}/ \mathfrak{m}$ son corps r\'esiduel suppos\'e de caract\'eristique $p>0$, $K$ son corps des fractions de caract\'eristique $0$.}\\

\noindent \textit{Soit}

$$
\xymatrix{
X'\ \ar[r]^{v} \ar[d]_{g} & X \ar[d]^{f}\\
S'\   \ar[r]_{u} &  S
}
$$

\noindent \textit{un carr\'e cart\'esien de $K$-espaces analytiques rigides, avec $f$ propre.}\\

 \noindent (1.2.1) \textit{Soit $E$ un $\mathcal{O}_{X}$-module coh\'erent. Alors, pour tout entier $i \geqslant 0$},
 
 \begin{enumerate}
 \item[(1)] \textit{$R^{i} f_{\ast}(E)$ est un $\mathcal{O}_{S}$-module coh\'erent.}
\item[(2)] \textit{Supposons de plus u plat; alors le morphisme de changement de base}

$$u ^{\ast} (R^{i} f_{\ast} (E)) \longrightarrow R^{i} g_{\ast} (v^{\ast}(E))$$
\textit{est un isomorphisme.}
\end{enumerate}

\noindent \textit{Plus g\'en\'eralement on a:}\\
\noindent (1.2.2)  \textit{Soient $E^{\bullet}$ un complexe born\'e de $\mathcal{O}_{X}$-modules coh\'erents et $E'^{\bullet} = E^{\bullet}    \otimes_{\mathcal{O}_{\mathcal{S}}} \mathcal{O}_{S'} $.  Alors pour tout entier $i \geqslant 0$}

\begin{enumerate}
\item[(1)] \textit{$R^i{} f_{\ast}(E^{\bullet})$ est un $\mathcal{O}_{S}$-module coh\'erent.}
\item[(2)] \textit{Supposons de plus u plat; alors le morphisme de changement de base
$$u^{\ast} R^{i} f_{\ast}(E^{\bullet})\ \longrightarrow R^{i} g_{\ast}(E'^{\bullet})$$
est un isomorphisme. }
\end{enumerate}

\vskip 10mm
\noindent \textit{D\'emonstration du th\'eor\`eme (1.2)}.\\

\noindent \textbf{1\`ere \'etape}. Supposons d'abord d\'emontr\'e le th\'eor\`eme dans le cas o\`u les $K$-espaces analytiques rigides sont tous quasi-compacts et quasi-s\'epar\'es. \\

Prouvons alors (1.2.1) dans le cas g\'en\'eral.\\

 L'assertion (1) est locale sur $S$ puisque $R^{i} f_{\ast}(E)$ est le faisceau associ\'e au pr\'efaisceau

$$V \mapsto\ H^i(f^{-1} (V), E)$$

\noindent o\`u $V$ parcourt les ouverts de $S$ [SGA 4, V, prop 5.1].

Soient $\psi : V \hooklongrightarrow S$ un ouvert affino¬\"{\i}de de $S$ et $W$ d\'efini par le carr\'e cart\'esien

$$
\begin{array}{c}
\xymatrix{
W\ \ar[r]^{\theta} \ar[d]_{f'} & X \ar[d]^{f}  &\\ 
V\   \ar[r]_{\psi} &  S  & \quad;
}
\end{array}
\leqno(1.2.1.1)  
$$

\noindent d'apr\`es [loc. cit.] on a alors un isomorphisme canonique

$$\psi^\ast  R^{i}  f_{\ast}(E)\   \cong\ R^{i}  f'_{\ast}(\theta^{\ast}(E)).$$

\noindent Or ici $V$ est quasi-compact, quasi-s\'epar\'e, et $W$ aussi car $f'$ est propre ; d'o\`u l'assertion (1) via le cas quasi-compact, quasi-s\'epar\'e.\\

Pour le (2) on reprend le carr\'e cart\'esien (1.2.1.1) : soient 

$$u' : V' = V\times_{S}  S' \longrightarrow\ V $$
$$\mbox{et}\ u'' : V'' \hookrightarrow V'\  \mbox{un ouvert affino¬\"{\i}} \mbox{de de}\  V',$$

\noindent et on consid\`ere le diagramme commutatif suivant

$$
\xymatrix{& X'\ar@{.>}[dd]^(.6){g} |\hole\ar@{=} [rr] & & X' \ar@{.>}[dd]_(.6){g} |\hole \ar[rr]^{v}& &
X \ar[dd]^f \\
W'' \ar[rr]^(.8){v''} \ar[dd]_{g''} \ar[ur]^{\theta''} & & W' \ar[rr]^(.8){v'} \ar[dd]_(.7){g'} \ar[ur]^{\theta'}  & & W \ar[dd]_(.7){f'} \ar[ur]_{\theta} & \\
& S' \ar@{==}[rr]  |\hole & & S' \ar@{.>}[rr]^(.8){u} |\hole & &  S  \\
V'' \ar@{^{(}->}[rr]^{u''} \ar@{.>}[ur]^{\psi''} & & V' \ar[rr]^{u'} \ar@{.>}[ur]_{\psi'} & & V \ar[ur]_{\psi} 
}
$$

\vskip 2mm
\noindent dont les faces verticales sont cart\'esiennes.\\

Vu le caract\`ere local de l'isomorphisme (2) cherch\'e, il suffit de montrer que l'on a un isomorphisme

$$\psi''^{\ast} u^{\ast} R^{i} f_{\ast}(E) \displaystyle \mathop{\longrightarrow}^\sim \psi''^{\ast} R^{i} g_{\ast}\  v^{\ast}(E).$$

\noindent Or ici $V$ et $V''$ sont quasi-compacts, quasi-s\'epar\'es, donc $W$ et $W''$ aussi car $f'$ et $g''$ sont propres et on peut appliquer l'isomorphisme de changement de base du cas quasi-compact, quasi-s\'epar\'e pour le carr\'e cart\'esien

$$
\xymatrix{
W'\ \ar[r]^{v'\circ v''} \ar[d]_{g''} & W \ar[d]^{f'} &\\
V''\   \ar[r]_{u' \circ u''}  & V & \quad;
}
$$

\noindent d'o\`u une suite d'isomorphismes

$$\psi''^{\ast}  u^{\ast} R^{i} f_{\ast}(E) =  (u' u'')^{\ast} \psi^{\ast} R^{i} f_{\ast}(E)$$

$\qquad  \qquad  \qquad \qquad  \qquad  \qquad \qquad = (u' u'')^{\ast}  R^{i} f'_{\ast}   (\theta^{\ast}(E))$ 

$\qquad  \qquad  \qquad \qquad  \qquad  \qquad \qquad  \simeq R^{i} g''_{\ast} ((v' v'')^{\ast}(\theta^{\ast}(E)))$

$\qquad  \qquad  \qquad \qquad  \qquad  \qquad \qquad  = R^{i} g''_{\ast} (\theta''^{\ast}(v^{\ast}(E))) $

$\qquad  \qquad  \qquad \qquad  \qquad  \qquad \qquad  = \psi''^{\ast} R^{i} g_{\ast} (v^{\ast}(E)).$

\noindent D'o\`u le (2).

Pour prouver (1.2.2) \`a partir du cas quasi-compact, quasi-s\'epar\'e la d\'emarche est analogue.\\

\noindent \textbf{2\`eme \'etape}. Prouvons le th\'eor\`eme dans le cas quasi-compact, quasi-s\'epar\'e.\\

Pour (1.2.1). Le (1) est un th\'eor\`eme de L¬\"utkebohmert [L¬\"u, theo 2.7].\\
Pour le (2) on adopte les notations de [Bo-L¬\"u 1, d\'emonstration de 4.1] : d'apr\`es [loc. cit.] il existe des mod\`eles formels $\stackrel{\circ}{\mathcal{X}}$, $\mathcal{S}$, $\stackrel{\circ}{\mathcal{S}'}$ de $X, S, S'$ respectivement et des \'eclatements admissibles

$$\tau_{\mathcal{X}} : \mathcal{X} \rightarrow\  \stackrel{\circ}{\mathcal{X}} \quad , \quad  \tau_{\mathcal{S}'} : \mathcal{S}'  \rightarrow\  \stackrel{\circ}{\mathcal{S}'} $$
et des morphismes

$$\varphi : \mathcal{X} \rightarrow \mathcal{S} \quad , \quad \theta : \mathcal{S}' \rightarrow \mathcal{S} $$
tels que

$$\varphi_{\mbox{rig}} = f\ \circ\  \tau_{\mathcal{X} \mbox{rig}} \quad \mbox{et} \quad  \theta_{\mbox{rig}} = u\ \circ\  \tau_{\mathcal{S}' \mbox{rig}},$$

\noindent avec $\theta$ plat [Bo-L\" u 2, 5.2, 5.10 (c)].\\

Remarquons que tous les sch\'emas formels pr\'ec\'edents sont admissibles au sens de [Bo - L¬\"u 1], donc sont plats sur $\mathcal{V}$ [Bo - L¬\"u 1, \S\ 1].

\noindent On dispose donc d'un carr\'e cart\'esien

$$
\xymatrix{
\mathcal{X}' \ar[r]^{\theta'} \ar[d]_{\varphi'} & \mathcal{X} \ar[d]^{\varphi} &\\
\mathcal{S'}  \ar[r]_{\theta} &  \mathcal{S} & \quad ,
}
$$

\noindent qui est un mod\`ele formel du carr\'e du th\'eor\`eme : de plus il existe un $\mathcal{O}_{\mathcal{X}}$-module coh\'erent $\mathcal{F}$ tel que $\mathcal{F}_{\mbox{rig}} = E$ [L¬\"u, 2.2], ou [Bo-L¬\"u 1, 5.6] et $\varphi$ est un morphisme propre [L¬\"u, 2.6]. D'apr\`es le th\'eor\`eme (1.1) le morphisme de changement de base

$$\theta^{\ast} R^{i} \varphi_{\ast}(\mathcal{F}) \rightarrow R^{i} \varphi '_{\ast}\ \theta '^{\ast}(\mathcal{F}) $$
est un isomorphisme ; par passage aux fibres g\'en\'eriques le (1.2.1) en r\'esulte.\\

Pour (1.2.2) : le complexe $E^\bullet$ provient, \`a multiplication pr\`es par $p^{-n}$, d'un complexe $\mathcal{E}^\bullet$ de $\mathcal{O}_{\mathcal{X}}$-modules coh\'erents. Il suffit alors d'appliquer le (1.1.2) du th\'eor\`eme (1.1) et de passer aux fibres g\'en\'eriques pour obtenir (1.2.2).\\
 
Nous rassemblons pour m\'emoire dans la proposition suivante quelques propri\'et\'es des immersions.
  
\vskip 3mm
\noindent \textbf{Proposition (1.3)}. \textit{Soient $\mathcal{V}$ et $K$ comme en (1.2). Alors}
\begin{enumerate}
\item[(1.3.1)]  \textit{Toute immersion (resp. immersion ouverte, resp. immersion ferm\'ee) $\alpha : X\hookrightarrow Y$ de $K$-espaces rigides analytiques quasi-compacts et quasi-s\'epar\'es admet un mod\`ele formel $\beta : \mathcal{X} \hookrightarrow \mathcal{Y}$ au sens de [Bo-L¬\"u 2, cor 5.10] qui est une immersion (resp. une immersion ouverte, resp. une immersion ferm\'ee) : en particulier $\mathcal{X}$ et $\mathcal{Y}$ sont plats sur $\mathcal{V}$.}

\item[(1.3.2)]  \textit{Pour tout $K$-espace analytique rigide $X$, le faisceau d'anneaux $\mathcal{O}_{X}$ est coh\'erent.}

\item[(1.3.3)] \textit{Soient $\alpha : X \hookrightarrow Y$ une immersion ferm\'ee de $K$-espaces analytiques rigides et $M$ un $\mathcal{O}_{X}$-module coh\'erent. Alors, pour tout entier $i > 0$ on a  }
$$R^{i} \alpha_{\ast}(M) = 0 $$

\noindent \textit{et le morphisme canonique}
$$\alpha^{\ast}\alpha_{\ast}(M) \rightarrow M $$

\noindent \textit{est un isomorphisme.}\\

\end{enumerate}

\vskip 3mm
\noindent \textit{D\'emonstration}. Le (1.3.1) n'est autre que [Bo-L¬\"u 2, cor 5.10].\\

Pour le (1.3.2) la propri\'et\'e est locale sur $X$ [B-G-R, 9.4.3] : on peut donc supposer $X$ affino¬\"{\i}de ; par suite $X$ est quasi-compact et quasi-s\'epar\'e et admet un mod\`ele formel $\mathcal{X}$.  Il suffit de prouver la coh\'erence du faisceau d'anneaux $\mathcal{O}_{\mathcal{X}}$ :  comme $\mathcal{X}$ est un sch\'ema formel de type fini sur l'anneau noeth\'erien  $\mathcal{V}$, le faisceau $\mathcal{O}_{\mathcal{X}}$ est un faisceau coh\'erent d'anneaux par [EGA I, (10.11.2)].\\

Pour le (1.3.3), les assertions sont locales sur $Y$ : on peut donc supposer $Y$ affino¬\"{\i}de, donc quasi-compact, quasi-s\'epar\'e et de m\^eme pour $X$ puisque $\alpha$ est propre [B-G-R, 9.5.3, prop 2, 9.6.2, prop 5]. On prend alors un mod\`ele formel $\beta : \mathcal{X} \hookrightarrow \mathcal{Y}$ de $\alpha$ [Bo-L\" u 2, 5.10(d)] et un $\mathcal{O}_{\mathcal{X}}$-module coh\'erent $\mathcal{M}$ tel que $\mathcal{M}_{K} = M$ [L\"u, lemma 2.2]. Soient $\beta_{n} : \mathcal{X}_{n} \hookrightarrow \mathcal{Y}_{n}$ la r\'eduction de $\beta$ mod $\mathfrak{m}^{n+1}$ et $\mathcal{M}_{n} = \mathcal{M} /\mathfrak{m}^{n+1}\ \mathcal{M}$ ; alors

$$ \beta^{\ast}  \beta_{\ast}(\mathcal{M}) /\mathfrak{m}^{n+1} = \beta^{\ast}_{n}\ \beta_{n\ast}(\mathcal{M}_{n})$$

\noindent et puisque $\mathcal{M}$ et $\beta_{\ast}(\mathcal{M})$ sont coh\'erents [th\'eo (1.1)] on a des isomorphismes

$$\mathcal{M} \overset{\sim}{\rightarrow} \displaystyle \mathop{\mbox{lim}}_{\longleftarrow \atop{n}}\  \mathcal{M}_{n} \quad \mbox{et} \quad  \beta^{\ast}  \beta_{\ast}(\mathcal{M}) \overset{\sim}{\rightarrow}  \displaystyle \mathop{\mbox{lim}}_{\longleftarrow \atop{n}}\   \beta^{\ast}_{n}\ \beta_{n\ast}(\mathcal{M}_{n}).$$

\noindent Or le morphisme canonique
$$\beta^{\ast}_{n}\ \beta_{n\ast}(\mathcal{M}_{n}) \longrightarrow \mathcal{M}_{n}$$

\noindent est un isomorphisme: en effet, $\mathcal{M}_{n}$ \'etant de pr\'esentation finie sur $\mathcal{O}_{\mathcal{X}_{n}}$, on est ramen\'e au cas de $\mathcal{O}_{\mathcal{X}_{n}}$ pour lequel l'assertion est claire puisque $\beta_{n}$ est une immersion ferm\'ee. En prenant la limite sur $n$ de ces isomorphismes et en passant aux fibres g\'en\'eriques on en d\'eduit  un isomorphisme canonique
$$\alpha^{\ast} \alpha_{\ast}(M) \overset{\sim}{\rightarrow} M \ .$$

\noindent  L'\'egalit\'e $R^{i} \alpha_{\ast}(M) = 0 $ pour $i > 0$ r\'esulte de l'isomorphisme
$$
 R^{i} \beta_{\ast}(\mathcal{M}) \overset{\sim}{\rightarrow}  \displaystyle \mathop{\mbox{lim}}_{\longleftarrow \atop{n}}\   R^{i}\beta_{n\ast}(\mathcal{M}_{n})\qquad [EGA III,(3.4.3)]
$$
\noindent et de l'\'egalit\'e $R^{i} \beta_{n\ast}(\mathcal{M}_{n}) = 0 $ pour $i > 0$, puisque $\beta_{n}$ est fini.
  $\square$

\newpage
\section*{2. Sorites sur les voisinages stricts}

\textbf{2.0.} Rappelons la d\'efinition de \guillemotleft voisinage strict\guillemotright\  dans un espace rigide analytique [G-K 2, 2.22].\\

Si $U$ est un ouvert admissible d'un espace rigide analytique $W$, un ouvert admissible $V \subset W$ est appel\'e voisinage strict de $U$ dans $W$ si $\{V, W \backslash U\}$ est un recouvrement admissible de $W$. Cette d\'efinition redonne celle de [B 3, (1.2.1)], [B 4, \S1], [LS, chap 2 et 3] dans le cas des tubes.\\

\textbf{2.1.} Consid\'erons un diagramme commutatif
$$
\begin{array}{c}
\xymatrix{
X \ar@{^{(}->}[r]^{i_{\mathcal{X}}} \ar[d]_{f}  \ar[d] & \mathcal{X} \ar[d]^{h} \ar[r]^{\psi} \ar @{} [dr] |{\square} & \mathcal{Y}\ar[d]^{\overline{h}}\\
S \ar@{^{(}->} [r]_{i_{\mathcal{S}}} & \mathcal{S} \ar[r]_{\varphi} & \mathcal{T}
} 
\end{array}
\leqno{(2.1.1)}
$$

\noindent dans lequel le carr\'e de droite est un carr\'e cart\'esien de $\mathcal{V}$-sch\'emas formels s\'epar\'es plats de type fini, $f$ est un morphisme de $k$-sch\'emas, $i_{\mathcal{X}}$, $i_{\mathcal{Y}}
 = \psi  \circ i_{\mathcal{X}}$, $i_{\mathcal{S}}$ et $i_{\mathcal{T}} = \varphi  \circ i_{\mathcal{S}}$ sont des immersions.\\
 
  \noindent On note
  $$
  h' :\ ]X[_{\mathcal{X}}\  \longrightarrow\  ]S[_{\mathcal{S}} ,\  \overline{h}' :\  ]X[_{\mathcal{Y}}\ 
\longrightarrow\  ]S[_{\mathcal{T}},
$$
$$ÊÊ \psi' :\ ]X[_{\mathcal{X}}\  \longrightarrow\  ]X[_{\mathcal{Y}},\ \varphi' :\ ]S[_{\mathcal{S}}\ \longrightarrow\  ]S[_{\mathcal{T}}
$$
 les morphismes induits sur les tubes respectivement par $h,\  \overline{h},\  \psi,\  \varphi$\ [B 3, (1.1.11) (i)], [B 4, 1.1], [LS, chap 2].

\vskip 3mm
\noindent \textbf{Proposition (2.1.2)}. \textit{ Avec les notations pr\'ec\'edentes, on a un diagramme commutatif \`a carr\'es cart\'esiens}

$$
\xymatrix{
]X[_{\mathcal{X}} \ar[r]^{\psi'} \ar[d]_{h'} & ]X[_{\mathcal{Y}} \ar[d]^{\overline{h}'}&\\
]S[_{\mathcal{S}} \ar[r]_{\varphi'}\ar@{^{(}->}[d]& ]S[_{\mathcal{T}}\ar@{^{(}->}[d]&\\
\mathcal{S}_{K} \ar[r]_{\varphi_{K}} & \mathcal{T}_{K}&.
}
$$

\vskip 3mm
\noindent \textit{D\'emonstration}. Pour le carr\'e du bas, \c ca r\'esulte de la d\'efinition des tubes et du morphisme de sp\'ecialisation. Pour le carr\'e du haut, il s'agit de v\'erifier que $]X[_{\mathcal{X}}$ satisfait la propri\'et\'e universelle du produit fibr\'e. Soit $Z$ un espace rigide analytique s'ins\'erant dans un diagramme commutatif

$$\xymatrix{
Z \ar[r]^{u}  \ar[d]  &  ]X[_{\mathcal{Y}} \ar[d] ^{\overline{h'}} \ar@{^{(}->}[r]  & \mathcal{Y}_{K} \ar[dd]^{\overline{h}_{K}} &\\
]S[_{\mathcal{S}} \ar[r]_{\varphi'} \ar@{^{(}->}[rd] &  ]S[ _{\mathcal{T}}\ar@{^{(}->}[rd] & &\\
& \mathcal{S}_{K} \ar[r]_{\varphi_{K}} & \mathcal{T}_{K} & \quad ;
}
$$

\noindent par propri\'et\'e universelle de $\mathcal{X}_{K} = \mathcal{S}_{K}  \times_{\mathcal{T}_{K}} \mathcal{Y}_{K}$ on en d\'eduit une fl\`eche $Z \displaystyle \mathop{\longrightarrow}^{v} \mathcal{X}_{K}$ qui s'ins\`ere dans le diagramme commutatif

$$
\xymatrix{
& X \ar@{=}[r]  \ar@{^{(}->}[d]^{i_{\mathcal{X}}} &  X  \ar@{^{(}->}[d]^{i_{\mathcal{Y}}} \\
& \mathcal{X} \ar[r]^{\psi}  & \mathcal{Y} \\
Z \ar[r]^{v}  \ar@/_4pc/[drr]_{u} & \mathcal{X}_{K} \ar[u]_{sp} \ar[r]^{\psi_{K}} & \mathcal{Y}_{K} \ar[u]_{sp} \\
& ]X[_{\mathcal{X}} \ar[r]^{\psi'} \ar@{^{(}->}[u] & ]X[_{\mathcal{Y}} \ar@{^{(}->}[u]
}
$$

\noindent o\`u $sp$ sont les morphismes de sp\'ecialisation. Puisque $sp (]X[_{\mathcal{Y}}) = X \displaystyle \mathop{\hookrightarrow}_{i_{\mathcal{Y}}} \mathcal{Y}$ et compte tenu de la commutativit\'e du diagramme pr\'ec\'edent le morphisme $v$ se factorise par $]X[_{\mathcal{X}}$ en

$$
\xymatrix{
Z \ar[r]^{v}  \ar[dr] & \mathcal{X}_{K} \\
& ]X[_{\mathcal{X}} \ar@{^{(}->}[u] &\quad ;
}
$$

\noindent d'o\`u la proposition. $\square$\\

$ \mathbf{2.2.}$  Consid\'erons \`a pr\'esent un diagramme commutatif
$$
\begin{array}{c}
$$\xymatrix{
X \ar@{^{(}->} [r]^{j_{X_{1}}} \ar [dr]_{f}& X_{1} \ar@{^{(}->}[r]^{j_{Y}} \ar[d]_{f_{1}}  \ar @{} [dr] |{\square} & Y \ar@{^{(}->}[r]^{i_{Y}} \ar[d]^{\overline{f}} & \mathcal{Y} \ar[d] ^{\overline{h}}\\
&S  \ar@{^{(}->}[r]_{j_{T}} & T \ar@{^{(}->}[r]_{i_{T}} & \mathcal{T} \ar[r]^{\rho}
 & \mathcal{W}\\
}
\end{array}
\leqno{(2.2.1)}
$$

\noindent dans lequel le carr\'e de gauche est cart\'esien $f, \ f_{1}$ et $\overline{f}$ sont des morphismes de $k$-sch\'emas, $\overline{h}$ et $\rho$ sont des morphismes de $\mathcal{V}$-sch\'emas formels s\'epar\'es plats de type fini, $ j_{X_{1}}, \  j_{Y}$ et $j_{T}$ sont des immersions ouvertes, $i_{Y}$ et $i_{T}$ sont des immersions ferm\'ees. Notons $X_{2}=\overline{h}^{-1}(S),\ Y_{2}=\overline{h}^{-1}(T)\ \mbox{et}\ f_{2}: X_{2}\rightarrow S,\  \overline{f}_{2}: Y_{2}\rightarrow T$\ les morphismes induits par $\overline{h}$. Soient $Y_{0}$ et $T_{0}$ les r\'eductions modulo $\pi$ de $\mathcal{Y}$ et $\mathcal{T}$ : les immersions ferm\'ees $i_{Y}$ et $i_{T}$ se factorisent respectivement via les immersions ferm\'ees $i_{2Y_{0}} : Y \hookrightarrow Y_{0}, \ i_{Y_{0}}: Y \hookrightarrow Y_{2}\hookrightarrow Y_{0}$ et $i_{T_{0}} : T \hookrightarrow T_{0}$. On d\'esigne par $\overline{h}_{X} :\  ]X[_{\mathcal{Y}} \ \longrightarrow\ ]S[_{\mathcal{T}},~ \overline{h}_{Y} :\   ]Y[_{\mathcal{Y}}\ \longrightarrow\ ]T[_{\mathcal{T}}$ les morphismes induits par $\overline{h}_{K} :\  \mathcal{Y}_{K}\ \longrightarrow\ \mathcal{T}_{K}$
[B 3, (1.1.11) (i)], [B 4, 1.1], [LS, chap 2] $j'_{Y} :\ ]X[_{\mathcal{Y}}\ \longrightarrow\ ]Y[_{\mathcal{Y}},~i'_{Y_{0}} :\ ]Y[_{\mathcal{Y}}\ \longrightarrow]Y_{0}[_{\mathcal{Y}}\  = \mathcal{Y_{K}}, \ i'_{2Y_{0}}: ]Y_{2}[\rightarrow \mathcal{Y}_{K}$ ceux induits par l'identit\'e de $\mathcal{Y}_{K}$ et $j'_{T}\ :\ ]S[_{\mathcal{T}}\ \longrightarrow\ ]T[_{\mathcal{T}},~i'_{T_{0}} :\ ]T[_{\mathcal{T}}\ \longrightarrow\ ]T_{0}[_{\mathcal{T}}\ =\ \mathcal{T}_{K}$ ceux induits par l'identit\'e de $\mathcal{T}_{K}$. Si $V$ est un voisinage strict de $]S[_{\mathcal{T}}$ dans $]T[_{\mathcal{T}}$, alors  $W := \overline{h}^{-1}_{Y}(V)$ est un voisinage strict de $ ]X_{1}[_{\mathcal{Y}}$ (donc de $ ]X[_{\mathcal{Y}}$) dans $ ]Y[_{\mathcal{Y}}$ [B 3, (1.2.7)], [LS, chap 3] et on note  $h_{V} := \overline{h}_{Y \mid W} : W \rightarrow V$.

\vskip 3mm
\noindent \textbf{Proposition (2.2.2)}.
\textit{Sous les hypoth\`eses (2.2) on a:}
\begin{enumerate}
	\item[(2.2.2.1)] \textit{Supposons que $\overline{f}^{-1}(S) = X$, alors le diagramme (2.2.1) induit un 	diagramme commutatif}

		$$
		\begin{array}{c}
		\xymatrix{
	]X[_{\mathcal{Y}} \ar@{^{(}->}[r]^{j'_{Y}} \ar[d]_{\overline{h}_{X}}  \ar@{} [dr] |{\square} & 		]Y[_{\mathcal{Y}}
	\ar@{^{(}->}[r]^{i'_{Y_{0}}} \ar[d]^{\overline{h}_{Y}} & \mathcal{Y}_{K}  \ar[d]^{\overline{h}_{K}} \\
	]S[_{\mathcal{T}} \ar@{^{(}->} [r]_{j'_{T}} & ]T[_{\mathcal{T}} \ar@{^{(}->} [r]_{i'_{T_{0}}} & 	\mathcal{T}_{K}\ ,
		} 
		\end{array}
		$$

	\textit{dans lequel le carr\'e de gauche est cart\'esien et les fl\`eches horizontales sont des 	immersions ouvertes .}\\

	\item[(2.2.2.2)] \textit{Supposons que $\overline{h}^{-1}(T) = Y $ et $\overline{f}^{-1}(S) = X$. Alors:\\
	Les deux carr\'es du diagramme pr\'ec\'edent sont cart\'esiens.\\
	Si de plus $V$ d\'ecrit un syst\`eme fondamental de voisinages stricts de $]S[_{\mathcal{T}}$ dans 	$]T[_{\mathcal{T}}$, alors $\overline{h}^{-1}_{Y}(V)$ d\'ecrit un syst\`eme fondamental de voisinages 	stricts de $]X[_{\mathcal{Y}}$ dans $]Y[_{\mathcal{Y}}$. }\\
\end{enumerate}

\vskip 3mm
\noindent \textit{D\'emonstration}.\\
 L'existence du diagramme commutatif dans (2.2.2.1) r\'esulte de (2.2.1) et [B 3, (1.1.11) (i)] ou [LS, 2.2]  via la d\'efinition des tubes [B 3, (1.1.1)], [LS, chap 2].\\
 Puisque $\overline{f}^{-1}(S) = X$ le diagramme (2.2.1) se d\'ecompose en
$$
\begin{array}{c}
$$\xymatrix{
X \ar@{^{(}->}[r]^{j_{Y}} \ar[d] \ar @{} [dr] |{\square} & Y \ar@{^{(}->}[r]^{i_{Y}} \ar[d]& \mathcal{Y} \ar@{=}[d] \\
 \overline{h}^{-1}(S) \ar[r] \ar[d]  \ar @{} [dr] |{\square} &  \overline{h}^{-1}(T) \ar[r] \ar[d] \ar @{} [dr] |{\square} & \mathcal{Y} \ar[d] ^{\overline{h}}\\
S  \ar@{^{(}->}[r]_{j_{T}} & T \ar@{^{(}->}[r]_{i_{T}} & \mathcal{T} \ar[r]^{\rho}
 & \mathcal{W}&.\\
}
\end{array}
\leqno{(2.2.2.3)}
$$
On est donc ramen\'e \`a \'etudier s\'epar\'ement le cas o\`u $\overline{h}^{-1}(T) = Y $ et $\overline{f}^{-1}(S) = X$ et celui o\`u $\overline{h}= id$: on montrera d'abord (2.2.2.2) et pour (2.2.2.1) il nous suffira de traiter le cas o\`u $h=id$.\\

 \textit{Pour (2.2.2.2)}. Les carr\'es du diagramme (2.2.2.1) sont cart\'esiens d'apr\`es la d\'efinition des tubes [B 3, (1.1.1)], [LS, chap 2] et le fait que  $\overline{h}^{-1}(T) = Y $ et $\overline{f}^{-1}(S) = X$.\\

Pour la deuxi\`eme assertion de (2.2.2.2), on va utiliser les voisinages standarts $V_{\underline{\eta}, \underline{\lambda}}= \bigcup_{n\in\mathbb{N}}V_{\eta_{n}, \lambda_{n}}$ de Berthelot [B 3, (1.2.4)], [LS, 3.4]: rappelons au passage que si $\lambda_{n}< \lambda'_{n}$ alors $V_{\eta_{n}, \lambda'_{n}}\subset V_{\eta_{n}, \lambda_{n}}$.\\

Si $V$ est un voisinage strict de $]S[_{\mathcal{T}}$ dans $]T[_{\mathcal{T}}$, alors $\overline{h}^{-1}_{Y}(V)$ est un voisinage strict de $]X[_{\mathcal{Y}}$ dans $]Y[_{\mathcal{Y}}$ [B 3, (1.2.7)], [LS, 3.1].  Soit $W$ un voisinage strict de $]X[_{\mathcal{Y}}$ dans $]Y[_{\mathcal{Y}}$ ; d'apr\`es [B 3, (1.2.2)] on se ram\`ene au cas o\`u $]Y[_{\mathcal{Y}}$ est affino¬\"{\i}de, et avec les notations de [loc. cit.] il existe $\lambda_{0} < 1$ tel que, pour $\lambda_{0} \leqslant \lambda < 1$, on ait $U_{\lambda} \subset W$. Avec les notations de [B 3, (1.2.4) (i)] si $V_{\underline{\eta},\underline{\lambda}}(S)$ parcourt un syst\`eme fondamental de voisinages stricts de $]S[_{\mathcal{T}}$ dans $]T[_{\mathcal{T}}$,  alors il existe $\eta_{n}$ et $\lambda_{n}$ assez proches de 1 tels que $(\overline{h}_{Y})^{-1}\ (V_{\eta_{n}, \lambda_{n}}(S))\ \subset\ U_{\lambda}$ car les \'equations locales de $Y$ (resp de $Z := Y\  \backslash\  X$) sont obtenues par image inverse par $\overline{h}$ des \'equations locales de $T$ (resp de $T\ \backslash\ S$) (cf. aussi [LS, 3.2]).\\

\textit{Pour (2.2.2.1)}. Puisque $]X[_{\mathcal{Y}}$ et $]Y[_{\mathcal{Y}}$ sont des ouverts de $\mathcal{Y}_{K}$ [B 3, (1.1.2)], [LS, chap 2] il en r\'esulte que $j'_{Y}$ et $i'_{Y_{0}}$ sont des immersions ouvertes ; de m\^eme pour $j'_{T}$ et $i'_{T_{0}}.$\\
  Pour montrer que le carr\'e de gauche du diagramme (2.2.2.1) est cart\'esien il nous suffit de traiter le cas o\`u $h=id$. Comme $\overline{f}^{-1}(S)=X$, les \'equations locales de $Z=Y\setminus X$ sont obtenues par image inverse par $\overline{f}$ des \'equations locales de $T\setminus S$: la d\'efinition des tubes [B 3, (1.1.1)], [LS, chap 2] fournit alors  l'\'egalit\'e $]X[_{\mathcal{Y}}=]S[_{\mathcal{Y}}\ \bigcap\ ]Y[_{\mathcal{Y}}$, d'o\`u le carr\'e cart\'esien de (2.2.2.1). \ $\square$\\

\noindent \textbf{Proposition (2.2.3)}.
\textit{Sous les hypoth\`eses de (2.2) on a:}
\begin{enumerate}
\item[(2.2.3.1)]  \textit{Supposons $\overline{h}^{-1}(T)=Y, \overline{f}^{-1}(S)=\overline{h}^{-1}(S)=X$ et $\overline{h}$ est propre. Alors $\overline{h}_{K}, \overline{h}_{Y}$ et $\overline{h}_{X}$ sont propres. Si $V$ est un voisinage strict de $]S[_{\mathcal{T}}$ dans $]T[_{\mathcal{T}}$ et $W = \overline{h}^{-1}_{Y}(V)$, alors $h_{V} = \overline{h}_{Y \mid_{W}} : W \rightarrow V$ est propre.}
\item[(2.2.3.2)] \textit{Supposons $\overline{h}$ lisse sur un voisinage de $X$ dans $\mathcal{Y}$. Alors}
	\begin{enumerate}
	\item[(i)] \textit{$\overline{h}_{X}$ est lisse et quel que soit $V$ un voisinage strict de $]S[_{\mathcal{T}}$ dans $]T[_{\mathcal{T}}$  il existe un voisinage strict $W$ de $]X[_{\mathcal{Y}}$ dans $]Y[_{\mathcal{Y}}$ tel que $\overline{h}_{K}$ induise un morphisme lisse $h_{V}: W \rightarrow V$. De plus  $\Omega^{i}_{W / V}$ est un $\mathcal{O}_{W}$-module coh\'erent et localement libre. }
	\item[(ii)] \textit{Si l'on suppose aussi que $\overline{h}^{-1}(T)=Y$, et $ \overline{h}^{-1}(S)=X$, alors il existe un voisinage strict $V$ de $]S[_{\mathcal{T}}$ dans $]T[_{\mathcal{T}}$ tel qu'en posant $W= \overline{h}^{-1}_{Y}(V)$ le morphisme $\overline{h}_{K}$ induise un morphisme lisse \\
	$h_{V}:W \rightarrow V$ .\\
	Si de plus $\overline{h}$ est propre, alors le morphisme propre $h_{V}$ est ouvert pour la topologie rigide et $h_{V}(W)$ est un ouvert admissible de $V$ et de $\mathcal{T}_{K}$}.
	\item[(iii)] \textit{Si en outre $g : \mathcal{T} \rightarrow Spf\ \mathcal{V}$ est lisse sur un voisinage de $S$ dans $\mathcal{T}$, alors on peut prendre le $V$ du (i)  lisse sur $K$ et ainsi $\Omega^{1}_{V /K}$ est localement libre de type fini sur le faisceau coh\'erent d'anneaux $\mathcal{O}_{V}$.Ê}
	\item[(iv)] \textit{Supposons $f$ surjectif, $\overline{h}^{-1}(T)=Y, \overline{h}^{-1}(S)=X$ et $\overline{h}$ est propre. Pour $V$ et $W$ comme en (ii) posons $V' = h_{V}(W)$.\\
Si $(W_{\lambda})_{\lambda}$ est un syst\`eme fondamental de voisinages stricts de $]X[_{\mathcal{Y}}$ dans $] Y [_{\mathcal{Y}}$ avec $W_{\lambda} \subset W$, alors $(h_{V}(W_{\lambda}))_{\lambda}$ est un syst\`eme fondamental de voisinages stricts de $] S [_{\mathcal{T}}$ dans $V'$.  \\
	De plus le morphisme $h_{V}$ induit un morphisme propre lisse et surjectif}
$$h_{V'} : W\ \longrightarrow\ V'\ .$$
	\end{enumerate}
\end{enumerate}

\vskip 3mm
\noindent \textit{D\'emonstration}. \\
\textit{Prouvons (2.2.3.1)}. Sous nos hypoth\`eses les deux carr\'es de (2.2.2.1) sont cart\'esiens [cf (2.2.2.2)]. D'apr\`es L¬\"utkebohmert [L¬\"u, theo 3.1] le morphisme propre $\overline{h}$ induit un morphisme propre d'espaces analytiques rigides $\overline{h}_{K} : \mathcal{Y}_{K} \rightarrow \mathcal{T}_{K} $. Comme la notion de morphisme propre est stable par changement de base en g\'eom\'etrie rigide [B-G-R, fin de 9.6.2, p 396], on en d\'eduit que $\overline{h}_{Y},\overline{h}_{X}$ et $h_{V}$  sont propres.\\

\newpage
\noindent \textit{Pour (2.2.3.2)}.
\begin{enumerate}
\item[(i)] L'ensemble $W'$ des points de $] Y [_{\mathcal{Y}}$ o\`u le morphisme $\overline{h}_{K}$ est lisse est un voisinage strict de $]X[_{\mathcal{Y}}$ dans $] Y [_{\mathcal{Y}}$ [B 3, (2.2.1)], [LS, 3.3]. Si $V$ est un voisinage strict quelconque de $]S[_{\mathcal{T}}$ dans $]T[_{\mathcal{T}}, W= \overline{h}_{Y}^{-1}(V)\cap W'$  est un voisinage strict de $]X[_{\mathcal{Y}}$ dans $]Y[_{\mathcal{Y}}$ [B 3, (1.2.7) et (1.2.10)], [LS, 3.1] et $\overline{h}_{K}$ induit donc un morphisme lisse $h_{V}: W \rightarrow V$; en particulier $\overline{h}_{X}: ]X[_{\mathcal{Y}} \rightarrow ]S[_{\mathcal{T}} \subset V$ est lisse.\\
Ê\quad  La lissit\'e de $h_{V}$ prouve que $\Omega^i_{W / V}$ est un $\mathcal{O}_{W}$-module localement libre de type fini, et comme $\mathcal{O}_{W}$ est un faisceau coh\'erent d'anneaux [prop (1.3)], il r\'esulte de [EGA $O_{I}$, (5.4.1)] que $\Omega^i_{W / V}$ est un $\mathcal{O}_{W}$-module coh\'erent.
\item[(ii)] Si $V$ d\'ecrit un syst\`eme fondamental de voisinages stricts de $]S[_{\mathcal{T}}$ dans $]T[_{\mathcal{T}}$, alors $\overline{h}_{Y}^{-1}(V)$ d\'ecrit un syst\`eme fondamental de voisinages stricts de $]X[_{\mathcal{Y}}$ dans $]Y[_{\mathcal{Y}}$ [(2.2.2.2)]: ainsi la premi\`ere assertion de (ii) r\'esulte de (i).\\
Si de plus $\overline{h}$ est propre alors le morphisme $h_{V} = \overline{h}_{Y \mid_{W}} : W \rightarrow V$ est propre et lisse; en particulier $h_{V}$ est ouvert pour la toplogie rigide: en effet on se ram\`ene, comme dans la premi\`ere \'etape de la preuve du th\'eor\`eme (1.2) (puisque $h_{V}$ est propre), au cas des espaces quasi-compacts quasi-s\'epar\'es et on applique alors [Bo-L¬\"u 2, 5.11]. Donc $V':= h_{V}(W)$ est un ouvert admissible de $V$, donc de $\mathcal{T}_{K}$ car $V$ est un ouvert admissible de $\mathcal{T}_{K}$.
\item[(iii)] Supposons en outre $g : \mathcal{T} \rightarrow Spf \mathcal{V}$ lisse sur un voisinage de $S$ dans $\mathcal{T}$ ; alors l'ensemble des points de $] T [_{\mathcal{T}}$ o\`u $g_{K}$ est lisse est un voisinage strict de $] S [_{\mathcal{T}}$ dans $] T [_{\mathcal{T}}$ [B 3, (2.2.1)], [LS, 3.3] : quitte \`a restreindre le $V$ du(i) [B 3, (1.2.10)], [LS, chap 3] on peut supposer que la restriction de $g_{K}$ \`a $V$ est lisse. Ainsi $\Omega^1_{V / K}$ sera un $\mathcal{O}_{V}$-module localement libre de type fini, donc coh\'erent sur $\mathcal{O}_{V}$  [(1.3)]. D'o\`u (iii).
\item[(iv)] Puisque $f$ est surjectif et $\overline{h}$ plat au voisinage de $X$, le morphisme
$$\overline{h}_{X} : ] X [_{\mathcal{Y}}\  \longrightarrow\  ] S [_{\mathcal{T}} $$

induit par $\overline{h}$ est surjectif [B 3, (1.1.12)] et lisse puisque c'est aussi la restriction de $h_{V}$. De plus $h_{V}(W)$ contient $] S [_{\mathcal{T}}$ par la surjectivit\'e de $\overline{h}_{X}$. Dans le diagramme commutatif
$$
\begin{array}{c}
 \xymatrix{
]X[_{\mathcal{Y}} \ar@{^{(}->}[r] \ar@{->>}[d]_{\overline{h}_{X}}   & W = h^{-1}_{V}(V') = h^{-1}_{V}(V)   \ar@{->>}[d]^{h_{V'}}  \ar@{=} [r] & W  \ar[d]^{h_{V}}\\
]S[_{\mathcal{T}} \ar@{^{(}->} [r] &V' := h_{V}(W) \ar@{^{(}->} [r] & V
} 
\end{array}
\leqno{(2.2.3.2.1)}
$$

les carr\'es sont cart\'esiens d'apr\`es (2.2.2.2)et $h_{V'}$, $\overline{h}_{X}$ sont lisses et surjectifs.\\

 \quad Soit $W_{\lambda}$ un voisinage strict de $]X[_{\mathcal{Y}}$ dans $]Y[_{\mathcal{Y}}$ avec $W_{\lambda} \subset  W$, o\`u $W$ est d\'efini ci-dessus : le morphisme ouvert $h_{V}$ [cf(ii)] envoie le recouvrement admissible  $\{ W_{\lambda} ;  W\  \backslash\  ]X[_{\mathcal{Y}} \}$ de $W$ sur le recouvrement admissible $\{ h_{V}(W_{\lambda}) ; \hfill\break h_{V}(W\  \backslash\   ]X[_{\mathcal{Y}}) \}$ de $V' = h_{V}(W)$. Par la surjectivit\'e de $h_{V'}$ et le fait que $h^{-1}_{V} (] S [_{\mathcal{T}}) =\  ] X [_{\mathcal{Y}}$ on a $h_{V}(W)\ \backslash\  ] X [_{\mathcal{Y}}) = V'\ \backslash\ ] S [_{\mathcal{T}} $ ; par suite $\ \{ h_{V}(W_{\lambda})) ; V'\ \backslash\ ] SÊ[_{\mathcal{T}} \}$ est un recouvrement admissible de $V'$, i.e. $h_{V}(W_{\lambda})$ est un voisinage strict de $] S [_{\mathcal{T}}$ dans $V'$.\\

Ê \quad Supposons maintenant que $(W_{\lambda})_{\lambda}$ d\'ecrive un syst\`eme fondamental de voisinages stricts de $]X[_{\mathcal{\mathcal{Y}}}$ dans $]Y[_{\mathcal{Y}}$ avec $W_{\lambda} \subset W$. Soit $V''$ un voisinage strict de $] S [_{\mathcal{T}}$ dans $V'$ : montrons que $h^{-1}_{V'}(V'')$ est un voisinage strict de $]X[_{\mathcal{Y}}$ dans $W$. D'abord $\{ V''; V' \backslash\ ] S [_{\mathcal{T}} \}$ est un recouvrement admissible de $V'$, donc par image inverse $\{ h^{-1}_{V'}(V'') ; h^{-1}_{V'}(V' \backslash ] S [_{\mathcal{T}}) \}$ est un recouvrement admissible de
 $h^{-1}_{V'}(V') = W$ ; en utilisant encore l'\'egalit\'e $h^{-1}_{V'} (] S [_{\mathcal{T}}) = ] X [_{\mathcal{Y}}$ on en d\'eduit que $h^{-1}_{V'}(V' \backslash\ ] S [_{\mathcal{T}}) = W \backslash\ ] X [_{\mathcal{Y}}$ : donc  $h^{-1}_{V'}(V'')$ est un voisinage strict de $] X [_{\mathcal{Y}}$ dans $W$. Ainsi il existe $\mu$ tel que $W_{\mu} \subset h^{-1}_{V'}(V'')$, d'o\`u 
 
 $$h_{V'}(W_{\mu}) \subset h_{V'}\  h^{-1}_{V'}(V'') = V''  ; $$
 
 par suite $(h_{V}(W_{\lambda}))_{\lambda}$ est bien un syst\`eme fondamental de voisinages stricts de $] S [_{\mathcal{T}}$ dans $V'$.\\
 
 \quad Si de plus $\overline{h}$ est propre alors $h_{V'}$ est de surcro\^\i t propre puisque (2.2.3.2.1) est \`a carr\'es cart\'esiens et $h_{V}$ est propre. $\square$
 
 \end{enumerate}

 \vskip5mm
 
\textbf{2.3.} Soient $S = Spec\ A_{0}$ un $k$-sch\'ema lisse et $f : X = Spec\ B_{0}\ \rightarrow\ S$ un $k$-morphisme fini. D\'esignons par $A = \mathcal{V} [t_{1},..., t_{n}]  / (f_{1},..., f_{r})$ une $\mathcal{V}$-alg\`ebre lisse relevant $A_{0}$, $P$ la fermeture projective de $Spec\ A$ dans $\mathbb{P}^{n}_{\mathcal{V}}$, $P'$ le normalis\'e de $P$ et notons $\mathcal{S} = Spf\ \hat{A}$, $\tilde{\mathcal{S}} = \hat{P}$ et $\overline{\mathcal{S}} = \widehat{P'}$  les compl\'et\'es formels de $Spec\ A$, $P$ et $P'$ comme dans le th\'eor\`eme [Et 5, th\'eo (3.1.3)]: on sait [loc. cit.] que $\tilde{\mathcal{S}}$ et $\overline{\mathcal{S}}$ sont propres sur $\mathcal{V}$, que $\overline{\mathcal{S}}$ est normal, qu'il existe une $A$-alg\`ebre finie $B$ et un carr\'e cart\'esien de $\mathcal{V}$-sch\'emas formels
$$
 \xymatrix{
Spf \hat{B} = \mathcal{X} \ar@{^{(}->}[r] \ar[d]_{h}  & \overline{\mathcal{X}} = \widehat{P''_{1}}  \ar[d]^{\overline{h}}  \\
Spf \hat{A} = \mathcal{S}  \ar@{^{(}->} [r]_{j}  &\overline{\mathcal{S}} = \widehat{P'}
} 
$$

 \noindent o\`u $P''_{1}$ est la fermeture int\'egrale de $P'$ dans Spec $B$, avec $h$ fini, $\overline{h}$ fini et $j$ une immersion ouverte. On note $\overline{f} :\overline{X}\ \rightarrow\ \overline{S}$ la r\'eduction de $\overline{h} : \overline{\mathcal{X}} \rightarrow \overline{\mathcal{S}}$ sur $k = \mathcal{V} / \mathfrak{m}$.\\
 
Rappelons [Et 5, th\'eo (3.1.3)] qu'il existe $a \in A$ et $f(t) \in A_{a}[t]$ tels que $B$ est la fermeture int\'egrale de $A$ dans $A_{a} [t]  /  (f)$. Fixons d'autre part une pr\'esentation de la $\mathcal{V}$-alg\`ebre $B$
$$B\ \simeq\ \mathcal{V} [t'_{1} ,..., t'_{n'}] / (g_{1} ,..., g_{s}).	$$

\noindent Soient $\overline{\mathcal{Y}}$ le compl\'et\'e formel de la fermeture projective $P''_{2}$ de Spec $B$ dans $\mathbb{P}^{n'}_{\mathcal{V}}$ et $\overline{Y}$ sa r\'eduction sur $k$.\\

\quad Comme $P'$ est le normalis\'e de $P$ on a un triangle commutatif
$$
\xymatrix{
&   \overline{\mathcal{S}}:=\hat{P'} \ar[d]^{v} \\
\mathcal{S} \ar@{^{(}->}[ur]^{j} \ar@{^{(}->}[r]_{\tilde{j}} & \tilde{\mathcal{S}}:=\hat{P}\\
}
$$

\noindent o\`u $v$ est fini et $\tilde{j}$ une immersion ouverte. Un syst\`eme fondamental de voisinages stricts de $\mathcal{S}_{K}$ dans $\tilde{\mathcal{S}}_{K}$ est fourni par les intersections $\tilde{V}_{\lambda}$ de $(\mbox{Spec}\ A)^{an}_{K}$ avec les boules ferm\'ees $B(0, \lambda^{+}) \subset \mathbb{A}^{n}_{K}$ pour $\lambda \rightarrow 1^+$ et $\tilde{V}_{\lambda} = \mbox{Spm}\ A_{\lambda}$, $A^{\dag}_{K} = \displaystyle \mathop{\lim}_{\rightarrow \atop{\lambda}} A_{\lambda}$  [B 3, (2.5.1)], [LS, 5.1]. Puisque $v$ est propre, et \'etale au voisinage de $\mathcal{S}$, il existe $\lambda_{0}  > 1$ tel que tout $\lambda, 1 < \lambda \leqslant \lambda_{0}$, $v$ induise un isomorphisme entre $\tilde{V}_{\lambda}$ et un voisinage strict $V_{\lambda}$ de $\mathcal{S}_{K}$ dans $\overline{S}_{K}$ [B 3, (1.3.5)], [LS, chap 3]: on identifiera $V_{\lambda}$ et $\tilde{V}_{\lambda}$ dans la suite. \\

Notons $P''_{3}$ l'adh\'erence sch\'ematique de $\mbox{Spec}\  B$ plong\'e diagonalement dans $ P''_{1}\ \times_{\mathcal{V}} P''_{2}$, $\overline{\mathcal{Z}} = \widehat{P''_{3}}$ le compl\'et\'e formel de $P''_{3}$ et $\overline{Z}$ sa r\'eduction mod $\mathfrak{m}$. On a un diagramme commutatif

$$
\xymatrix{
&   \overline{Y} \ar@{^{(}->}[r] & \overline{\mathcal{Y}}  \\
X \ar@{^{(}->}[ur] \ar@{^{(}->}[r]  \ar@{^{(}->}[rd] & \overline{Z} \ar@{^{(}->}[r] \ar[u]_{v_{2}} \ar[d]^{v_{1}}& \overline{\mathcal{Z}} \ar[u]_{u_{2}} \ar[d]^{u_{1}}\\
& \overline{X} \ar@{^{(}->}[r] & \overline{\mathcal{X}}
}
$$

\noindent o\`u les $u_{i}, v_{i}$ sont propres et les $u_{i}$ sont \'etales au voisinage de $X$. D'apr\`es [B 3, (1.3.5)], [LS, chap 3] $u_{1K}$ induit un isomorphisme  entre un voisinage strict de $] X [_{\overline{\mathcal{Z}}}$ dans $\overline{\mathcal{Z}}_{K}$ et un voisinage strict de $] X [_{\overline{\mathcal{X}}}\ \simeq\  ] X [_{\mathcal{X}}\  =\  \mathcal{X}_{K}$ dans $\overline{\mathcal{X}}_{K}$ et par suite un isomorphisme entre des syst\`emes fondamentaux de tels voisinages stricts. De m\^eme $u_{2K}$ induit un isomorphisme entre un syst\`eme fondamental de voisinages stricts de $] X [_{\overline{\mathcal{Z}}}$ dans $\overline{\mathcal{Z}}_{K}$ et un syst\`eme fondamental $(W'_{\lambda'})_{\lambda'} $ de voisinages stricts de $\mathcal{X}_{K}$ dans $\overline{\mathcal{Y}}_{K}$. Par composition il en r\'esulte pour $\lambda'\rightarrow 1^+$ un isomorphisme entre les $W'_{\lambda'}$ et un syst\`eme fondamental de voisinages stricts $(W''_{\lambda''})$ de $\mathcal{X}_{K}$ dans $\overline{\mathcal{X}}_{K}$ identifi\'es ci-apr\`es. Pour $\lambda > 1$, il existe donc $\lambda' > 1$ et des immersions ouvertes
$$\mbox{Spm}\ \hat{B}_{K} = \mathcal{X}_{K} \hookrightarrow W'_{\lambda'} \displaystyle \ \mathop{\hooklongrightarrow}^{j'_{\lambda \lambda'}}\  \overline{h}^{-1}_{K}(V_{\lambda}) =: W_{\lambda} .$$

\noindent \textbf{Proposition (2.3.1)}.
\textit{Avec les notations de 2.3 on a : }

\begin{enumerate}
\item[(1)] \textit{Si $(V_{\lambda})_{\lambda}$ est un syst\`eme fondamental de voisinages stricts de $\mathcal{S}_{K}$ dans $\overline{\mathcal{S}}_{K}$, alors $(W_{\lambda})_{\lambda} :=  (\overline{h}^{-1}_{K}(V_{\lambda}))_{\lambda} $ est un syst\`eme fondamental de voisinages stricts de $\mathcal{X}_{K}$ dans $\overline{\mathcal{X}}_{K}$. \\
Si $V_{\lambda} = \mbox{Spm}\ A_{\lambda}$, alors $W_{\lambda} = \mbox{Spm}\ B_{\lambda}$ o\`u $B_{\lambda}$ est la fermeture int\'egrale de $A_{\lambda}$ dans $\hat{B}_{K}$, et $B_{\lambda}$ est une $A_{\lambda}$-alg\`ebre finie.}

\item[(2)] \textit{Supposons de plus $f$ fini et plat (resp. fini et fid\`element plat, resp.fini \'etale, resp. fini \'etale galoisien de groupe G) et $V_{\lambda} = \mbox{Spm}\ A_{\lambda}$. Alors il existe $\lambda_{0} > 1$ tel que pour tout $\lambda, 1 < \lambda \leqslant \lambda_{0}$, et $W_{\lambda} := \overline{h}^{-1}_{K}(V_{\lambda}) $, le morphisme induit par $\overline{h}_{K}$ }
$$h_{\lambda}	:= \overline{h}_{K \mid W_{\lambda}} : W_{\lambda} \longrightarrow V_{\lambda}$$

\textit{soit fini et plat (resp. fini et fid\`element plat, resp. fini \'etale, resp. fini \'etale galoisien de groupe G), avec $V_{\lambda}$ lisse sur $K$ et $\Omega^1_{V_{\lambda}/K}$ localement libre de type fini sur le faisceau coh\'erent d'anneaux $\mathcal{O}_{V_{\lambda}}$.}

\end{enumerate}

\newpage
\noindent \textit{D\'emonstration}.  On utilise les notations du 2.3.

\begin{enumerate}

\item[(1)] On a d\'ej\`a prouv\'e la premi\`ere assertion du (1) en (2.2.2.2) : prouvons la seconde assertion du (1). On peut supposer $A, B, P'$ et $P''_{1}$ int\'egralement clos. Comme $\overline{h}_{K}$ est fini, 
$W_{\lambda} = \overline{h}^{-1}_{K}(V_{\lambda})$ est un affino\" ide  not\'e $\mbox{Spm}\ B_{\lambda}$ et $B_{\lambda}$ est une $A_{\lambda}$-alg\`ebre finie [B-GR, 9.4.4 cor 2]; on note $h_{\lambda} := \overline{h}_{K \mid W_{\lambda}} : W_{\lambda} = \mbox{Spm}\ B_{\lambda} \rightarrow V_{\lambda} = \mbox{Spm}\ A_{\lambda}$. Puisque $P''_{1}$ est la fermeture int\'egrale de $P'$ dans $Spec\ B$ , il r\'esulte de [EGA IV, 6.14.4] que $B_{\lambda}$ est la fermeture int\'egrale de $A_{\lambda}$ dans $\hat{B}_{K}$.\\
D'o\`u le (1).\\

\item[(2)] L\`a encore on peut supposer $A$ et $B$ int\'egralement clos : on traitera \`a part le cas fini \'etale galoisien.\\

\ \quad Notons $\varphi : A \rightarrow B$ le morphisme fini tel que Spec $\varphi : \mbox{Spec}\  B \rightarrow \mbox{Spec}\  A$ rel\`eve $f$ [Et 5, th\'eo (3.1.3)], et $\hat{\varphi} : \hat{A} \rightarrow \hat{B}$ (resp. $\varphi^{\dag}
 : A^{\dag} \rightarrow B^{\dag})$ le morphisme induit sur les s\'epar\'es compl\'et\'es (resp. sur les compl\'et\'es faibles). D'apr\`es [loc. cit.], $\varphi^{\dag}$ et $\hat{\varphi}$ sont finis et plats (resp. finis et fid\`element plats, resp. finis \'etales) si et seulement si $f$ l'est. Avec les notations du (2.3) on a :

 $$
 A^{\dag}_{K} = \displaystyle{\lim_{\rightarrow \atop{> \atop{\lambda \rightarrow 1}}}}\  A_{\lambda}, B^{\dag}_{K} = \displaystyle{\lim_{\rightarrow \atop{> \atop{\lambda \rightarrow 1}}}} B_{\lambda}\  \textrm{et}\   \varphi^{\dag}_{K} : A^{\dag}_{K} \rightarrow B^{\dag}_{K}
 $$
 
\noindent est la limite inductive des $\varphi_{\lambda} : A_{\lambda} \rightarrow B_{\lambda}$ [EGA IV, (8.5.2.1)] avec 
 
 $$h_{\lambda} = \mbox{Spm}\ (\varphi_{\lambda}) : W_{\lambda} = \mbox{Spm}\ B_{\lambda}\ \rightarrow V_{\lambda} = \mbox{Spm}\ A_{\lambda}. $$

 Si $f$ est fini et plat (resp. ...) alors $\varphi^{\dag}_{K}$ l'est et pour $\lambda
$ assez proche de 1, $\varphi_{\lambda}$ l'est aussi par [EGA IV, 11.2.6, 8.10.5, 17.7.8], de m\^eme  pour $h_{\lambda}$. D'o\`u le (2) hormis la cas galoisien.\\

\quad Consid\'erons \`a pr\'esent le cas galoisien.\\

\quad Puisque $f$ est galoisien il est surjectif; ainsi \\
$$h : \mathcal{X} = \mbox{Spf}\ \hat{B}\ \rightarrow\ \mathcal{S} = \mbox{Spf}\ \hat{A}$$
 est fini \'etale surjectif et galoisien de groupe $G$, d'o\`u en particulier une injection : $\hat{A} \hookrightarrow \hat{B}$. Par suite $h_{K} : \mathcal{X}_{K} \rightarrow \mathcal{S}_{K}$ est fini \'etale surjectif et galoisien de groupe $G$.\\

\quad Remarquons ensuite que puisque $B$ est la fermeture int\'egrale de $A$ dans $A_{a}[t] / (f)$ et que Spec $A^{\dag} \rightarrow \mbox{Spec}\ A $ est un morphisme normal [Et 5, prop (1.1)], il r\'esulte de [EGA IV, (6.14.4)] que $B^{\dag} = B \otimes_{A} A^{\dag} $ est la fermeture int\'egrale de $A^{\dag}$ dans $(B^{\dag})_{a} = (A^{\dag})_{a}\  [t] / (f) $ : par suite $B^{\dag}_{K}$ est la fermeture int\'egrale de $A^{\dag}_{K}$ dans $(A^{\dag})_{a,K}\  [t] / (f)$. L'anneau $A^{\dag}$ est r\'eduit par [Et 5, prop (1.7)], car $A$ est r\'eduit, et $A^{\dag} \rightarrow B^{\dag}$ est fini \'etale car $\hat{A} \rightarrow \hat{B}$ l'est : donc $B^{\dag}$ est r\'eduit car $A^{\dag}$ est r\'eduit [Et 5, lemme (1.6)]. Ainsi $B^{\dag}$ est int\'egralement ferm\'e dans $\hat{B}$ [Et 5, th\'eo (2.2)(2)ii]; d'o\`u $B^{\dag}_{K}$ est la fermeture int\'egrale de $A^{\dag}_{K}$ dans $\hat{B}_{K}$.\\
 
 \quad On a vu ci-dessus que, pour $\lambda$ suffisamment proche de 1, 
 
 $$h_{\lambda} : W_{\lambda} = \mbox{Spm}\ B_{\lambda} \rightarrow V_{\lambda} = \mbox{Spm}\ A_{\lambda} $$
 
 est fini \'etale. Compte tenu du diagramme commutatif

  $$
 \xymatrix{
 \hat{B}_{K} &    & B_{\lambda}  \ar@{_{(}->} [ll]  \\
 \hat{A}_{K} \ar@{^{(}->} [u] & A^{\dag}_{K} \ar@{_{(}->} [l] & A_{\lambda} \ar@{_{(}->} [l] \ar[u]
  }
 $$

 \noindent la fl\`eche $A_{\lambda} \rightarrow B_{\lambda} $ induite par $h_{\lambda}$ est injective, donc $h_{\lambda}$ est surjectif. Choisissons un ensemble fini $\{x_{i} \}$ de g\'en\'erateurs de $B_{\lambda}$ sur $A_{\lambda}$.  Comme chaque $x_{i}$ est entier sur $A_{\lambda}$, les \'el\'ements $g_{\hat{B}_{K}}(x_{i}) \in \hat{B}_{K}$, pour $g$ d\'ecrivant $G$ et $g_{\hat{B}_{K}} : \hat{B}_{K} \rightarrow \hat{B}_{K}$ induit par $g$, sont aussi entiers sur $A_{\lambda} \subset A^{\dag}_{K}$, donc a fortiori sur $A^{\dag}_{K} = \displaystyle \mathop{\lim}_{\rightarrow \atop{\mu}} A_{\mu}$. Or on a vu que $B^{\dag}_{K}$ est la fermeture int\'egrale de $A^{\dag}_{K}$ dans $\hat{B}_{K}$ et que $B^{\dag}_{K}$ est int\'egralement ferm\'e dans $\hat{B}_{K}$ : il existe donc $\lambda'$, $1 < \lambda' \leqslant \lambda$ tel que pour tout $i$ et tout $g \in G$ on ait $g_{\hat{B}_{K}}(x_{i}) \in B_{\lambda'}$. Ainsi l'action de $G$ s'\'etend de $\mathcal{X}_{K}$ \`a $W_{\lambda'} = \mbox{Spm}\ B_{\lambda'}$ : en effet $g \in G$ d\'efinit un morphisme $g_{\lambda \lambda'}: W_{\lambda'} \rightarrow W_{\lambda}$ s'ins\'erant dans le diagramme commutatif \`a carr\'e cart\'esien

$$
\xymatrix{W_{\lambda'} \ar@/^/[rrd]^{g_{\lambda \lambda'}} \ar@/_/[rdd]_{h_{\lambda'}} \ar@{.>}[rd]_{g_{\lambda'}} \\
& W_{\lambda'} \ar[d]^{h_{\lambda'}} \ar@{^{(}->}[r] & W_{\lambda}  \ar[d]^{h_{\lambda}}\\
& V_{\lambda'} \ar@{^{(}->}[r]_{\alpha_{\lambda \lambda'}} & \ V_{\lambda}\ ;}
$$

  \noindent d'o\`u la factorisation de $g_{\lambda \lambda'}$ par $W_{\lambda'}$.\\
  
  Montrons que le morphisme fini \'etale surjectif
    $$h_{\lambda'} : W_{\lambda'} \rightarrow V_{\lambda'}  $$
  
  est galoisien de groupe $G$. On a
    $$(B_{\lambda'})^G \subset (\hat{B}_{K})^G = \hat{A}_{K} ;	$$
  
  d'o\`u
    $$A_{\lambda'} \subset (B_{\lambda'})^G \subset B_{\lambda'} \cap \hat{A}_{K} $$
  
  et on dispose d'un carr\'e commutatif 
$$
\xymatrix{
A_{\lambda'} \ar@{^{(}->}[r] \ar[d] & \hat{A}_{K} \ar[d]\\
B_{\lambda'} \ar@{^{(}->}[r]& \hat{B}_{K}
}
$$
  
  \noindent avec $A_{\lambda'} \rightarrow B_{\lambda'}$ fid\`element plat et $\hat{B}_{K} = \hat{A}_{K} \otimes_{A_{\lambda'}} B_{\lambda'}$: d'apr\`es [Et  2, prop 2] on en d\'eduit 
  
  $$A_{\lambda'} = B_{\lambda' } \cap \hat{A}_{K} = (B_{\lambda'})^G , $$
  
  d'o\`u la proposition (2.3.1). $\square$
 \end{enumerate}

 \vskip 3mm
 \section*{3. Images directes d'isocristaux}
 
\subsection*{3.1 Sections surconvergentes}
 
 On suppose donn\'e un diagramme commutatif tel que (2.2.1). Pour un voisinage strict $W$ (resp. un couple de voisinages stricts $W'\subset W$) de $]X[_{\mathcal{Y}}$ dans $]Y[_{\mathcal{Y}}$ on note $\alpha_{W}$ (resp.$\alpha_{W W'})$ l'immersion ouverte de $W$ dans $]Y[_{\mathcal{Y}}$(resp. de $W'$ dans $W$). Si $\mathcal{A}$ est un faisceau d'anneaux sur $W$ et $E$ un $\mathcal{A}$-module, on pose [B 3,(2.1.1.1)], [LS, chap 5]:\\
 
 \noindent (3.1.1) $\qquad\qquad j^{\dag}_{W}E:=\displaystyle\mathop{\mbox{lim}}_{\longrightarrow\atop{W'\subset W}} \alpha^{}_{W W'^{\ast}} \alpha^{\ast}_{W W'} E$ ,\\
 \noindent la limite \'etant prise sur les voisinages $W'\subset W$.\\
 De m\^eme [B 3,(2.1.1.3)], [LS, chap 5] :\\
 
 \noindent (3.1.2) $\qquad\qquad j^{\dag}_{Y} E:= \alpha^{}_{W ^{\ast}} j^{\dag}_{W} E$ .\\
 
 \noindent(3.1.3) Si $V$ est un voisinage strict de $]S[_{\mathcal{T}}$ dans $]T[_{\mathcal{T}}$, alors $W = \overline{h}^{\ -1}_{K} (V)\ \cap\ ]Y[_{\mathcal{Y}}\ =\overline{h}^{\ -1}_{Y} (V)\ $ est un voisinage strict de $]X[_{\mathcal{Y}}$ dans $]Y[_{\mathcal{Y}}$ [B 3, (1.2.7)], [LS, chap 3] et on note $h_{V}$ la restriction de $\overline{h}^{}_{Y} $ \`a $W$, et $R^{i} \overline{h}_{K^{\ast}} j^{\dag}_{Y} E:= R^{i}Ê\overline{h}_{Y^{\ast}} j^{\dag}_{Y} E.$\\
 
 \vskip 3mm
 \noindent \textbf{Proposition (3.1.4)}.
\textit {Avec les hypoth\`eses et notations de (3.1.3) supposons que}  $\overline{h}_{Y} :\  ]Y[_{\mathcal{Y}}{\longrightarrow}]T[_{\mathcal{T}} $  \textit{soit quasi-compact et quasi-s\'epar\'e; soit $E$ un faisceau ab\'elien sur $W$.}
 \begin{enumerate}
 \item[(a)] \textit{Supposons que $\overline{h}^{\ -1}(T) = Y$ et $\overline{h}^{\ -1}(S) = X$; alors, pour tout entier $ i \geqslant 0 $, on a des isomorphismes canoniques}
   \begin{enumerate}
   \item [(3.1.4.1)] $ \qquad\qquad R^{i}h_{V^{\ast}}(j^{\dag}_{W}E ) \overset{\sim}{\longrightarrow} j^{\dag}_{V}R^{i}h_{V^{\ast}} (E ).  $
   \item [(3.1.4.2)] $ \qquad\qquad R^{i}\overline{h}_{K^{\ast}}(j^{\dag}_{Y}E ) \overset{\sim}{\longrightarrow} j^{\dag}_{T}R^{i}h_{V^{\ast}} (E ).  $\\
   
   \textit{Si de plus $\overline{h}$ est une immersion ferm\'ee, alors}
  \item [(3.1.4.3)] $ \qquad\qquad R^{i}h_{V^{\ast}}(j^{\dag}_{W}E )\ =\ 0$ \textit{pour $ i \geqslant 1 $}\\
  
   \textit{et le morphisme canonique}\\
   \item [(3.1.4.4)] $ \qquad\qquad \overline{h}_{K}^{\ast}\overline{h}_{K^{\ast}}j^{\dag}_{Y}E  \overset{\sim}{\longrightarrow} j^{\dag}_{Y} E $\\
      \item[]\textit{est un isomorphisme.}
    \end{enumerate}
  \item[(b)] \textit{Si l'on ne suppose plus que $\overline{h}^{\ -1}(T) = Y$ et $\overline{h}^{\ -1}(S) = X$, alors, pour tout entier $ i \geqslant 0 $, on a des isomorphismes canoniques}
   \begin{enumerate}
   \item [(3.1.4.5)] $ \qquad\qquad R^{i}h_{V^{\ast}}(j^{\dag}_{W}E ) \overset{\sim}{\longrightarrow} j^{\dag}_{V}R^{i}h_{V^{\ast}} (j^{\dag}_{W}E ).$
   \item [(3.1.4.6)] $ \qquad\qquad R^{i}\overline{h}_{K^{\ast}}(j^{\dag}_{Y}E ) \overset{\sim}{\longrightarrow} j^{\dag}_{T}R^{i}h_{V^{\ast}} (j^{\dag}_{W}E ) .$
    \end{enumerate}

 \end{enumerate}
  \vskip15mm
 \noindent\textit{D\'emonstration}
 \begin{enumerate}
   \item[(a)] Les deux foncteurs
 $$\mathcal{F}: E \longmapsto h_{V^{\ast}}j^{\dag}_{W}E$$ et $$\mathcal{G}: E \longmapsto j^{\dag}_{V}h_{V^{\ast}}E$$ de la cat\'egorie $\mathcal{C}$ des faisceaux ab\'eliens sur $W$ dans la cat\'egorie  des faisceaux ab\'eliens sur $V$ sont exacts \`a gauche [B 3,(2.1.3)(iii)], [LS, 5.3]. Comme la cat\'egorie ab\'elienne $\mathcal{C}$ admet suffisamment d'injectifs, les foncteurs d\'eriv\'es droits $R^{i}\mathcal{F}$ et  $R^{i}\mathcal{G}$ existent et  $ (R^{i} \mathcal{F})_{i}$ , $(R^{i} \mathcal{G})_{i}$ sont des $ \delta$-foncteurs universels:   puisque $j^{\dag}_{W}$ et $j^{\dag}_{V}$ sont exacts [loc. cit.],  et que $\alpha_{W}:  W \hooklongrightarrow ]Y[_{\mathcal{Y}}$  (resp. $\alpha_{V} : V \hooklongrightarrow ]T[_{\mathcal{T}})$   est exact sur la cat\'egorie des  $j^{\dag}_{W}\mathbb{Z}$-modules (resp. des  $j^{\dag}_{V}\mathbb{Z}$-modules) [B 3, d\'em. de (2.1.3)], on est ramen\'e pour le (a) \`a prouver que $\mathcal{F}= \mathcal{G}$.\\
 
 Comme $\overline{h}_{Y}$ est quasi-compact et quasi-s\'epar\'e, il en est de m\^eme par changement de base pour $h_{V}$, donc $h_{V^{\ast}}$ commute aux limites inductives filtrantes: de plus les hypoth\`eses entra\^inent que si $V'$ d\'ecrit un syst\`eme fondamental de voisinages stricts de $]S[_{\mathcal{T}}$ dans $]T[_{\mathcal{T}}$, alors $h^{-1}_{V}(V')= W'$ d\'ecrit un syst\`eme fondamental de voisinages stricts de $]X[_{\mathcal{Y}}$ dans $]Y[_{\mathcal{Y}}$ [(2.2.2.2)]; d'o\`u
   \begin{eqnarray*}
   h_{V^{\ast}}(j^{\dag}_{W}E )&=& \displaystyle\mathop{\mbox{lim}}_{\longrightarrow\atop{W'\subset W}} h_{V^{\ast}} \alpha^{}_{WW'^{\ast}} \alpha^{-1}_{W W'} E\\
     &=&\displaystyle\mathop{\mbox{lim}}_{\longrightarrow\atop{V' \subset V}} \alpha^{}_{VV'^{\ast}} h_{V'^{\ast}} \alpha^{-1}_{W W'} E\\
     &=& \displaystyle\mathop{\mbox{lim}}_{\longrightarrow\atop{V'\subset V}} \alpha^{}_{VV'^{\ast}}\alpha^{-1}_{VV'} h_{V^{\ast}}E\\
     &=&  j^{\dag}_{V}h_{V^{\ast}} (E ), 
      \end{eqnarray*}
ce qui prouve (3.1.4.1) et (3.1.4.2).\\

Si $ \overline{h}$ est une immersion ferm\'ee, alors $ \overline {h}_{K }$ en est une aussi [B 3, (0.2.4)(iv)] , ou [Bo-L\"u 1, \S4] et [B-G-R, 9.5.3 prop 2 et 7.1.4 def 3], de m\^eme que $h_{V}$ par changement de base: en particulier $h_{V}$ est quasi-compact et quasi-s\'epar\'e. Ainsi (3.1.4.3) r\'esulte de (3.1.4.1) et (1.3.3). Pour (3.1.4.4) on utilise la suite d'isomorphismes
 \begin{eqnarray*}
   \overline{h}_{K}^{\ast} \overline{h}_{K^{\ast}}j^{\dag}_{Y}E & \overset{\sim}{\longrightarrow}&\overline{h}_{K}^{\ast}j^{\dag}_{T} h_{V^{\ast}}E \qquad\qquad (3.1.4.2)    \\
    & \overset{\sim}{\longleftarrow} & j_{Y}^{\dag} h^{\ast}_{V} h_{V^{\ast}}(E)\  \quad\qquad [\mbox{B}\ 3, (2.1.4.8)], [\mbox{LS}, 5.3]
     \\
     & \overset{\sim}{\longrightarrow} & j_{Y}^{\dag}(E) \qquad\qquad\qquad (1.3.3). 
   \end{eqnarray*}
     \item[(b)]L'ouvert W = $\overline{h}_{Y}^{\ -1}(V)$ est un voisinage strict de $]X_{1}[_{\mathcal{Y}}$ (donc de $]X[_{\mathcal{Y}}$) dans $]Y[_{\mathcal{Y}}$; la d\'efinition de $j^{\dag}_{W}(E)$ fait intervenir une limite inductive sur les voisinages stricts $W^{'}$ de $]X[_{\mathcal{Y}}$ dans $]Y[_{\mathcal{Y}}$: si cette fois la limite inductive est prise sur les voisinages stricts $W^{'}_{1}$de $]X_{1}[_{\mathcal{Y}}$ dans $]Y[_{\mathcal{Y}}$ nous noterons $j^{\dag}_{W_{1}}(E)$ le r\'esultat, et on a [B 3, (2.1.7)], [LS, chap 5]
        $$j^{\dag}_{W}(E) = j^{\dag}_{W} \circ j^{\dag}_{W_{1}}(E) = j^{\dag}_{W_{1}}\circ j^{\dag}_{W}(E).$$
 D'o\`u, en appliquant (3.1.4.1):
  \begin{eqnarray*}
   R^{i}h_{V^{\ast}}(j^{\dag}_{W}E )&=& R^{i}h_{V^{\ast}}j^{\dag}_{W_{1}} (j^{\dag}_{W}E )\\
     &=&j^{\dag}_{V} R^{i}h_{V^{\ast}}(j^{\dag}_{W}E);
     \end{eqnarray*} 
 de m\^eme pour (3.1.4.6). $\square$
 \end{enumerate}

\subsection*{3.2  D\'efinition des images directes}
   
  \noindent \textbf{(3.2.1)} On suppose fix\'e un diagramme commutatif tel que (2.2.1) et on fait l'hypoth\`ese suppl\'ementaire que $\overline{h}$ (resp $\rho$) est lisse sur un voisinage de $X$ dans $\mathcal{Y}$ (resp de $S$ dans $\mathcal{T}$)
$$\xymatrix{
X \ar@{^{(}->}[r]^{j_{Y}} \ar[d]_{f} & Y \ar@{^{(}->}[r]^{i_{Y}} \ar[d]^{\overline{f}} & \mathcal{Y} \ar[d] ^{\overline{h}}\\
S  \ar@{^{(}->}[r]_{j_{T}} & T \ar@{^{(}->}[r]_{i_{T}} & \mathcal{T} \ar[r]^{\rho} & \mathcal{W}.\\
}
$$

Soient $E\in Isoc^{\dag} ((X,Y)/ \mathcal{W})$, $W$ un voisinage strict de $]X[_{\mathcal{Y}}$ dans $]Y[_{\mathcal{Y}}$ et $E_{W}$ un $\mathcal{O}_{W}$-module coh\'erent muni d'une connexion int\'egrable tel que $j_{Y}^{\dag}E_{W}=: E_{\mathcal{Y}}$ soit une r\'ealisation de $E$ [B 3, (2.3.2)], [LS, 7.1]. Lorsque $Y$ est propre sur $\mathcal{W}$ on notera $Isoc^{\dag} ((X,Y)/ \mathcal{W})=Isoc^{\dag} (X/ \mathcal{W})$ [B 3, (2.3.6)]: si de plus $\mathcal{W}= Spf\mathcal{V}$ on notera $Isoc^{\dag} (X/ Spf\mathcal{V})=: Isoc^{\dag} (X/ K)$.\\
 
Pour $E\in Isoc^{\dag} ((X,Y)/ \mathcal{W})$ Berthelot a d\'efini dans [B 5,(3.1.11)] les images directes en cohomologie rigide (cf. aussi [LS, (7.4)] et [C-T, 10]) par la formule
 \begin{eqnarray*}
 (3.2.1.1)\qquad\qquad\mathbb{R}\overline{f}_{rig^{\ast}}((X,Y)/ \mathcal{T};E)&:=& \mathbb{R}\overline{h}_{K^{\ast}}(j^{\dag}_{Y}E_{W}\otimes_{\mathcal{O}_{]Y[_{\mathcal{Y}}}}\Omega^{^{\bullet}}_{{]Y[_{\mathcal{Y}}}/ \mathcal{T}_{K}} ),\\
     &:=& \mathbb{R}\overline{h}_{Y^{\ast}}(j^{\dag}_{Y}E_{W}\otimes_{\mathcal{O}_{]Y[_{\mathcal{Y}}}}\Omega^{^{\bullet}}_{{]Y[_{\mathcal{Y}}}/ ]T[_{\mathcal{T}}} );
  \end{eqnarray*}
  \noindent et la cohomologie de ces complexes est ind\'ependante du $\mathcal{Y}$ choisi [B 5, (3.1.2)],  [LS, 7.4.2].\\
 Lorsque $X = Y$, alors $\overline{f}:X \rightarrow T$ et on obtient la cohomologie convergente:\\
 
  \noindent(3.2.1.2)$\qquad\qquad\mathbb{R}\overline{f}_{conv^{\ast}}((X,Y)/ \mathcal{T};E):= \mathbb{R}\overline{f}_{rig^{\ast}}((X,X)/ \mathcal{T};E).$\\
  
  \noindent\textbf{(3.2.2)} Sous les hypoth\`eses (3.2.1) supposons de plus $\overline{h}$ propre: ainsi $Y$ est une compactification $\overline{X}=Y$ de $X$ au-dessus de $\mathcal{T}$. Berthelot d\'efinit alors $\mathbb{R}f_{rig^{\ast}}(X/ \mathcal{T};E)$ par la formule [B 5, (3.2.3)] ( cf. aussi [LS, 8.2])\\
   
\noindent$ (3.2.2.1)\qquad\qquad\mathbb{R}f_{rig^{\ast}}(X/ \mathcal{T};E):=\mathbb{R}\overline{f}_{rig^{\ast}}((X,\overline{X})/ \mathcal{T};E);$\\

  \noindent et la cohomologie de ce complexe est ind\'ependante du $\overline{X}$ choisi [B 5,(3.2.2)] [LS, 8.2.1] [C-T, 10.5.3].\\

     \subsection*{3.3 Changement de base}
  
  \noindent\textbf{(3.3.1)} Avec les notations de (2.2) on consid\`ere un parall\'el\'epip\`ede commutatif\\  
  
  $$
\begin{array}{c}
\xymatrix{
& X \ar@{.>}[dd]^(.3){f} |\hole \  \ar@{^{(}->}[rr]^{j_{Y}} & & Y\ar@{.>}[dd]^(.3){\overline{f}} |\hole \ \ar@{^{(}->}[rr]^{i_{Y}} & & \mathcal{Y} \ar[dd]^{\overline{h}} \\
X^{'} \ \ar@{^{(}->}[rr]^(.8){j_{Y'}} \ar[dd]_{f^{'}} \ar[ur]^{\varphi'} & & Y' \ \ar@{^{(}->}[rr]^(.8){i_{Y'}} \ar[dd]_(.7){\overline{f}'} \ar[ur]^{{\overline{\varphi}}'} & & \mathcal{Y}' \ar[dd]_(.7){\overline{h}'} \ar[ur]^{\overline{g}'} \\
& S\ \ar@{^{(}.>}[rr]^(.7){j_{T}} |\hole & &T\ \ar@{^{(}.>}[rr]^(.8){i_{T}}  |\hole & & \mathcal{T}\ar[rr]^{\rho} && \mathcal{W}\\
S^{'} \  \ar@{^{(}->}[rr]_{j_{T'}} \ar@{.>}[ur]^{\varphi}& & T' \ \ar@{^{(}->}[rr]_{i_{T'}} \ar@{.>}[ur]_{\overline{\varphi}} & & \mathcal{T'} \ar[ur]_{\overline{g}}\ar[rr]_{\rho '}&&\mathcal{W'}\ar[ur]_{\theta} 
}
\end{array}
\leqno{(3.3.1.1)}
  $$
  \noindent dans lequel le cube de gauche est form\'e de $k$-sch\'emas s\'epar\'es de type fini, les morphismes $\overline{g}$, $\overline{h}$, $\overline{g}'$, $\overline{h}'$, $\rho $, $\rho '  $, $\theta$ sont des morphismes de $\mathcal{V}$-sch\'emas formels ($\mathcal{V}$-sch\'emas formels que l'on supposera s\'epar\'es, plats et de type fini), les $j$ (resp. les $i$) sont des immersions ouvertes (resp. ferm\'ees). On suppose $\theta$ lisse et $X'=X\times_{S}S'$, $Y'=Y\times_{T}T'$, et $\mathcal{Y}'=\mathcal{Y}\times_{\mathcal{T}}\mathcal{T'}$.\\
  On v\'erifie alors facilement que l'on a\\
  
\noindent (3.3.1.2)
 $ \qquad\qquad  ]X'[_{\mathcal{Y'}}=]X[_{\mathcal{Y}}\times_{]S[_{\mathcal{T}}}]S'[_{\mathcal{T'}},  
     \qquad\qquad  ]Y'[_{\mathcal{Y'}}=]Y[_{\mathcal{Y}}\times_{]T[_{\mathcal{T}}}]T'[_{\mathcal{T'}}.$\\
     
       Soient $V$ un voisinage strict de $]S[_{T}$ dans $]T[_{T}$ et $V'$ un voisinage strict de $]S'[_{T'}$ dans $]T'[_{T'}$ tel que $V'\subset{\overline{g}_{K}^{\ -1}}(V)\cap]T'[_{\mathcal{T'}}$; alors $W:={\overline{h}_{K}^{\ -1}}(V)\cap]Y[_{\mathcal{Y}}$ est un voisinage strict de $]X[_{Y}$ dans $]Y[_{Y}$ et $W':={\overline{h}_{K}^{ \prime-1}}(V')\cap]Y'[_{\mathcal{Y'}}$ est un voisinage strict de $]X'[_{Y'}$ dans $]Y'[_{Y'}$ tel que $W'\subset{\overline{g}_{K}^{  \prime -1}}(W)\cap]Y'[_{\mathcal{Y'}}$.\\
 Notons \\
     
          $h_{V}: W\rightarrow V, \ g_{V}: V'\rightarrow V, \ h'_{V'}: W'\rightarrow V', \ \mbox{et} \ g'_{W}: W'\rightarrow W $\\
          
\noindent les morphismes induits respectivement par $\overline{h}_{K},\ \overline{g}_{K},\ \overline{h}'_{K},\  \mbox{et}\  \overline {g}'_{K}.$\\
  Ainsi on dispose d'un cube commutatif
  
  $$
\begin{array}{c}
\xymatrix{
& W \ \ar@{.>} [dd]^(.3){h_{V}} |\hole  \ar@{^{(}->} [rr] & & {\mathcal{Y}_{K}} \ar [dd]^{\overline{h}_{K}} \\
W' \ \ar@{^{(}->} [rr] \ar [dd]_{h'_{V'}} \ar[ur]^{g'_{W}} & & {\mathcal{Y}'_{K}} \ar[dd]_(.7){\overline{h}'_{K}} \ar [ur]^{{\overline{g}'_{K}}} \\
& V \ \ar@{^{(}.>}[rr] |\hole & &{\mathcal{T}_{K}}\\
V'  \ \ar@{^{(}->}[rr] \ar@{.>}[ur]^{g_{V}}& & \mathcal{T}'_{K}  \ar[ur]_{\overline{g}_{K}} 
}
\end{array}
\leqno{(3.3.1.3)}
  $$
  
  \noindent dans lequel  $W'=W\times_{V}V'$ et $\mathcal{Y}'_{K}=\mathcal{Y}_{K}\times_{\mathcal{T}_{K}}\mathcal{T'}_{K}$. \\
  
  \noindent\textbf{Lemme (3.3.1.4)}.
   \textit{Avec les hypoth\`eses et notations de (3.3.1) supposons de plus que $\overline{g}$ soit plat (resp. que $\overline{g}$ soit lisse sur un voisinage de $S'$ dans $\mathcal{T'}$). Alors il existe un voisinage strict $V'$ de $]S'[_{\mathcal{T'}}$ dans $]T'[_{\mathcal{T'}}$ tel que $g_{V}$ soit plat (resp. que $g_{V}$ soit lisse).}\\
   
   \noindent\textit{D\'emonstration}. Dans le cas plat  c'est clair puisque $\overline{g}_{K}$est plat; dans le cas lisse, c'est [B 3, (2.2.1)], [LS, chap 3]. $\square$\\
   
   \noindent\textbf{D\'efinition(3.3.1.5)}
   \textit{Soit} $\mathcal{E}$   \textit{un} $j_{T}^{\dag}\mathcal{O}_{V}$-\textit{module; son image inverse surconvergente est d\'efinie par la formule}
   
   $$(\varphi, \overline{\varphi}, \overline{g})^{\dag}(\mathcal{E}):= j^{\dag}_{T'}(\overline{g}_{K}^{\ast}\mathcal{E});$$
   \textit{lorsque} $\mathcal{E}$ \textit{est une r\'ealisation d'un isocristal E $\in Isoc^{\dag}((S,T)/\mathcal{W})$ on \'ecrira aussi}
   $$(\varphi, \overline{\varphi})^{\ast}(\mathcal{E})=(\varphi, \overline{\varphi}, \overline{g})^{\dag}(\mathcal{E})= j^{\dag}_{T'}(\overline{g}_{K}^{\ast}\mathcal{E}).$$\\

   \noindent\textbf{Th\'eor\`eme (3.3.2)}
   \textit{Sous les hypoth\`eses (3.3.1) supposons que $\overline{h}^{\ -1}(T)=Y, \overline{h}^{\ -1}(S)=X$ et $h_{V}$ propre (cette derni\`ere hypoth\`ese est v\'erifi\'ee si $\overline h$ est propre).}
\begin{enumerate}
      \item[(3.3.2.1)] \textit{
      Soit $E_{W}$ un $\mathcal{O}_{W}$-module coh\'erent. Alors, pour tout entier $i\geqslant 0$, on a:
      				    }
         \begin{enumerate}
         \item[(1)]$R^{i}\overline{h}_{K^{\ast}}(j^{\dag}_{Y}E_{W} )$ \textit{est un $j_{T}^{\dag}\mathcal{O}_{]T[_{\mathcal{T}}}$-module coh\'erent et on a un isomorphisme}
         
         $ \qquad\qquad R^{i}\overline{h}_{K^{\ast}}(j^{\dag}_{Y}E_{W} ) \overset{\sim}{\longrightarrow} j^{\dag}_{T}R^{i}h_{V^{\ast}} (E_{W} ).  $
         
         \item[(2)] \textit{Supposons de plus $g_{V}$ plat; alors}
   		\begin{enumerate}
   			\item[(i)]\textit{On a des isomorphismes de changement de base au sens surconvergent}				\begin{eqnarray*}
				(\varphi, \overline{\varphi}, \overline{g})^{\dag}(R^{i}\overline{h}_{K^{\ast}}j^{\dag}_{Y}E_{W} )&\simeq & R^{i}\overline{h}'_{K^{\ast}}(\overline{g}'^{\ast}_{K}j^{\dag}_{Y}E_{W} )\\
		&\simeq &	 R^{i}\overline{h}'_{K^{\ast}}(j^{\dag}_{Y'}g'^{\ast}_{W}E_{W} )\\
		&\simeq & j^{\dag}_{T'}R^{i}{h}'_{V'^{\ast}}(g'^{\ast}_{W}E_{W} ).
				\end{eqnarray*}   			
			\item[(ii)]\textit{ Si de plus $\overline{\varphi}^{\ -1}(S)= S'$, les isomorphismes pr\'ec\'edents deviennent}
				\begin{eqnarray*}
				\overline{g}_{K}^{\ast}R^{i}\overline{h}_{K^{\ast}}j^{\dag}_{Y}E_{W} &\simeq & R^{i}\overline{h}'_{K^{\ast}}(j^{\dag}_{Y'}g'^{\ast}_{W}E_{W} )\\		
						&\simeq & j^{\dag}_{T'}R^{i}{h}'_{V'^{\ast}}(g'^{\ast}_{W}E_{W} ),
				\end{eqnarray*} 
				\textit{et l'on a des isomorphismes
				\begin{eqnarray*}
				\overline{g}_{V}^{\ast}R^{i}h_{V^{\ast}}j^{\dag}_{W}E_{W} &\simeq & R^{i}h'_{V'^{\ast}}(g'^{\ast}_{W}j_{W}^{\dag}E_{W} )\\		
						&\simeq & R^{i}h'_{V'^{\ast}}j_{W'}^{\dag}g'^{\ast}_{W}E_{W} .
				\end{eqnarray*} 
				}
   		  \end{enumerate}
   	     \end{enumerate}
      \item[(3.3.2.2)]
      \textit{
      Soient $E_{W}^{^{\bullet}}$ un complexe born\'e de $\mathcal{O}_{W}$-module coh\'erents et 
      $$  E_{W'}^{^{\bullet}}=E_{W}^{^{\bullet}} \otimes_{\mathcal{O}_{V}} \mathcal{O}_{V'} = g_{W}^{\prime \ast}E_{W}^{^{\bullet}}.$$ 
      Alors, pour tout entier $i\geqslant 0$, on a:
      				    }
         \begin{enumerate}
         \item[(1)]$R^{i}\overline{h}_{K^{\ast}}(j^{\dag}_{Y}E_{W}^{^{\bullet}} )$ \textit{est un $j_{T}^{\dag}\mathcal{O}_{]T[_{\mathcal{T}}}$-module coh\'erent et on a un isomorphisme}
         
         $ \qquad\qquad R^{i}\overline{h}_{K^{\ast}}(j^{\dag}_{Y}E_{W}^{^{\bullet}} ) \overset{\sim}{\longrightarrow} j^{\dag}_{T}R^{i}h_{V^{\ast}} E_{W}^{^{\bullet}}.  $   
               
         \item[(2)] \textit {Supposons de plus $g_{V}$ plat; alors}
   		\begin{enumerate}
   			\item[(i)]\textit{On a des isomorphismes de changement de base au sens surconvergent}
				\begin{eqnarray*}
				(\varphi, \overline{\varphi}, \overline{g})^{\dag}(R^{i}\overline{h}_{K^{\ast}}j^{\dag}_{Y}E_{W}^{^{\bullet}} )&\simeq &  R^{i}\overline{h}'_{K^{\ast}}(j^{\dag}_{Y'}g'^{\ast}_{W}E_{W}^{^{\bullet}} )\\
				&\simeq & j^{\dag}_{T'}R^{i}{h}'_{V'^{\ast}}(E_{W'}^{^{\bullet}} ).
				\end{eqnarray*}   			
			\item[(ii)]\textit{ Si de plus $\overline{\varphi}^{\ -1}(S)= S'$, les isomorphismes pr\'ec\'edents deviennent}
				\begin{eqnarray*}
				\overline{g}_{K}^{\ast}R^{i}\overline{h}_{K^{\ast}}j^{\dag}_{Y}E_{W}^{^{\bullet}} &\simeq & R^{i}\overline{h}'_{K^{\ast}}(j^{\dag}_{Y'}E_{W'}^{^{\bullet}})\\		
						&\simeq & j^{\dag}_{T'}R^{i}h'_{V'^{\ast}}(E_{W'}^{^{\bullet}} ),
				\end{eqnarray*} 
				\textit{et l'on a des isomorphismes				
				\begin{eqnarray*}
				\overline{g}_{V}^{\ast}R^{i}h_{V^{\ast}}j^{\dag}_{W}E_{W}^{^{\bullet}} &\simeq & R^{i}h'_{V'^{\ast}}g'^{\ast}_{W}j_{W}^{\dag}E_{W}^{^{\bullet}} \\		
						&\simeq & R^{i}h'_{V'^{\ast}}j_{W'}^{\dag}g'^{\ast}_{W}E_{W}^{^{\bullet}} .
				\end{eqnarray*} 
				}
   		  \end{enumerate}
   	     \end{enumerate}
 \end{enumerate}
 
 \noindent\textit{D\'emonstration}\\
 \textit{Le (1) de (3.3.2.1)} r\'esulte de (3.1.4.2) et (1.2.1).\\
 \textit{Pour le (i) de (3.3.2.1)(2)} on consid\`ere le diagramme commutatif
 
   $$
\begin{array}{c}
\xymatrix{
S\  \ar@{^{(}->}[r]^{j_{T}} \ar @{} [dr] |{\square} & T\\
S'_{1} \ \ar@{^{(}->} [r] _{j_{1T'}} \ar [u] ^{\varphi _{1}} & T' \ar[u] _{\overline{\varphi}} \ar@{=}[d]\\
S' \ \ar@{^{(}->} [r] _{j_{T'}} \ar[u] ^{j} & T'
		} 
\end{array}
$$

\noindent dans lequel le carr\'e du haut est cart\'esien et $\varphi= \varphi_{1}\circ j$. On a alors une suite d'isomorphismes
$$ \begin{array}{lll@{\qquad}l}
(\varphi, \overline{\varphi}, \overline{g})^{\ \dag}(R^{i}\overline{h}_{K^{\ast}}j^{\dag}_{Y}E_{W} ) &\simeq & j_{T'}^{\dag} \overline{g}_{K}^{ \ast}R^{i}\overline{h}_{K^{\ast}}j^{\dag}_{Y}E_{W} & [(3.3.1.5)]\\
		&\overset{\sim}{\rightarrow}  & j^{\dag}_{T'} j^{\dag}_{1T'}g^{ \ast}_{V}R^{i}h_{V^{\ast}}E_{W} & [\mbox{B}\ 3,(2.1.4.8)], [\mbox{LS}, (5.3)]\\
		&\overset{\sim}{\rightarrow}  & j^{\dag}_{T'} g_{V}^{\ast} R^{i}{h}_{V^{\ast}}E_{W} & [\mbox{B}\ 3, (2.1.7)], [\mbox{LS}, (5.1)]\\
		& \overset{\sim}{\rightarrow}  & j^{\dag}_{T'}  R^{i} {h}'_{{V'}^{ \ast}}{g'}_{W}^{\ast}E_{W} & [(1.2.1)(2)]\\
		& \overset{\sim}{\leftarrow}  & R^{i}\overline{h}'_{K^{\ast}}j^{\dag}_{Y'}{g'}_{W}^{\ast}E_{W} &  [(3.1.4.2)]\\
		& \overset{\sim}{\leftarrow}  & R^{i}\overline{h}_{K^{\ast}} {\overline{g}'}_{K}^{\ast} j^{\dag}_{Y}E_{W} & [\mbox{B}\ 3,(2.1.4.8)], [\mbox{LS}, (5.3)].
\end{array}$$

\textit{Pour le (ii) de (3.3.2.1)(2)} il suffit de remarquer que les hypoth\`eses impliquent que
$(\varphi, \overline{\varphi}, \overline{g})^{\ \dag}(\mathcal{E})=\overline{g}_K^{\ast}(\mathcal{E})$ pour tout faisceau ab\'elien $\mathcal{E}$ sur $\mathcal{T}_K$; la derni\`ere assertion r\'esulte de (3.1.4.1), [B 3,(2.1.4.7)] ou [LS, 5.3] et (1.2.1) comme ci-dessus.
\\
\\
\textit{Pour (3.3.2.2)} on proc\`ede de m\^eme en utilisant cette fois le (1.2.2) du th\'eor\`eme (1.2). $\square$ \\
\\
\textit{Remarque (3.3.3)} En fait, dans le (3.3.2.1) (2) (ii) du th\'eor\`eme pr\'ec\'edent, si l'on ne suppose plus l'existence de $\overline{g}$, mais que l'on suppose toujours l'existence du carr\'e cart\'esien
$$
\begin{array}{c}
\xymatrix{
W'\ \ar@{->}[r]^{g'_W} \ar@{->}[d]_{h'_{V'}} & W \ar@{->}[d]^{h_V}\\
V'\ \ar@{->}[r]_{g_V} &V,
		}
\end{array}
$$
on obtient, pour tout $\mathcal{O}_{W}$-module coh\'erent $E_{W}$, un isomorphisme de changement de base
$$\begin{array}{ccc}
 g_{V}^{\ast}R^{i}h_{V^{\ast}}j^{\dag}_{W}E_{W} &\overset{\sim}{\rightarrow} & R^{i}h'_{V'^{\ast}}g'^{\ast}_{W}j_{W}^{\dag}E_{W} \\		
						&\stackrel{\sim}{\rightarrow} & R^{i}h'_{V'^{\ast}}j_{W'}^{\dag}g'^{\ast}_{W}E_{W} .
\end{array} 
\leqno{(3.3.3.1)}$$
De m\^eme, pour le (3.3.2.2)(2)(ii) du th\'eor\`eme, on a un isomorphisme
$$
 g_{V}^{\ast}R^{i}h_{V^{\ast}}{j_{W}}^{\dag}(E_{W}^{^{\bullet}}) \overset{\sim}{\rightarrow} R^{i}h'_{V'^{\ast}}j^{\dag}_{W'}(E_{W'}^{^{\bullet}}).
\leqno{(3.3.3.2)}
$$\\

\subsection*{3.4 Surconvergence des images directes.}

Nous allons consid\'erer dans le prochain th\'eor\`eme l'une des trois situations suivantes.\\

\textbf{(3.4.1) \textit{Dans le premier cas}}, nous consid\'erons un diagramme commutatif satisfaisant aux hypoth\`eses de (2.2)
$$
\begin{array}{c}
\xymatrix{
X\ \ar@{^{(}->}[r]^{j_{Y}} \ar @{}[dr] |{\square} \ar[d]_{f} & Y \ \ar@{^{(}->}[r]^{i_{Y}} \ar[d]^{\overline{f}} \ar @{}[dr] |{\square}& \mathcal{Y} \ar[d] ^{\overline{h}}\\
S \ \ar@{^{(}->}[r]_{j_{T}} & T\ \ar@{^{(}->}[r]_{i_{T}} & \mathcal{T} \ar[r]_{\rho} & \mathcal{W},\\
}
\end{array}
\leqno{(3.4.1.1)}
$$
et un diagramme commutatif
$$
\begin{array}{c}
\xymatrix{
S\ \ar@{^{(}->}[r]^{j_{T}} & T\ \ar@{^{(}->}[r]^{i_{T}}  & \mathcal{T} \ar[r] ^\rho&\mathcal{W}\\
S'\  \ar@{^{(}->}[r]_{j_{T'}} \ar@{->}[u]^{\varphi} & T'\ \ar@{^{(}->}[r]_{i_{T'}} \ar@{->}[u]^{\overline{\varphi}}& \mathcal{T}' \ar[r]_{\rho'} \ar@{->}[u]_{\overline{g}}& \mathcal{W}' \ar@{->}[u]_{\theta}\\
}
\end{array}
\leqno{(3.4.1.2)}
$$
\noindent tel qu'en prenant l'image inverse de (3.4.1.1) par (3.4.1.2) on obtienne un parall\'el\'epip\`ede commutatif tel que (3.3.1.1). On suppose de plus les carr\'es de (3.4.1.1) cart\'esiens ($\overline{h}^{\ -1}(T)=Y\ ,\ \overline{h}^{\ -1}(S)=X$), $ \overline{h}$ propre, $\overline{h}$ lisse sur un voisinage de X dans $\mathcal{Y}$, $\theta$ lisse, $\rho$ (resp $ \rho' $) lisse sur un voisinage de S dans $\mathcal{T}$ (resp de $S'$ dans $\mathcal{T}'$). On suppose \'egalement satisfaite l'une des deux hypoth\`eses suivantes: $\overline{g}$ est lisse sur un voisinage de $S' $ dans $\mathcal{T}'$, ou $\overline{g}$ est plat.\\

\textbf{(3.4.2) \textit{Dans le deuxi\`eme cas}}, nous consid\'erons un diagramme commutatif tel que (3.4.1.1) et un diagramme commutatif
$$
\begin{array}{c}
\xymatrix{
S\ \ar@{^{(}->}[r]^{j_{T}} & T\ \ar@{^{(}->}[r]^{i_{T}}  & \mathcal{T} \ar[r] ^\rho&\mathcal{W}\\
S'\  \ar@{^{(}->}[r]_{j_{T'}} \ar@{->}[u]^{\varphi} & T'\ \ar@{^{(}->}[r]_{i_{T'}} \ar@{->}[u]^{\overline{\varphi}}& \mathcal{T}' \ar[r]_{\rho'} & \mathcal{W}' \ar@{->}[u]_{\theta}\\
}
\end{array}
\leqno{(3.4.2.1)}
$$
satisfaisant aux m\^emes propri\'et\'es que (3.4.1.2) except\'e l'existence de $\overline{g}$, mais en supposant $ \rho\  propre$, et nous noterons
$$\xymatrix{X'=X \times_{S} S'\quad \ar@{^{(}->}[r]^{j_{Y'}}&\quad Y'=Y \times _{T}T'}$$
et
  $$
\begin{array}{c}
\xymatrix{
&X\ \ar@{.>} [dd]^(.3){f} |\hole  \ar@{^{(}->} [rr]^{j_{Y}} & & Y \ar [dd]^(.5){\overline{f}} \\
X' \ \ar@{^{(}->} [rr]^(.7){j_{Y'}} \ar [dd]_{f'} \ar[ur]^{\varphi'} & & Y' \ar[dd]_(.7){\overline{f}'} \ar [ur]^{\overline{\varphi}'} \\
& S \ \ar@{^{(}.>}[rr]^(.7){j_{T}} |\hole & &{T}\\
S'  \ \ar@{^{(}->}[rr]_{j_{T'}} \ar@{.>}[ur]_{\varphi}& & {T}' \ar[ur]_{\overline{\varphi}} 
}
\end{array}
\leqno{(3.4.2.2)}
  $$
l'image inverse par ($\varphi$, $\overline{\varphi}$) de (3.4.1.1).\\

\textbf{(3.4.3) \textit{Dans le troisi\`eme cas}}, qui g\'en\'eralise le premier, nous consid\'erons des diagrammes commutatifs tels que (3.4.1.1) et (3.4.2.1) mais sans supposer $\rho$ propre. Par contre nous supposons de plus l'existence d'un diagramme commutatif
$$
\begin{array}{c}
\xymatrix{
X' \ \ar@{^{(}->}[r]^{j_{Y'}} \ar @{}[dr] |{\square} \ar[d]_{f'} & Y' \ \ar@{^{(}->}[r]^{i_{Y'}} \ar[d]^{\overline{f}'} \ar @{}[dr] |{\square}& \mathcal{Y'} \ar[d] ^{\overline{h}'}\\
S' \ \ar@{^{(}->}[r]_{j_{T'}} & T'\ \ar@{^{(}->}[r]_{i_{T'}} & \mathcal{T}' \ar[r]^{\rho'} & \mathcal{W}'\\
}
\end{array}
\leqno{(3.4.3.1)}
$$
satisfaisant aux m\^emes hypoth\`eses que (3.4.1.1), et dans lequel $X', Y' $ satisfont aux propri\'et\'es de (3.4.2.2).\\

Dans le cas relevable le th\'eor\`eme suivant r\'esout une conjecture de Berthelot [B 2, (4.3)] et g\'en\'eralise le th\'eor\`eme 5 de loc. cit.\\
\\ 
\noindent\textbf{Th\'eor\`eme (3.4.4)}
\textit{Pour tout entier i $\geqslant$ 0, on a:
\begin{enumerate}
   \item[(3.4.4.1)] Sous les hypoth\`eses (3.4.3), on a:
      \begin{enumerate}
         \item[(i)]$\overline{f}$ induit un foncteur\\
         \\
                           		      $R^{i}\overline{f}_{rig \ast}((X,Y)/ \mathcal{T};-):{Isoc}^{\dag}((X,Y)/				\mathcal{W}) \longrightarrow {Isoc}^{\dag}((S,T)/W).$\\
		\item[(ii)] il existe un morphisme de changement de base\\
		\\
				${(\varphi,\overline{\varphi})}^{\ast}R^{i}\overline{f}_{rig \ast}((X,Y)/				\mathcal{T};E) \longrightarrow R^{i}\overline{f'}_{rig \ast}((X',Y')/				\mathcal{T}';{(\varphi',\overline{\varphi}')}^{\ast}(E))$\\
				\\
					et celui-ci est un isomorphisme dans ${Isoc}^{\dag}((S',T')/ \mathcal{W}')$.
	\end{enumerate}
     \item[(3.4.4.2)] Sous les hypoth\`eses (3.4.2) on a :
       \begin{enumerate}
       	\item[(i)] f induit un foncteur\\
	\\
		$R^{i}f_{rig \ast}(X/\mathcal{T};-):{Isoc}^{\dag}(X/\mathcal{W})  					\longrightarrow {Isoc}^{\dag}(S/\mathcal{W}).$\\
		\item[(ii)]L'isomorphisme de changement de base de (3.4.4.1)(ii) existe et devient\\
		\\
		${(\varphi,\overline{\varphi})}^{\ast}R^{i}f_{rig \ast}(X/				\mathcal{T};E) \overset{\sim}{\longrightarrow} R^{i}\overline{f'}_{rig \ast}((X',Y')/				\mathcal{T}';{(\varphi',\overline{\varphi}')}^{\ast}(E)).$\\
	\item[(iii)] Si de plus $\rho'$ est propre, alors l'isomorphisme de (ii) pr\'ec\'edent devient un isomorphisme dans ${Isoc}^{\dag}(S'/\mathcal{W}')$:\\
	\\
	${\varphi}^{\ast}R^{i}f_{rig \ast}(X/\mathcal{T};E) \overset{\sim}	{\longrightarrow} 	R^{i}{f}'_{rig \ast}(X'/\mathcal{T}';{\varphi'}^{\ast}(E)).$\\
	\item[(iv)] Si $S'=T'$, l'isomorphisme du (ii) pr\'ec\'edent devient un isomorphisme dans $ Isoc(S'/ \mathcal{W}'):$\\
	\\
	${\varphi}^{\ast}j_{T}^{\ast}R^{i}f_{rig \ast}(X/\mathcal{T};E) \overset{\sim}		{\longrightarrow} R^{i}{f}'_{conv \ast}(X'/\mathcal{T}'; {\varphi'}^{\ast}(\hat{E}))$\\
	\\
	o\`u $\hat{E}$ $\in$ $Isoc(X/\mathcal{W})$ est l'isocristal convergent associ\'e \`a E $\in$ ${Isoc}^{\dag}(X/\mathcal{W})$ par le foncteur d'oubli ${Isoc}^{\dag}(X/\mathcal{W})\ \rightarrow\ {Isoc}(X/\mathcal{W}).$\\
	En particulier si $S'=T'=S$, on a un isomorphisme \\
	\\
	$j_{T}^{\ast}R^{i}f_{rig \ast}(X/\mathcal{T};E) \overset{\sim}					{\longrightarrow} R^{i}{f}_{conv \ast}(X/\mathcal{T};\hat{E})$\\
	\end{enumerate}
    \item[(3.4.4.3)] Sous les hypoth\`eses (3.4.3) avec $S=T\  et \ S'=T'$ on a :
    	\begin{enumerate}
	    \item[(i)] f induit un foncteur\\
	    \\
	    $R^{i}f_{conv \ast}: Isoc(X/\mathcal{W}) \longrightarrow \ Isoc(S/\mathcal{W}).$\\
	    \item[(ii)] Il existe un isomorphisme de changement de base dans $Isoc(S'/ \mathcal{W}')$\\
	    \\
	    ${\varphi}^{\ast}R^{i}f_{conv \ast}(X/\mathcal{T};\mathcal{E}) \overset{\sim}				      {\longrightarrow} 	R^{i}{f}'_{conv \ast}(X'/\mathcal{T}';{\varphi'}^{\ast}		(\mathcal{E})).$
	 \end{enumerate}
 \end{enumerate}
	}
	
	\vskip10mm
Avant de donner la preuve du th\'eor\`eme, faisons quelques remarques:\\

\noindent\textit{$\underline{Remarques \ (3.4.4.4)}$}:\\

\begin{enumerate}
   \item[(i)] Sous les hypoth\`eses de (3.4.4.2)(iii), et en supposant de plus que $\varphi 	$ est l'identit\'e de S et $\theta$ l'identit\'e de $\mathcal{W}$, l'isomorphisme de 	changement de base prouve que $R^{i}f_{rig \ast}(X/\mathcal{T};E)$ est ind\'ependant du sch\'ema formel $\mathcal{T}$ dans lequel S est plong\'e (avec 		bien s\^ur $\mathcal{T}$ propre sur $\mathcal{W}$, $\rho$ lisse sur un voisinage 	de S dans $\mathcal{T}$ et $\overline{h}$ v\'erifiant (3.4.4)).\\
         La m\^eme remarque s'applique \`a $R^{i}f_{conv \ast}(X/\mathcal{T};\mathcal	{E}).$\\
         
   \item[(ii)] D'apr\`es [Et 5, th\'eo (3.2.1)] les hypoth\`eses de (3.4.4.2)(iii) sont v\'erifi\'ees pour $	\mathcal{W}=Spf \mathcal{V}, S$ affine et lisse sur $k$ et certains morphismes $f$ projectifs et lisses.\\
   
   \item[(iii)] En conjuguant (i) et (ii) nous en d\'eduirons plus loin [Th\'eor\`eme (3.4.8.2)] que les constructions se recollent pour certains morphismes $f$ projectifs lisses et $S$ un $k$-sch\'ema lisse ( plus n\'ecessairement affine).\\
     
\end{enumerate}

\noindent\textit{D\'emonstration de (3.4.4)}.\\

Le (3.4.4.3) est cons\'equence directe de (3.4.4.1).\\

La preuve de (3.4.4.1) et (3.4.4.2) va se faire en six \'etapes que l'on pr\'ecise ici:\\

Dans les quatre premi\`eres \'etapes on va faire la preuve de (3.4.4.1) sous l'hypoth\`ese plus particuli\`ere (3.4.1) (i.e. existence de $\overline{g}$):
\begin{enumerate}
  \item[-]\textit{\'etape $\  \rondI$}: montrer que $R^{i}\overline{f}_{rig \ast}((X,Y)/\mathcal{T};E)$ est un $j_{T}^{\dag}\mathcal{O}_{]T[_{\mathcal{T}}}$-module coh\'erent.
  \item[-]\textit{\'etape $\  \rondII$}: montrer l'isomorphisme de changement de base (3.4.4.1)(ii) 	lorsque $\overline{\varphi}^{-1}(S)=S'$.
  \item[-]\textit{\'etape $\  \rondIII$}: achever la preuve de (3.4.4.1)(i).
  \item[-]\textit{\'etape $\  \rondIV$}: achever la preuve de l'isomorphisme de changement de base (3.4.4.1)(ii).
  
  \item[-]\textit {\'etape $\  \rondV$}: on prouve (3.4.4.1) sous les hypoth\`eses (3.4.3).
  \item[-]\textit {\'etape $\  \rondVI$}: on prouve (3.4.4.2).\\
\end{enumerate}
Pla\c cons-nous d'abord sous les hypoth\`eses (3.4.1).\\
On se donne donc un diagramme tel que (3.3.1.1) avec $\overline{h}^{-1}(T)=Y$, $\overline{h}^{-1}(S)=X$: on note $S'_{1}=\overline{\varphi}^{-1}(S) \ \mbox{et}\ S' \xrightarrow[]{j}S_{1}^{\prime}\xrightarrow[]{\varphi_{1}}S $ la factorisation de $\varphi$ o\`u $j$ est une immersion ouverte; on en d\'eduit un diagramme commutatif \`a carr\'es verticaux cart\'esiens\\
\\
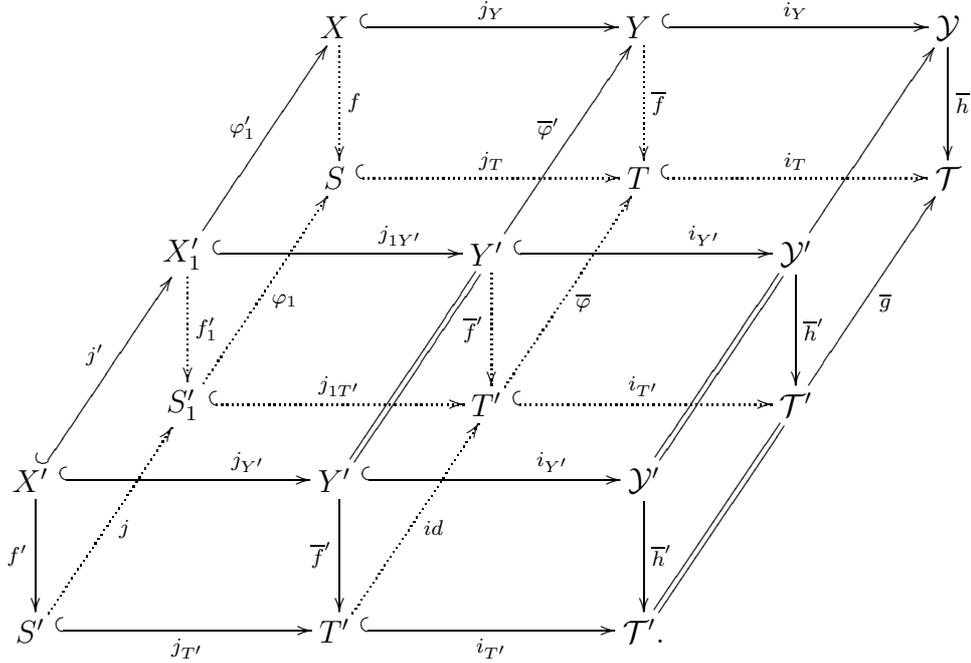
$$
\begin{array}{c}
 \shorthandoff{;:!?}
 \xymatrix@!0 @R=1cm @C=2cm{
&&X\ar@{.>}[dd]^{f}\  \ar @{^{(}->}[rr]^{j_{Y}}&&Y  \ar@{.>}[dd]^{\overline{f}}\  \ar@{^{(}->}[rr]^{i_{Y}}&&\mathcal{Y} \ar[dd]^{\overline{h}}\\
&&&&&&\\
&&S\  \ar@{^{(}.>}[rr]^{j_{T}}  && T\  \ar@{^{(}.>}[rr]^{i_{T}}&& \mathcal{T}\\
&X'_{1}\ar@{.>}[dd]^{f_{1}^{\prime}} \ar[uuur]^{\varphi'_{1}}\  \ar@{^{(}->}[rr]^(.7){j_{1Y'}}&&Y'\ar@{.>}[dd]_{\overline{f}'} \  \ar@{^{(}->}[rr]^(.7){i_{Y'}} \ar[uuur]^{\overline{\varphi}'}&& \mathcal{Y'} \ar[dd]^{\overline{h}'} \ar[uuur]\\
&&&&&&\\
& S'_{1}\  \ar@{^{(}.>}[rr]^{j_{1T'}} \ar@{.>}[uuur]_{\varphi_{1}} &&T' \  \ar@{^{(}.>}[rr] ^{i_{T'}} \ar@{.>}[uuur]_{\overline{\varphi}}&& \mathcal{T'} \ar[uuur]_{\overline{g}} &\\
X' \ar@{^{(}->}[rr]^(.7){j_{Y'}} \ar[dd]_{f'}\   \ar@{^{(}->}[uuur]^{j'}&&Y'\ar[dd]_{\overline{f}'}\  \ar@{^{(}->}[rr]^(.7){i_{Y'}} \ar@{=}[uuur]&& \mathcal{Y'}\ar[dd]^{\overline{h}'} \ar@{=}[uuur]&&\\
&&&&&&\\
S' \  \ar@{^{(}->}[rr]_{j_{T'}}  \ar@{.>}[uuur]_{j}&&T' \ar@{^{(}->}[rr]_{i_{T'}} \  \ar@{.>}[uuur]_{id}&&\mathcal{T'}. \ar@{=}[uuur]&&
}
\end{array}
$$
Pour $E \in Isoc^{\dag}((X,Y)/ \mathcal{W})$ on notera $(\varphi, \overline{\varphi})^{\ast}(E)$ son image inverse par le couple $(\varphi, \overline{\varphi})$ [B 3, (2.3.2)(iv)], [LS, chap 7] pour pr\'eciser la d\'ependance en $\varphi$ et $\overline{\varphi}$: d'autres images inverses seront utilis\'ees, que le contexte pr\'ecisera (cf. d\'ef. (3.3.1.5)).\\

\noindent\textbf{Etape $\  \rondI$.} Avec les notations de (3.2) il existe un voisinage strict  $W$ de $]X[_{\mathcal{Y}}$ dans $]Y[_{\mathcal{Y}}$ tel que $E_{\mathcal{Y}}:=j_{Y}^{\dag}E_{W}$ soit une r\'ealisation de E, et d'apr\`es (2.2.2.2) on peut supposer que $W$ est de la forme $W=\overline{h}_{Y}^{-1}(V)$ pour un voisinage strict $V$ de $]S[_{\mathcal{T}}$ dans $]T[_{\mathcal{T}}$; de plus on a un diagramme commutatif
$$
\begin{array}{c}
\xymatrix{
]X[_{\mathcal{Y}}\ \ar@{^{(}->}[r] \ar[d]_{\overline{h}_{X}} & W\
\ar@{^{(}->}[r]^{\alpha_{W}} \ar[d]_{{h}_{V}} & \mathcal{Y}_{K}  \ar[d]^{\overline{h}_{K}} \\
]S[_{\mathcal{T}}\ \ar@{^{(}->} [r]&V\ \ar@{^{(}->} [r]_{\alpha_{V}} & \mathcal{T}_{K}
} 
\end{array}
$$
dans lequel les carr\'es sont cart\'esiens et o\`u les fl\`eches horizontales sont des immertions ouvertes, ${\overline{h}_{K}}$ est propre et ${\overline{h}_{X}}$ est propre et lisse [(2.2.3.2)(i)]: quitte \`a restreindre $V$ on peut supposer via [(2.2.3.2)(ii)] que ${h}_{V}$ est propre et lisse.\\
Or $R^{i+j}\overline{f}_{rig^{\ast}}((X,Y)/\mathcal{T};E)$ est l'aboutissement d'une suite spectrale de terme $E_{1}^{i,j}$ donn\'e par
$$E_{1}^{i,j}=R^{j}{\overline{h}_{K^{\ast}}}(j_{Y}^{\dag}E_{W}\otimes_{\mathcal{O}_{]Y[_{\mathcal{Y}}}}\Omega^{i}_{]Y[_{\mathcal{Y}}/\mathcal{T}_{K}})$$
avec filtration
\begin{eqnarray*}
Fil^{i}&:= &Fil^{i}(j_{Y}^{\dag}E_{W}\otimes\Omega^{\bullet}_{]Y[_{\mathcal{Y}}/\mathcal{T}_{K}})\\
	&=&j_{Y}^{\dag}E_{W}\otimes\Omega^{\geqslant i}_{]Y[_{\mathcal{Y}}/\mathcal{T}_{K}},
\end{eqnarray*}
et on a une suite d'isomorphismes
$$ \begin{array}{lll@{\ }l}
j_{Y}^{\dag}E_{W}\otimes_{\mathcal{O}_{]Y[_{\mathcal{Y}}}}\Omega^{i}_{]Y[_{\mathcal{Y}}/\mathcal{T}_{K}}&=&(\alpha_{W^{\ast}}j_{W}^{\dag}E)\otimes_{\mathcal{O}_{]Y[_{\mathcal{Y}}}}\Omega^{i}_{]Y[_{\mathcal{Y}}/\mathcal{T}_{K}}& \\
&\simeq&\alpha_{W^{\ast}}(j_{W}^{\dag}E_{W}\otimes_{\mathcal{O}_{W}}\alpha_{W}^{\ast}(\Omega^{i}_{]Y[_{\mathcal{Y}}/\mathcal{T}_{K}}))&\\
&\simeq&\alpha_{W^{\ast}}(j_{W}^{\dag}E_{W}\otimes_{\mathcal{O}_{W}}\Omega^{i}_{W/\mathcal{T}_{K}})&\mbox{car}\ \alpha_{W}\ \mbox{est \'etale}\\
&\simeq&\alpha_{W^{\ast}}(j_{W}^{\dag}E_{W}\otimes_{\mathcal{O}_{W}}\Omega^{i}_{W/V})&\mbox{car}\ \alpha_{V}\  \mbox{est \'etale}\\
&\simeq&\alpha_{W^{\ast}}(j_{W}^{\dag}\mathcal{O}_{W}\otimes_{\mathcal{O}_{W}}E_{W}\otimes_{\mathcal{O}_{W}}\Omega^{i}_{W/V})&[\mbox{B}\ 3,(2.1.3)], [\mbox{LS, chap 5}]\\
&\simeq&\alpha_{W^{\ast}}j_{W}^{\dag}(E_{W}\otimes_{\mathcal{O}_{W}}\Omega^{i}_{W/V})&[\mbox{loc. cit.}]\\
&=&j_{Y}^{\dag}(E_{W}\otimes_{\mathcal{O}_{W}}\Omega^{i}_{W/V}),&\\
\end{array}$$
o\`u $\Omega^{i}_{W/V}$ est un $\mathcal{O}_{W}$-module coh\'erent et localement libre [(2.2.3.2(i)]. D'apr\`es (3.1.4) on en d\'eduit un isomorphisme
\begin{eqnarray*}
E_{1}^{i,j}&=&R^{j}{\overline{h}_{K^{\ast}}}({j}_{Y}^{\dag}E_{W}\otimes_{\mathcal{O}_{]Y[_{\mathcal{Y}}}}\Omega^{i}_{]Y[_{\mathcal{Y}}/\mathcal{T}_{K}})\\
&\overset{\sim}{\rightarrow}&j_{T}^{\dag}R^{j}{{h}_{V^{\ast}}}(E_{W}\otimes_{\mathcal{O}_{W}}\Omega^{i}_{W/V});
\end{eqnarray*}
or $R^{j}{{h}_{V^{\ast}}}(E_{W}\otimes_{\mathcal{O}_{W}}\Omega^{i}_{W/V})$ est un $\mathcal{O}_{V}$-module coh\'erent d'apr\`es le th\'eor\`eme (1.2), donc $E_{1}^{i,j}$ est un $j_{T}^{\dag}\mathcal{O}_{]T[_{\mathcal{T}}}$-module coh\'erent. Comme la filtration $Fil^{i}$ est finie, il en r\'esulte que l'aboutissement $R^{i+j}\overline{f}_{rig \ast}((X,Y)/\mathcal{T};E)$ est un $j_{T}^{\dag}\mathcal{O}_{]T[_{\mathcal{T}}}$-module coh\'erent. Remarquons que pour prouver cette coh\'erence on aurait pu appliquer le (1) de (3.3.2.2).\\
\\
\noindent\textbf{Etape \  \rondII}. On a vu \`a l'\'etape \  \rondI \ que 
\begin{eqnarray*}
R^{i}\overline{f}_{rig \ast}((X,Y)/\mathcal{T};E)&:=&R^{i}{\overline{h}_{K^{\ast}}}({j}_{Y}^{\dag}E_{W}\otimes_{\mathcal{O}_{]Y[_{\mathcal{Y}}}}\Omega^{\bullet}_{]Y[_{\mathcal{Y}}/\mathcal{T}_{K}})\\
&\simeq &R^{i}{\overline{h}_{{K}^{\ast}}}({j}_{Y}^{\dag}(E_{W}\otimes_{\mathcal{O}_{W}}\Omega^{\bullet}_{W/V})),
\end{eqnarray*}
et les $E_{W}\otimes_{\mathcal{O}_{W}}\Omega^{j}_{W/V}$ sont des $\mathcal{O}_{W}$-module coh\'erents.\\
Puisque ou bien $\overline{g}$ est plat, ou bien $\overline{g}$ est lisse sur un voisinage de $S'$ dans $\mathcal{T}'$, on peut d'apr\`es (3.3.1.4) et quitte \`a restreindre $V$ supposer $g_{V}$ plat: le (2) de (3.3.2.2) nous fournit alors l'isomorphisme de changement de base
$$\overline{g}_{K}^{\ast}R^{i}\overline{f}_{rig \ast}((X,Y)/\mathcal{T};E)\ \overset{\sim}{\rightarrow}\ R^{i}\overline{f'}_{rig \ast}((X'_{1},Y')/\mathcal{T}';{(\varphi'_{1},\overline{\varphi}')}^{\ast}(E)).
\leqno{(3.4.4.1.1)}
$$\\
\\
\noindent\textbf{Etape \  \rondIII }.  En reprenant la construction de la connexion de Gau\ss -Manin d\'ecrite explicitement par Katz dans [K, 3.4 et 3.5] on prouve que cette connexion agit aussi bien sur le terme $E_{1}^{i,j}$ que sur l'aboutissement [loc. cit., theo. 3.5]: ainsi $R^{i}{f}_{rig \ast}(X/\mathcal{T};E)$ est muni d'une connexion $\nabla^{i}$ (de Gau\ss -Manin) int\'egrable. Pour construire cette connexion de Gau\ss -Manin $\nabla^{i}$ on peut aussi \'etablir des isomorphismes [B 3, (2.2.5.1)], [LS, chap 7] v\'erifiant la condition de cocycles: cette construction co\"{\i}ncide avec la pr\'ec\'edente [B 1,V, 3.6.3, 3.6.4, 3.6.5], et cette deuxi\`eme construction va nous permettre d'\'etablir la surconvergence de $\nabla^{i}$.\\
On a un diagramme commutatif \`a carr\'es cart\'esiens
$$
\begin{array}{c}
\xymatrix{
{\mathcal{Y}}\ar@{->}[rr]^(.45){\Gamma_{\overline{h}}=(1_{\mathcal{Y}}\times\overline{h})} \ar[d]_{\overline{h}} & &{\mathcal{Y}}\times_{\mathcal{W}}\mathcal{T}\ar@{->}[rr]^(.55){p'_{1}} \ar[d]^{\overline{h}\times1_{\mathcal{T}}} & &\mathcal{Y} \ar[d]^{\overline{h}} \\
{\mathcal{T}}\ar@{->} [rr]_{\Delta_{\mathcal{T}}}& &\mathcal{T}\times_{\mathcal{W}}\mathcal{T}\ar@{->} [rr]_{p_{1\mathcal{T}}} & &\mathcal{T},
}
\end{array}
$$
o\`u $p_{1\mathcal{T}},\ p'_{1}$ sont les premi\`eres projections, $\Gamma_{\overline{h}}$ est le morphisme graphe de $\overline{h}$ et $\Delta_{\mathcal{T}}$ est le morphisme diagonal: $\Gamma_{\overline{h}}$ est une immersion ferm\'ee puisque $\mathcal{T}$ est s\'epar\'e sur $\mathcal{W}$; on note alors $]Y[_{\mathcal{Y}\times_{\mathcal{W}}\mathcal{T}}$ le tube de $Y$ dans ${\mathcal{Y}}\times_{\mathcal{W}}\mathcal{T}$ pour l'immersion ferm\'ee compos\'ee ${Y}\hooklongrightarrow\mathcal{Y}\overset{\Gamma_{\overline{h}}}{\hooklongrightarrow}{\mathcal{Y}}\times_{\mathcal{W}}\mathcal{T}$, et $ ]T[_{\mathcal{T}^{2}}$ celui de $T$ dans ${\mathcal{T}^{2}}$ pour l'immersion ferm\'ee $T\hooklongrightarrow\mathcal{T}\underset{\Delta_{\mathcal{T}}}{\hooklongrightarrow}{\mathcal{T}}\times_{\mathcal{W}}\mathcal{T}$.\\
D'apr\`es (2.1.2) on a alors un diagramme commutatif \`a carr\'e cart\'esien
$$
\begin{array}{c}
\xymatrix{
 &  & &  ]Y[_{\mathcal{Y}^{2}}\ar@/^1.5pc/[ddl]^{\overset{\sim}{h}} \ar@/_1.5pc/[dlll]_{p_{1 \mathcal{Y}}} \ar@{->}[dl]^{h''_{1}}\\
]Y[_{\mathcal{Y}}\ar@{->}[d]_{\overline{h}_{Y}} & &]Y[_{\mathcal{Y}\times\mathcal{T}} \ar@{->}[d]^{h'_{1}} \ar@{->}[ll]_{p'_{1}}&\\
]T[_{\mathcal{T}} & & ]T[_{\mathcal{T}^{2}} \ar@{->}[ll]^(.45){p_{1_{\mathcal{T}}}}&
	}
\end{array}
\leqno{(3.4.4.1.2)}
$$
o\`u $h'_{1}$, $h''_{1}$, $ \overset{\sim}{h}$ sont induits respectivement par $\overline{h}_{K}\times1_{\mathcal{T}_{K}}$, $1_{\mathcal{Y}_{K}}\times\overline{h}_{K}$ et $\overline{h}_{K}\times\overline{h}_{K}$.\\
On a de m\^eme un diagramme analogue $(3.4.4.1.2)' $ avec $p_{2}$ \`a la place de $p_{1}$ et $h'_{2},h''_{2}$...\\
D'apr\`es (3.2.2.1) et [B 5, (3.2.3) et (3.1.11)(i)] ou [LS, 7.4] on a les isomorphismes canoniques
 \begin{eqnarray*}
R^{i}\overline{f}_{rig \ast}((X,Y)/\mathcal{T}^{2};E)&=&R^{i}{\overset{\sim}{h}_{\ast}}({p}_{1\mathcal{Y}}^{\ast}(E_{\mathcal{Y}})\otimes\Omega^{\bullet}_{]Y[_{\mathcal{Y}^{2}}/\mathcal{T}_{K}^{2}})\\
&\simeq &R^{i}{{h}'_{1\ast}}({p'}_{1}^{\ast}(E_{\mathcal{Y}})\otimes\Omega^{\bullet}_{]Y[_{\mathcal{Y}\times\mathcal{T}}/\mathcal{T}_{K}^{2}})\\
&\overset{\sim}{\leftarrow}&{p}_{1\mathcal{T}}^{\ast}R^{i}\overline{h}_{Y}^{\ast}(E_{\mathcal{Y}}\otimes\Omega^{\bullet}_{]Y[_{\mathcal{Y}}/\mathcal{T}_{K}})\qquad [\mbox{\'etape \  \rondII} ]\\
&=&{p}_{1\mathcal{T}}^{\ast}R^{i}\overline{f}_{rig \ast}((X,Y)/\mathcal{T};E);
\end{eqnarray*}
l'avant-dernier isomorphisme ci-dessus est r\'ealisable via l'\'etape \ \rondII \ car $\rho$ \'etant lisse sur un voisinage de $S$ dans $\mathcal{T}$, $p_{1\mathcal{T}}:\mathcal{T}\times_{\mathcal{W}}\mathcal{T}\rightarrow\mathcal{T}$ est lisse sur un voisinage de $S$ dans $\mathcal{T}\times_{\mathcal{W}}\mathcal{T}$.\\
Or l'isomorphisme [B 3, (2.2.5.1)], [LS, 7.2]
$$
p_{2\mathcal{Y}}^{\ast}(E_{\mathcal{Y}})\overset{\sim}{\longrightarrow}p_{1\mathcal{Y}}^{\ast}(E_{\mathcal{Y}})
$$
assurant l'existence d'une connexion (surconvergente) sur $E_{\mathcal{Y}}$ fournit aussi les isomorphismes
$$
R^{i}{\overset{\sim}{h}_{\ast}}({p}_{2\mathcal{Y}}^{\ast}(E_{\mathcal{Y}})\otimes\Omega^{\bullet}_{]Y[_{\mathcal{Y}^{2}}/\mathcal{T}_{K}^{2}})\overset{\sim}{\longrightarrow}R^{i}{\overset{\sim}{h}_{\ast}}({p}_{1\mathcal{Y}}^{\ast}(E_{\mathcal{Y}})\otimes\Omega^{\bullet}_{]Y[_{\mathcal{Y}^{2}}/\mathcal{T}_{K}^{2}}),
$$
c'est-\`a-dire, par les m\^emes arguments que ci-dessus, des isomorphismes\\
$$
(3.4.4.1.3)\ {p}_{2\mathcal{T}}^{\ast}R^{i}\overline{f}_{rig\ast}((X,Y)/\mathcal{T};E)\overset{\sim}{\rightarrow}R^{i}\overline{f}_{rig\ast}((X,Y)/\mathcal{T}^{2};E)\overset{\sim}{\leftarrow}{p}_{1\mathcal{T}}^{\ast}R^{i}\overline{f}_{rig\ast}((X,Y)/\mathcal{T};E)\\
$$
satisfaisant aux conditions de [B 3, (2.2.5)], [LS, 7.2]. De plus la connexion surconvergente sur $R^{i}\overline{f}_{rig\ast}((X,Y)/\mathcal{T};E)$ obtenue par l'isomorphisme (3.4.4.1.3) est bien celle de Gau\ss-Manin [B 1,V, 3.6.4].\\
 Ceci ach\`eve la preuve de l'\'etape \ \  \rondIII.  
  \\
  
\noindent\textbf{Etape \  \rondIV }. Puisque $R^{i}\overline{f}_{rig\ast}((X,Y)/\mathcal{T};E)$ est un isocristal surconvergent le long de $T\setminus S$, son image inverse par $(\varphi, \overline{\varphi})$ est, d'apr\`es l'\'etape\  \  \rondII , l'image inverse par $(j, id_{T'})$ de l'isocristal surconvergent (le long de $T'\backslash S^{\prime}_{1}$) \\
$$R^{i}\overline{f'}_{rig\ast}((X'_{1},Y')/\mathcal{T}';{(\varphi'_{1}, \overline{\varphi}')}^{\ast}(E))=j_{1T'}^{\dag}(R^{i}h'_{{V'}^{\ast}}({g}_{W}^{\prime\ast}(E_{W})\otimes_{\mathcal{O}_{W'}}\Omega^{\bullet}_{W'/V'}))\\\mbox{(cf (3.3.1.3))}.$$
Le $\mathcal{O}_{V'}$-module
$$\mathcal{E}_{V'}^{i}:=R^{i}h'_{{V'}^{\ast}}({g}_{W}^{\prime\ast}(E_{W})\otimes_{\mathcal{O}_{W'}}\Omega^{\bullet}_{W'/V'})
$$
est coh\'erent, et par d\'efinition des images inverses [B 3, (2.3.2)(iv)], [LS, chap 5] on a
$$
{(j, id_{T'})}^{\ast}(j_{1T'}^{\dag}(\mathcal{E}_{V'}^{i}))=j_{{T'}}^{\dag}(\mathcal{E}_{V'}^{i});
$$
or d'apr\`es [(3.3.2.2)(1)] on a:
\begin{eqnarray*}
j_{{T'}}^{\dag}(\mathcal{E}_{V'}^{i})&\simeq&R^{i}\overline{h}'_{{K}^{\ast}}(j_{Y'}^{\dag}({g}_{W}^{\prime\ast}(E_{W}))\otimes_{\mathcal{O}_{W'}}\Omega^{\bullet}_{W'/V'})\\
&=:&R^{i}\overline{f}'_{rig^{\ast}}((X',Y')/\mathcal{T}';{(\varphi, \overline{\varphi})}^{\ast}(E)),
\end{eqnarray*}
d'o\`u l'isomorphisme de changement de base
$$
(3.4.4.1.4)\ {(\varphi, \overline{\varphi})}^{\ast}R^{i}{\overline{f}}_{rig\ast}((X,Y)/\mathcal{T};E)\simeq R^{i}\overline{f}'_{rig\ast}((X',Y')/\mathcal{T}';{(\varphi', \overline{\varphi'})}^{\ast}(E)),
$$
qui est celui de (3.4.4.1)(ii).\\
\\
\noindent\textbf{Etape  \  \rondV}. Pla\c{c}ons-nous sous les hypoth\`eses (3.4.3) et prouvons (3.4.4.1) dans ce cas.\\
Consid\'erons les parall\'el\'epip\`edes commutatifs suivants dans lesquels les faces verticales sont cart\'esiennes

 $$
\shorthandoff{;:!?}
\begin{array}{l}
\xymatrix@!0 @R=1,4cm @C=1,4cm{
& X \ar@{.>}[dd]^(.3){f}  \ar@{^{(}->}[rr]^{j_{Y}} & & Y\ar@{.>}[dd]^(.3){\overline{f}}  \ar@{^{(}->}[rr]^{i_{Y}} & &\mathcal{Y} \ar[dd]^{\overline{h}} \\
X^{'}  \ar@{^{(}->}[rr]^(.8){j_{Y'}} \ar[dd]_{f^{'}} \ar[ur]^{\varphi'} & & Y' \ar@{^{(}->}[rr]^(.8){i''_{Y'}} \ar[dd]_(.7){\overline{f}'} \ar[ur]^{{\overline{\varphi}}'} & & \mathcal{Y}'' \ar[dd]_(.7){\overline{h}''} \ar[ur]^{p_{1\mathcal{Y}}} \\
& S\ar@{^{(}.>}[rr]^(.7){j_{T}}  & &T\ar@{^{(}.>}[rr]^(.8){i_{T}}   & &  \mathcal{T}\ar[rr]_{\rho}&& \mathcal{W}\\
S^{'}  \ar@{^{(}->}[rr]_{j_{T'}} \ar@{.>}[ur]^{\varphi}& & T' \ar@{^{(}->}[rr]_{i''_{T'}} \ar@{.>}[ur]_{\overline{\varphi}} & & **[r]\mathcal{T''}=\mathcal{T}\times_{\mathcal{W}}\mathcal{T}' \ar[ur]_{p_{1 \mathcal{T}''}} 
}
\end{array}
\leqno{(3.4.4.1.5)}
  $$
 $$
  \begin{array}{c}
\xymatrix{
& X' \ar@{.>}[dd]^(.3){f^{'}}  \ar@{^{(}->}[rr]^{j_{Y'}} & & Y'\ar@{.>}[dd]^(.3){\overline{f}'}  \ar@{^{(}->}[rr]^{i_{Y'}'''} & &\mathcal{Y}''' \ar[dd]^{\overline{h}\ '''} \ar[dl]  \\
X^{'}  \ar@{^{(}->}[rr]^(.8){j_{Y'}} \ar[dd]_{f^{'}} \ar@{=}[ur] & & Y' \ar@{^{(}->}[rr]^(.8){i_{Y'}} \ar[dd]_(.7){\overline{f}'} \ar@{=}[ur] & & \mathcal{Y}' \ar[dd] _(.7){\overline{h}'}\\
& S'\ar@{^{(}.>}[rr]^(.7){j_{T'}}  \ar@{.>}[dl]_{id}& &T'\ar@{^{(}.>}[rr]^(.8){i''_{T'}}  \ar@{.>} [dl] ^{id}& &  \mathcal{T}'' \ar[dl]^{p_{2 \mathcal{T}''}}\\
S^{'}  \ar@{^{(}->}[rr]_{j_{T'}} & & T' \ar@{^{(}->}[rr]_{i_{T'}} & &\mathcal{T}'  \ar[rr]_{\rho'}&& \mathcal{W}' 
}
\end{array}
\leqno{(3.4.4.1.6)}
  $$
o\`u $p_{i}$ est la projection sur le $i \ieme $  facteur et $i''_{T'}=(i_{T}\circ\overline{\varphi},i_{T'})$.\\
On forme aussi le carr\'e cart\'esien
$$
\begin{array}{c}
\xymatrix{
\overset{\sim}{\mathcal{Y}} \ar[rr]^{u_2} \ar[d]_{u_{1}}& &\mathcal{Y}'' \ar[d]^{\overline{h}''}\\
\mathcal{Y}''' \ar[rr]_{\overline{h}'''}& &\mathcal{T}'' \   .
}
\end{array}
\leqno{(3.4.4.1.7)}
$$
Comme $u_{1}$ et $u_{2}$ sont propres, et lisses sur un voisinage de $X'$ dans $\widetilde{\mathcal{Y}}$, on peut d'apr\`es [B 5, (3.1.2)] ou [LS, 7.4] faire le calcul de 
$$R^{i}\overline{f}'_{rig\ast}((X',Y')/\mathcal{T}'';{(\varphi', \overline{\varphi'})}^{\ast}(E))$$
\`a l'aide de la cohomologie de de Rham, aussi bien avec $\overline{h}'''$, $\overline{h}''' \circ u_{1} = \overline{h}'' \circ u_{2}$ ou $\overline{h}''$, qui sont tous les trois propres, et lisses sur un voisinage de $X'$ respectivement dans $\mathcal{Y}''', \widetilde{\mathcal{Y}}$ et $\mathcal{Y}''$. Or $p_{1 \mathcal{T}''}$ (resp $ p_{2\mathcal{T}''}$) \'etant lisse sur un voisinage de $S'$ dans $\mathcal{T}''$, on a d'apr\` es l'\'etape\  \  \rondIV  \ des isomorphismes de changement de base (le premier calcul\'e via $\overline{h}''$ et le second via $\overline{h}'''$)
$$ {(\varphi, \overline{\varphi}, p_{1\mathcal{T}''})}^{\ast}R^{i}{\overline{f}}_{rig\ast}((X,Y)/\mathcal{T};E)\simeq R^{i}\overline{f}'_{rig\ast}((X',Y')/\mathcal{T}'';{(\varphi', \overline{\varphi'})}^{\ast}(E)) \ ,$$
$$ {(p_{2\mathcal{T}''})}^{\ast}\ : \ R^{i}\overline{f}'_{rig\ast}((X',Y')/\mathcal{T}' ;{(\varphi', \overline{\varphi'})}^{\ast}(E))\simeq R^{i}\overline{f}'_{rig\ast}((X',Y')/\mathcal{T}'';{(\varphi', \overline{\varphi'})}^{\ast}(E))$$
et les seconds membres co\"{\i}ncident en faisant le calcul de la cohomologie de de Rham via $\overline{h}''' \circ u_{1} = \overline{h}'' \circ u_{2}$. Ceci ach\`eve la preuve de l'isomorphisme de changement de base: d'o\`u, ici encore, l'isomorphisme (3.4.4.1.3). Comme l'\'etape \ \rondI\ s'applique on a prouv\'e (3.4.4.1).\\

\noindent\textbf{Etape \  \rondVI }. L'\'etape \ \rondIII \ reste inchang\'ee puisqu'on peut faire les changements de base par $p_{1\mathcal{T}}$ et $p_{2\mathcal{T}}$ car ils sont lisses sur un voisinage de $S$ dans $\mathcal{T}\times_{\mathcal{W}}\mathcal{T}$: en particulier $R^{i}{\overline{f}}_{rig\ast}((X,Y)/\mathcal{T};E)$ est un isocristal surconvergent le long de $T\setminus S$. Soit $\mathcal{W}_{0}$ la r\'eduction sur $k$ de $\mathcal{W}$. On note \\
$$
\begin{array}{c}
\xymatrix{
X'' \ar@{^{(}->}[r]^{j_{Y''}} \ar[d]_{f''}&Y'' \ar@{^{(}->}[r]^{i_{Y''}} \ar[d]_{\overline{f}''} & **[r] \mathcal{Y}''=\mathcal{Y}\times_{\mathcal{W}}\mathcal{T}' \ar[d]_{\overline{h}''}\\
**[l]S''=S\times_{\mathcal{W}_o}S' \ar@{^{(}->}[r]_(.3){j_{T''}=j_{T}\times j_{T'}}&T'':=T\times_{\mathcal{W}_o}T'\ar@{^{(}->}[r]_(.6){i_{T''}=i_{T}\times i_{T'}}&**[r] \mathcal{T}'':=\mathcal{T}\times_{\mathcal{W}}\mathcal{T}'
}
\end{array}
\leqno{(3.4.4.2.1)}
$$
l'image inverse de (2.2.1) par le diagramme commutatif
$$
\begin{array}{c}
\xymatrix{
S\ \ar@{^{(}->}[rr]^{j_{T}} &&T\ \ar@{^{(}->}[rr]^{i_{T}}&&\mathcal{T}\\
S''\ \ar@{^{(}->}[rr]_{j_{T''}} \ar[u]^{p_{1S}} &&T'' \ar@{^{(}->}[rr]_{i_{T''}} \ar[u]_{p_{1T}}&&\mathcal{T}''\ \ar[u]_{p_{1\mathcal{T}''}}
}
\end{array}
\leqno{(3.4.4.2.2)}
$$
o\`u les $p_{1}$ sont les projections sur le premier facteur: $\overline{h}''$ est propre et $\overline{h}''$ est lisse sur un voisinage de $X''$ dans $\mathcal{Y}''$.\\
De m\^eme 
$$
\begin{array}{c}
\xymatrix{
X'\ \ar@{^{(}->}[rr]^{j_{Y'}} \ar[d]_{f'}&&Y'\ \ar@{^{(}->}[rr]^{i_{Y''}\circ\psi_Y} \ar[d]^{\overline{f}'} &&\mathcal{Y}'' \ar[d]^{\overline{h}''}\\
S'\ \ar@{^{(}->}[rr]_{j_{T'}}&&T'\ \ar@{^{(}->}[rr]_{i_{T''}\circ\psi_T}&&\mathcal{T}''
}
\end{array}
\leqno{(3.4.4.2.3)}
$$
est l'image inverse de (3.4.4.2.1) par le diagramme commutatif
$$
\begin{array}{c}
\xymatrix{
S''\ \ar@{^{(}->}[rr]^{j_{T''}} &&T''\ \ar@{^{(}->}[rr]^{i_{T''}}&&\mathcal{T}''\\
S'\ \ar@{^{(}->}[rr]_{j_{T'}} \ar[u]^{\psi_{S}=(\varphi,1_{S'})
} &&T'\ \ar@{^{(}->}[rr]_{i_{T''}\circ\psi_T} \ar[u]_{\psi_T=(\overline{\varphi},1_{T'})}&&\mathcal{T}''\ \ar@{=}[u]_{id_{\mathcal{T}''}}
}
\end{array}
\leqno{(3.4.4.2.4)}
$$
o\`u $\psi_{S}, \psi_{T}$ sont des immersions ferm\'ees et $\psi_{Y}$ est l'immersion ferm\'ee d\'efinie par le carr\'e cart\'esien
$$
\xymatrix{
Y'\ \ar[d]_{\overline{f}'} \ar@{^{(}->}[rr]^{\psi_Y} && Y'' \ar[d]^{\overline{f}''}\\
T'\ \ar@{^{(}->}[rr]_{\psi_T} && T''\ .
}
$$
En appliquant (3.4.4.1)(ii) au changement de base par $p_{1\mathcal{T}}$ et $\psi_{T}$ on obtient un isomorphisme 
$$ \psi_{T}^{\ast}p_{1\mathcal{T}''}^{\ast}R^{i}{\overline{f}}_{rig\ast}((X,Y)/\mathcal{T};E)\simeq R^{i}\overline{f}'_{rig\ast}((X',Y')/\mathcal{T}'';{(\varphi', \overline{\varphi'})}^{\ast}(E)) .$$
Puisque le carr\'e suivant commute
$$
\xymatrix{
T\ \ar@{^{(}->}[rr]^{i_T} && \mathcal{T} \\
T'\ \ar@{^{(}->}[rr]_{i_{T''}\circ\psi_T} \ar[u]^{p_{1T}\circ\psi_T= \overline{\varphi}}&& \mathcal{T}'' \ar[u]_{p_{1\mathcal{T}''}} \ ,
}
$$
pour prouver (3.4.4.2) il suffit d'apr\`es [B 3, (2.3.2)(iv)] ou [LS, chap 7] de prouver que la projection sur le second facteur $p_{2\mathcal{T}''}: \mathcal{T}'' \rightarrow \mathcal{T}'$ induit un isomorphisme 
$$ p_{2\mathcal{T}''}^{\ast}R^{i}\overline{f}'_{rig\ast}((X',Y')/\mathcal{T}';{(\varphi', \overline{\varphi'})}^{\ast}(E))\simeq R^{i}\overline{f}'_{rig\ast}((X',Y')/\mathcal{T}'';{(\varphi', \overline{\varphi'})}^{\ast}(E)) .$$
Consid\'erons alors le diagramme commutatif \`a carr\'e cart\'esien

$$
\begin{array}{c}
\xymatrix{
X'\ar@{^{(}->}[r]^{j_{Y'}}&Y'\ar@{^{(}->}[r]^{\psi_{Y}}&Y''\ar@{^{(}->}[r]^{i_{Y''}}&\mathcal{Y}'' \ar@/_1.5pc/[dddr]_{\overline{h}''} \ar[dr]_{i_{\mathcal{Y}''}} \ar@{=}[drrr]^{Id_{\mathcal{Y}''}}\\
&&&&\mathcal{Y}'''\ar[rr]_{p_{3\mathcal{Y}}} \ar[dd]^{\overline{h}'''}&&\mathcal{Y}''\ar[d]^{\overline{h}''}\\
&&&&&&\mathcal{T}'' \ar[d]^{p_{2\mathcal{T}''}}\\
&&&&\mathcal{T}''\ar[rr]_{p_{2\mathcal{T}''}}&&\mathcal{T}'\ .
}
\end{array}
$$
D'apr\`es les hypoth\`eses faites sur $\rho$, $p_{2\mathcal{T}''}$ est propre et $p_{2\mathcal{T}''}$ est lisse sur un voisinage de $S''$ dans $\mathcal{T}''$: ainsi $p_{2\mathcal{T}''}\circ \overline{h}''$ est propre et $p_{2\mathcal{T}''}\circ \overline{h}''$ est lisse sur un voisinage de $X'$ dans $\mathcal{Y}''$; par suite en appliquant (3.4.4.1) on en d\'eduit que 
$$R^{i}\overline{f}'_{rig\ast}((X',Y')/\mathcal{T}';{(\varphi', \overline{\varphi'})}^{\ast}(E)) \ \mbox{et} \ R^{i}\overline{f}'_{rig\ast}((X',Y')/\mathcal{T}'';{(\varphi', \overline{\varphi'})}^{\ast}(E))$$
sont des isocristaux surconvergents et puisque $  i_{\mathcal{Y}''} $ est propre, que l'on a des isomorphismes (o\`u $p_{1\mathcal{Y}}: \mathcal{Y}'' \rightarrow \mathcal{Y}$ est la premi\`ere projection)
$$
\begin{array}{rcl}
p_{2\mathcal{T}''}^{\ast}R^{i}\overline{f}'_{rig\ast}((X',Y')/\mathcal{T}';{(\varphi', \overline{\varphi'})}^{\ast}(E))& \simeq& R^{i}\overline{h}_{\ast}'''{p}_{3\mathcal{Y}}^{\ast}((p_{1\mathcal{Y}}^{\ast}E_{\mathcal{Y}})\otimes\Omega^{\bullet}_{]Y'[_{\mathcal{Y}''}/ \mathcal{T}'_{K}})\\
&\simeq & R^{i}\overline{h}'''_{\ast}({p}_{3\mathcal{Y}}^{\ast}p_{1\mathcal{Y}}^{\ast}(E_{\mathcal{Y}})\otimes\Omega^{\bullet}_{]Y'[_{\mathcal{Y}^{'''}}/ \mathcal{T}_{K}^{''}}) \ [(2.1.2)]\\
&\simeq & R^{i}\overline{h}''_{\ast}( {i}_{\mathcal{Y}^{''}}^{\ast}{p}_{3\mathcal{Y}}^{\ast}{p}_{1\mathcal{Y}}^{\ast}(E_{\mathcal{Y}})\otimes\Omega^{\bullet}_{]Y'[_{\mathcal{Y}^{''}}/ \mathcal{T}_{K}^{''}})\  [\mbox{B}\ 5, (3.2.2)] \\
&= & R^{i}\overline{f}'_{rig\ast}((X',Y')/\mathcal{T}^{''};{(\varphi', \overline{\varphi'})}^{\ast}(E))\  [(3.2.2)];
\end{array}
$$
au lieu de [B 5, (3.2.2)] ci-dessus on peut utiliser [LS, 7.4]. D'o\`u (3.4.4.2), compte tenu des d\'efinitions (3.2). Ceci ach\`eve la preuve de (3.4.4). $\square$  \\

\noindent\textbf{(3.4.5)}\\

(3.4.5.1) Supposons donn\'e un diagramme commutatif
$$
\xymatrix{
X\ \ar@{^{(}->}[rr]^{j_{\mathcal{Y}}} \ar[d]_f \ar @{}[drr] |{\square}&& \mathcal{Y}\ar[d]^{\overline{h}}&&\\
S\ \ar@{^{(}->}[rr]_{j_{\mathcal{T}}} &&\mathcal{T}\ar[rr]^{\rho}&&\mathcal{W},
}
$$
dans lequel le carr\'e est cart\'esien, $f$ est un morphisme propre de $k$-sch\'emas s\'epar\'es de type fini, $\overline{h}$ et $\rho$ sont des morphismes de $\mathcal{V}$-sch\'emas formels s\'epar\'es plats de type fini, $\overline{h}$ est propre, $\overline{h}$ (resp $\rho$) est lisse sur un voisinage de $X$ (resp $S$) dans $\mathcal{Y}$ (resp $\mathcal{T}$), $j_{\mathcal{Y}} \ \mbox{et} \ j_{\mathcal{Y}}$ sont des immersions. Soit $T$ l'adh\'erence sch\'ematique de $S$ dans $\mathcal{T}$ et 
$$S\underset{j_{T}}{\hooklongrightarrow} T\underset{i_{T}}{\hooklongrightarrow}{\mathcal{T}}$$
la factorisation de $j_{\mathcal{T}}$: $j_{T}$ est une immersion ouverte dominante et $i_{T}$ est une immersion ferm\'ee. On note $Y= \overline{h}^{-1}(T)$ et $\overline{f}:= \overline{h}_{|_{Y}}: Y \rightarrow T.$ Avec les notations pr\'ec\'edentes on est donc dans la situation (2.2.1) avec $\overline{f}^{-1}(S) = X.$\\

(3.4.5.2) Supposons de plus donn\'e un diagramme commutatif
$$
\begin{array}{c}
\xymatrix{
S\ \ar@{^{(}->}[rr]^{j_{\mathcal{T}}} &&\mathcal{T} \ar[rr]^{\rho}&&\mathcal{W}\\
S'\ \ar@{^{(}->}[rr]^{j_{\mathcal{T}'}} \ar[u]^{\varphi} &&\mathcal{T}' \ar[rr]^{\rho'} \ar[u]_{\overline{g}}&&\mathcal{W}'\ \ar[u]_{\theta}
}
\end{array}
$$
dans lequel $\varphi$ est un morphisme de $k$-sch\'emas s\'epar\'es de type fini, $\overline{g}, \ \rho'$ sont des morphismes s\'epar\'es de $\mathcal{V}$-sch\'emas formels s\'epar\'es plats de type fini, $\rho'$ est lisse sur un voisinage de $S'$ dans $\mathcal{T}'$, $\theta$ est lisse et $j_{\mathcal{T}'}$ est une immersion.\\
On d\'efinit $f':X'\rightarrow S'$ par le carr\'e cart\'esien
$$
\xymatrix{
X' \ar[rr]^{\varphi'} \ar[d]_{f'} && X \ar[d]^{f}\\
S' \ar[rr]_{\varphi}&&S\ .
}
$$
On note $T'$ l'adh\'erence sch\'ematique de $S'$ dans $\mathcal{T'}$ et 
$$S'\underset{j_{T'}}{\hooklongrightarrow} T' \underset{i_{T'}}{\hooklongrightarrow}{\mathcal{T}'}$$
la factorisation de $j_{\mathcal{T}'}$ par l'immersion ouverte dominante $j_{T'}$ et l'immersion ferm\'ee $i_{T'}$: $\overline{g}$ induit un $k$-morphisme $\overline{\varphi}:T' \rightarrow T\ .$\\
Par image inverse de (3.4.5.1) par (3.4.5.2) on obtient un parall\'el\'epip\`ede commutatif tel que (3.3.1.1), dont on reprend les notations; on en d\'eduit:\\ 

\noindent\textbf{Corollaire (3.4.5.3).}
\textit { Sous les hypoth\`eses (3.4.5.1) et (3.4.5.2) les parties (3.4.4.1) et (3.4.4.3) du th\'eor\`eme (3.4.4) sont valides.}\\

(3.4.5.4) Supposons donn\'e cette fois un diagramme commutatif
$$
\xymatrix{
S\ \ar@{^{(}->}[rr]^{j_{\mathcal{T}}}&&\mathcal{T}\ar[rr]^{\rho}&&\mathcal{W}\\
S'\ \ar[u]^{\varphi} \ar@{^{(}->}[rr]_{j_{\mathcal{T}'}}&&\mathcal{T}'\ar[rr]_{\rho'}&&\mathcal{W}' \ar[u]_{\theta}
}
$$
dans lequel $\varphi$ est un morphisme de $k$-sch\'emas s\'epar\'es de type fini, $\rho, \rho' ,\theta$ sont des morphismes (s\'epar\'es) de $\mathcal{V}$-sch\'emas formels s\'epar\'es plats de type fini, $\theta$ est lisse, $\rho$ (resp $\rho'$) est lisse sur un voisinage de $S$ dans $\mathcal{T}$ (resp de $S'$ dans $ \mathcal{T}' $), $ j_{\mathcal{T}} $ et $ j_{\mathcal{T}' }$ sont des immersions. Et on suppose de plus que $ \rho $ est propre. D'apr\`es l'\'etape \ \rondVI \  \ de la d\'emonstration du th\'eor\`eme (3.4.4) on a montr\'e:\\

\noindent\textbf{Corollaire (3.4.5.5).}
\textit { Sous les hypoth\`eses (3.4.5.1) et (3.4.5.4) la partie (3.4.4.2) du th\'eor\`eme (3.4.4) s'applique.}\\

\noindent\textbf{3.4.6.} Sous les hypoth\`eses de (3.4.5.1) avec $k$ parfait supposons que $\mathcal{W}= Spf\ \mathcal{V}$ et que $\rho:\  \mathcal{T} \rightarrow Spf\ \mathcal{V}$ est lisse sur un ouvert $\mathcal{S}$ de $\mathcal{T}$, avec $S$ contenu dans $\mathcal{S}$.\\

Soient $i_{s}: s= Spec\ k(s) \hookrightarrow S$ un point ferm\'e de $S$ et $f_{s}:X_{s} \rightarrow s$ la fibre de $f$ en $s$. On note $\mathcal{V}(s)= W(k(s))\otimes_{W} \mathcal{V}$, o\`u $W=W(k),\ W(k(s)) $ sont les anneaux de vecteurs de Witt \`a coefficients dans $k$ et $k(s)$ respectivement, et $K(s)$ le corps des fractions de $\mathcal{V}(s)$. Le morphisme $i_{s}$ d\'efinit des foncteurs images inverses [B 3, (2.3.2)(iv) et (2.3.6)], [LS, chap 7]
$$
\xymatrix{
\  i_{s}^{\dag \ast}: \ Isoc^{\dag}((S,T)/ \mathcal{V})\ \ar[rr]&&\ Isoc^{\dag} (Spec\ k(s) / K(s))\ar [d]_{\simeq}\\
{}&&\ Isoc (Spec\ k(s) / K(s)),
}
\leqno (3.4.6.1)
$$
qui, pour $\rho$ propre, devient
$$
\xymatrix{
\  i_{s}^{\dag \ast}: \ Isoc^{\dag}(S /  K)\ \ar[rr]&&\ Isoc^{\dag} (Spec\ k(s) / K(s))\ar [d]_{\simeq}\\
{}&&\ Isoc (Spec\ k(s) / K(s)),
}
\leqno (3.4.6.2)
$$
et qui, pour $S=T$, devient
$$
\xymatrix{
\  \hat{i}_{s}^{\ast}: \ Isoc (S /  K)\ \longrightarrow Isoc (Spec\ k(s) / K(s)).
}
\leqno (3.4.6.3)
$$
Nous allons donner maintenant une r\'ealisation explicite de ces foncteurs: dans la notation $i_{s}^{\dag \ast}$, le $()^{\dag}$ n'est pas utile, sauf pour insister que l'on travaille avec la cat\'egorie surconvergente, et pas seulement la convergente.\\

Dans le carr\'e cart\'esien de $\mathcal{V}$-sch\'emas formels plats et s\'epar\'es 
$$
\xymatrix{
\mathcal{U}(s) \ar[rr]^{u'} \ar[d]_{v'} && \mathcal{S} \ar[d]^{v}\\
Spf\ \mathcal{V}(s)\ar[rr]_{u}&&Spf\ \mathcal{V}\ 
}
$$
les morphismes $u$ et $u'$ (resp $v$ et $v'$) sont finis \'etales (resp lisses et s\'epar\'es). Par lissit\'e de $\mathcal{S}$ sur $\mathcal{V}$ le morphisme compos\'e 

$$Spec\ k(s)\ \overset{i_{s}}{\hooklongrightarrow} S\ \hooklongrightarrow \mathcal{S}$$

\noindent se rel\`eve en un morphisme s\'epar\'e 
$$g_{s}\ :\ Spf\ \mathcal{V}(s) \longrightarrow \ \mathcal{S}\ .$$
 Par la propri\'et\'e universelle du produit fibr\'e, $g_{s}$ se factorise en 
$$g_{s}\ :\ Spf\ \mathcal{V}(s)\ \overset{\tilde{g}_{s}}{\longrightarrow} \mathcal{U}(s)\ \overset{u'}{\hooklongrightarrow} \mathcal{S}$$
avec $v'\circ \tilde{g}_{s} = Id_{Spf\ \mathcal{V}(s)}$: on en d\'eduit que $\tilde{g}_{s}$ est une immersion  ferm\'ee. En notant $\overline{g}_{s}$ le morphisme compos\'e
$$\overline{g}_{s}\ :\ Spf\ \mathcal{V}(s)\ \overset{g_{s}}{\longrightarrow} \mathcal{S}\ \hooklongrightarrow \mathcal{T},$$
 les foncteurs $ i_{s}^{\dag \ast}$ et $\hat{i}_{s}^{\ast}$ sont induits par $\overline{g}_{s}^{\ast}$ [B 3, (2.3.6) et (2.3.2)(iv)], [LS, chap 7], ce qui fournit une r\'ealisation explicite de ces foncteurs comme nous l'annoncions.\\

\noindent\textbf{Th\'eor\`eme (3.4.7).}
\textit{ On se place sous les hypoth\`eses et notations de (3.4.6).
\begin{enumerate}
\item [(3.4.7.1)] Soit $E \in Isoc^{\dag}((X,Y) / \mathcal{V})$. Alors, pour tout entier $i\geqslant 0$ et tout point ferm\'e $s$ de $S$, on a des isomorphismes
$$
\begin{array}{rcl}
i_{s}^{\dag \ast}R^{i}\overline{f}_{rig\ast}((X,Y)/\mathcal{T};\ E)& \simeq& R^{i}f_{s\ rig\ast}(X_{s}/\mathcal{V}(s); E_{X_{s}})\\
&\simeq & R^{i}f_{s\ conv\ast}(X_{s}/\mathcal{V}(s);\widehat{E_{X_{s}}})\\
&\simeq &H^{i}_{rig}( X_{s}/ K(s); E_{X_{s}})\\
&= & H^{i}_{conv}( X_{s}/ K(s); \widehat{E_{X_{s}}}) 
\end{array}
$$
o\`u l'isocristal convergent $\widehat{E_{X_{s}}}\in Isoc\ (X_{s}/K(s))= Isoc^{\dag}\ (X_{s}/K(s))$ co\"\i ncide avec l'isocristal surconvergent $E_{X_{s}}\in Isoc^{\dag}\ (X_{s}/K(s)).$
\item[(3.4.7.2)] Supposons $\rho\:\ \mathcal{T}\rightarrow\ Spf\ \mathcal{V}$ propre et soit $E \in Isoc^{\dag}\ (X/ K).$ Alors, pour tout entier $i\geqslant 0$ et tout point ferm\'e $s$ de $S$, on a des isomorphismes
$$
\begin{array}{rcl}
i_{s}^{\dag \ast}R^{i}f_{rig\ast}(X/\mathcal{T};\ E)& \simeq& R^{i}f_{s\ rig\ast}(X_{s}/\mathcal{V}(s);E_{X_{s}})\\
&\simeq & R^{i}f_{s\ conv\ast}(X_{s}/\mathcal{V}(s); \widehat{E_{X_{s}}})\\
&\simeq &H^{i}_{rig}( X_{s}/ K(s); E_{X_{s}})\\
&= & H^{i}_{conv}( X_{s}/ K(s); \widehat{E_{X_{s}}}) 
\end{array}
$$
o\`u l'isocristal convergent $\widehat{E_{X_{s}}}\in Isoc\ (X_{s}/K(s))= Isoc^{\dag}\ (X_{s}/K(s))$ co\"\i ncide avec l'isocristal surconvergent $E_{X_{s}}\in Isoc^{\dag}\ (X_{s}/K(s)).$
\item[(3.4.7.3)] Supposons que $S=T$ et soit $\mathcal{E}\in Isoc\ (X/K)$.  Alors, pour tout entier $i\geqslant 0$ et tout point ferm\'e $s$ de $S$, on a des isomorphismes
$$
\begin{array}{rcl}
\hat{i}_{s}^{ \ast}R^{i}f_{conv\ast}(X/\mathcal{T};\ \mathcal{E})& \simeq& R^{i}f_{s\ conv\ast}(X_{s}/\mathcal{V}(s); \mathcal{E}_{X_{s}})\\
&= & H^{i}_{conv}( X_{s}/ K(s); \mathcal{E}_{X_{s}}). 
\end{array}
$$
\end{enumerate}
}
\noindent\textit{D\'emonstration.} On applique le changement de base par ($i_{s},\ \overline{g}_{s}$) et les corollaires (3.4.5.3) et (3.4.5.5). $\square$\\

\noindent\textbf{ 3.4.8.}
Soient $S$ un $k$-sch\'ema lisse et s\'epar\'e et $f\ :  X\rightarrow S$ un $k$-morphisme projectif et lisse. Notons $S= \underset{\alpha}{\bigcup} \ S_{\alpha,0}$ une d\'ecomposition de $S$ en r\'eunion d'ouverts connexes affines $S_{\alpha,0}= Spec(A_{\alpha,0}), \ A_{\alpha}= \mathcal{V}[t_{1},...,t_{d_{\alpha}}]/J_{\alpha}$ une $\mathcal{V}$-alg\`ebre lisse relevant $A_{\alpha,0}$ dont on a fix\'e une pr\'esentation et $S_{\alpha}= Spec (A_{\alpha})$. On d\'esigne par $\overline{S}_{\alpha}$ l'adh\'erence sch\'ematique de $S_{\alpha}$ dans $\mathbb{P}_{\mathcal{V}}^{d_{\alpha}}$, par $\mathcal{S}_{\alpha}$ (resp $\overline{\mathcal{S}}_{\alpha}$) le compl\'et\'e formel $\mathfrak{m}$-adique de $S_{\alpha}$ (resp $\overline{S}_{\alpha}$) et par 
$$
f_{\alpha}:X_{\alpha,0}=X\times_{S}S_{\alpha,0}\longrightarrow S_{\alpha,0}
$$
la restriction de $f$. Quitte \`a d\'ecomposer $X_{\alpha,0}$ en somme disjointe de ses composantes connexes on peut supposer $X_{\alpha,0}$ connexe. D'apr\`es [Et 5, (3.2.1.1)] il existe un carr\'e cart\'esien
$$
\begin{array}{c}
\xymatrix{
X_{\alpha} \ar@{^{(}->} [rr] \ar [d]_{h_{\alpha}} && \overline{X}_{\alpha}\ar[d]^{\overline{h}_{\alpha}}\\
S_{\alpha} \ar@{^{(}->} [rr] && \overline{S}_{\alpha} 
}
\end{array}
$$
dans lequel $\overline{h}_{\alpha}$ est projectif, $h_{\alpha}$ est un rel\`evement projectif de $f_{\alpha}$ et les fl\`eches horizontales sont des immersions ouvertes. Le compl\'et\'e formel de ce carr\'e est d'apr\`es [Et 5, th\'eo (3.2.1)] un carr\'e cart\'esien de $\mathcal{V}$-sch\'emas formels 
$$
\begin{array}{c}
\xymatrix{
\mathcal{X}_{\alpha} \ar@{^{(}->} [rr] \ar [d]_{\hat{h}_{\alpha}} && \overline{\mathcal{X}}_{\alpha}\ar[d]^{\hat{\overline{h}}_{\alpha}}\\
\mathcal{S}_{\alpha} \ar@{^{(}->} [rr] && \overline{\mathcal{S}}_{\alpha} 
}
\end{array}
\leqno (3.4.8.1)
$$
dans lequel $\hat{\overline{h}}_{\alpha}$ est projectif, $\hat{h}_{\alpha}$ est un rel\`evement projectif de la restriction $f_{\alpha}$ de $f$ et les fl\`eches horizontales sont des immersions ouvertes.\\

\noindent\textbf{Th\'eor\`eme (3.4.8.2).}
\textit{Sous les hypoth\`eses (3.4.8) supposons que, pour tout $\alpha$, $X_{\alpha}$ est plat sur $\mathcal{V}$. Alors, pour tout entier $i \geqslant 0$, on a un diagramme commutatif de foncteurs naturels induits par $f$ et d\'efinis ci-apr\`es en (3.4.8.5)
$$
\xymatrix{
 Isoc^{\dag}(X/ K) \ar[rr]^{R^{i}f_{rig\ast}}\ar[d]&&Isoc^{\dag}(S/K)\ar [d]\\
Isoc (X/ K)\ar[rr]^{R^{i}f_{conv\ast}}&&Isoc(S/ K)
}
$$
o\`u les fl\`eches verticales sont les foncteurs d'oubli.\\
De plus ces foncteurs commutent \`a tout changement de base $S'\rightarrow S$ entre $k$-sch\'emas lisses et s\'epar\'es; en particulier ils commutent aux passages aux fibres en les points ferm\'es de $S$.
}\\

\noindent\textit{D\'emonstration.} Puisque, pour tout $\alpha$, $X_{\alpha}$ est plat sur $\mathcal{V}$, le th\'eor\`eme [Et 5, th\'eo (3.2.1)] prouve que $\hat{h}_{\alpha}$ est un rel\`evement projectif et lisse de la restriction $f_{\alpha}$ de $f$. Soient $E \in Isoc^{\dag}(X/K)$ et $E_{\alpha}$ sa restriction \`a $X_{\alpha}$: gr\^ace \`a l'existence du carr\'e cart\'esien (3.4.8.1) dans lequel $\hat{h}_{\alpha}$ est un rel\`evement projectif et lisse de  $f_{\alpha}$ on conclut \`a l'aide du th\'eor\`eme (3.4.4) que
$$
R^{i}f_{\alpha rig\ast}(X_{\alpha,0}/\overline{\mathcal{S}}_{\alpha}; E_{\alpha}) \in Isoc^{\dag}(S_{\alpha,0}/K).
$$
Montrons que celui-ci ne d\'epend que de $S_{\alpha,0}$ et non de $\overline{\mathcal{S}}_{\alpha}$.\\

Supposons donn\'es un $k$-sch\'ema propre $\overline{S}'_{\alpha,0}$, un $\mathcal{V}$-sch\'ema formel propre $\overline{\mathcal{S}}'_{\alpha}$, une immersion ouverte dominante $j'_{\alpha,0}: S_{\alpha,0}\hookrightarrow \overline{S}'_{\alpha,0} $ et une immersion ferm\'ee $i'_{\alpha,0}: \overline{S}'_{\alpha,0}\hookrightarrow \overline{\mathcal{S}}'_{\alpha,0} $ tels que le morphisme $\overline{\mathcal{S}}'_{\alpha}\longrightarrow Spf \mathcal{V}$ soit lisse sur un voisinage de $S_{\alpha,0}$ dans $ \overline{\mathcal{S}}'_{\alpha}$. Notons $\overline{T}_{\alpha,0}$ l'adh\'erence sch\'ematique de $S_{\alpha,0}$ plong\'e diagonalement dans $\overline{S}''_{\alpha}=\overline{S}_{\alpha,0}\times_{k} \overline{S}'_{\alpha,0}$ et $\overline{\mathcal{S}}''_{\alpha}=\overline{\mathcal{S}}_{\alpha}\hat{\times}_{\mathcal{V}} \overline{\mathcal{S}}'_{\alpha}$. On a un diagramme commutatif
$$
\begin{array}{c}
\xymatrix{
S_{\alpha,0}\ \ar@{^{(}->}[rr]^{j_{T}} \ar[d]_{id}&& \overline{T}_{\alpha,0}\ \ar@{^{(}->}[rr]^{i_{T}} \ar[d]^{p_{\alpha}} &&\overline{\mathcal{S}}''_{\alpha} \ar[d]^{\overline{p}_{1,\alpha}}\\
S_{\alpha,0}\ \ar@{^{(}->}[rr]_{j_{S}}&&\overline{S}_{\alpha,0}\ \ar@{^{(}->}[rr]_{i_{S}} &&\overline{\mathcal{S}}
}
\end{array}
\leqno{(3.4.8.3)}
$$
dans lequel les fl\`eches verticales sont induites par la premi\`ere projection, les $i$ (resp les $j$) sont des immersions ferm\'ees (resp ouvertes). Comme $p_{\alpha}$ est propre, le foncteur $(p_{\alpha},\overline{p}_{1,\alpha})^{\ast}$ est une \'equivalence de cat\'egorie [B 3, (2.3.5)], [LS, 7.1] de la cat\'egorie des isocristaux sur $S_{\alpha,0}/K$ surconvergents le long de $Z_{\alpha,0}=\overline{S}_{\alpha,0}\setminus S_{\alpha,0}$ dans la cat\'egorie des isocristaux sur $S_{\alpha,0}/K$ surconvergents le long de $Z''_{\alpha,0}=\overline{T}_{\alpha,0}\setminus S_{\alpha,0}$. De plus $\overline{p}_{1,\alpha}$ est propre, et lisse sur un voisinage de $S_{\alpha,0}$ dans $\overline{\mathcal{S}}''_{\alpha}$, donc par le th\'eor\`eme (3.4.4) on a isomorphisme de changement de base
$$\overline{p}_{1\alpha}^{\ \ast}R^{i}f_{\alpha rig \ast}(X_{\alpha,0}/\overline{\mathcal{S}}_{\alpha};E_{\alpha}) \overset{\sim}	{\longrightarrow} 	R^{i}{f}_{\alpha rig \ast}(X_{\alpha,0}/\overline{\mathcal{S}}''_{\alpha};E_{\alpha}),$$
et de m\^eme
$$\overline{p}_{2\alpha}^{\ \ast}R^{i}f_{\alpha rig \ast}(X_{\alpha,0}/\overline{\mathcal{S}}'_{\alpha};E_{\alpha}) \overset{\sim}	{\longrightarrow} 	R^{i}{f}_{\alpha rig \ast}(X_{\alpha,0}/\overline{\mathcal{S}}''_{\alpha};E_{\alpha}),$$
d'o\`u un isomorphisme canonique
$$R^{i}f_{\alpha rig \ast}(X_{\alpha,0}/\overline{\mathcal{S}}_{\alpha};E_{\alpha}) \overset{\sim}	{\longrightarrow} 	R^{i}{f}_{\alpha rig \ast}(X_{\alpha,0}/\overline{\mathcal{S}}'_{\alpha};E_{\alpha}).
\leqno (3.4.8.4)
$$
Ainsi $R^{i}f_{\alpha rig \ast}(X_{\alpha,0}/\overline{\mathcal{S}}_{\alpha};E_{\alpha})$ ne d\'epend que de $S_{\alpha,0}$; de plus sur $S_{\alpha,0}\cap S_{\beta,0}$ ces isocristaux se recollent car on a des isomorphismes analogues \`a (3.4.8.4) et v\'erifiant la condition de cocycles. Ces donn\'ees qui se recollent fournissent un isocristal surconvergent sur $S/K$ [B 3, (2.3.2)], [LS, 8.1] not\'e
$$
R^{i}f_{rig \ast}(E),
\leqno (3.4.8.5)
$$
et c'est le seul possible d'apr\`es le th\'eor\`eme de pleine fid\'elit\'e de Kedlaya  [Ked 3, Theo 5.2.1] qui prouve que l'extension de l'isocristal surconvergent de $S_{\alpha,0}\cap S_{\beta,0}$  \`a $S_{\alpha,0}$ est unique.\\
On raisonne de mani\`ere analogue pour $R^{i}f_{conv \ast}(E)$ et le carr\'e du th\'eor\`eme est clairement commutatif.\\
L'assertion sur les changements de base $S'\rightarrow S$ r\'esulte du th\'eor\`eme (3.4.4). $\square$\\

\noindent\textbf{Corollaire (3.4.8.6).}
\textit{Sous les hypoth\`eses (3.4.8), supposons que $f$ d\'efinit $X$ comme une intersection compl\`ete relativement \`a $S$ dans un espace projectif sur $S$ [Et 5, d\'ef (3.2.5)]. Alors les conclusions du th\'eor\`eme (3.4.8.2) demeurent valides.} \\
 
\noindent \textit{D\'emonstration}. Avec les notations de [Et 5, (3.2.5.2)], chaque $f_{\alpha, \beta} :\  X_{\alpha, \beta} \rightarrow S_{\alpha, \beta}$ se rel\`eve d'apr\`es [Et 5, cor (3.2.6)] en un morphisme projectif et lisse au-dessus de $\mathcal{V}$ et donne lieu \`a un diagramme commutatif tel que (3.4.8.1) en rempla\c cant $\mathcal{S}_{\alpha}$ par $\mathcal{S}_{\alpha, \beta}$, et on conclut comme pour (3.4.8.2). $\square$\\

On a la variante suivante du th\'eor\`eme (3.4.8.2):\\

\noindent\textbf{Th\'eor\`eme (3.4.9).}
\textit{Supposons que le th\'eor\`eme (3.4.8.2) est vrai pour tout morphisme projectif et lisse $f :  X\rightarrow S$ avec $S$ un $k$-sch\'ema affine et lisse (donc sans faire l'hypoth\`ese $X_{\alpha}$ plat sur $\mathcal{V}$). Soient $S$ un $k$-sch\'ema lisse et s\'epar\'e et $f :  X\rightarrow S$ un $k$-morphisme projectif et lisse. Alors, pour tout entier $i \geqslant 0$, on a un diagramme commutatif de foncteurs naturels induits par $f$ et d\'efinis \`a la mani\`ere de(3.4.8.5)
$$
\xymatrix{
 Isoc^{\dag}(X/ K) \ar[rr]^{R^{i}f_{rig\ast}}\ar[d]&&Isoc^{\dag}(S/K)\ar [d]\\
Isoc (X/ K)\ar[rr]^{R^{i}f_{conv\ast}}&&Isoc(S/ K)
}
$$
o\`u les fl\`eches verticales sont les foncteurs d'oubli.\\
De plus ces foncteurs commutent \`a tout changement de base $S'\rightarrow S$ entre $k$-sch\'emas lisses et s\'epar\'es; en particulier ils commutent aux passages aux fibres en les points ferm\'es de $S$.
}\\

\noindent\textit{D\'emonstration.} Il suffit de d\'ecomposer $S$ en sch\'emas affines et lisses sur $k$ et d'op\'erer les recollements comme dans la preuve de (3.4.8.2).
 $\square$\\

 \end{document}